\documentclass[12pt,leqno]{article}
\usepackage{amsmath,amssymb,graphicx,MnSymbol,mathtools,url}

\DeclareMathOperator{\graph}{graph}
\DeclareMathOperator{\projection}{proj}
\DeclareMathOperator{\spt}{spt}
\DeclareMathOperator{\divergence}{div}
\DeclareMathOperator{\reg}{reg}
\DeclareMathOperator{\sing}{sing}
\DeclareMathOperator{\dist}{dist}
\DeclareMathOperator{\clos}{clos}

\def\res{\hbox{ {\vrule height .22cm}{\leaders\hrule\hskip.2cm} }}

\newcommand{\C}{\mathbb{C}}

\newcommand{\D}{\overline{D}}

\newcommand{\x}{\mathbf{x}}
\newcommand{\z}{\mathbf{z}}
\newcommand{\E}{\mathbf{E}}
\newcommand{\I}{\mathbf{I}}
\newcommand{\N}{\mathbf{N}}
\newcommand{\R}{\mathbf{R}}
\newcommand{\TI}{\mathbf{TI}}

\newcommand{\HH}{\mathcal{H}}

\newcommand{\proj}[2]{\projection_{#1} #2}
\newcommand{\dive}[2]{\divergence_{#1} #2}
\newcommand{\gph}[2]{\graph_{#1} #2}

\numberwithin{equation}{section}
\newtheorem{theorem}[equation]{Theorem}
\newtheorem{lemma}[equation]{Lemma}
\newtheorem{definition}[equation]{Definition}

\begin{document}
\begin{flushleft}

TITLE: Co-Dimension One Area-Minimizing Currents with $C^{1,\alpha}$ Tangentially Immersed Boundary. 

\medskip

AUTHOR: Leobardo Rosales, Keimyung University

\medskip

ABSTRACT: We introduce and study co-dimension one area-minimizing locally rectifiable currents $T$ with $C^{1,\alpha}$ tangentially immersed boundary: $\partial T$ is locally a finite sum of orientable co-dimension two submanifolds which only intersect tangentially with equal orientation. We show that any such $T$ is supported in a smooth hypersurface near any point on the support of $\partial T$ where $T$ has tangent cone which is a hyperplane with constant orientation but non-constant multiplicity. We also introduce and study co-dimensional one area-minimizing locally rectifiable currents $T$ with boundary having co-oriented mean curvature: $\partial T$ has generalized mean curvature $H_{\partial T} = h \nu_{T}$ with $h$ a real-valued function and $\nu_{T}$ the generalized outward pointing unit normal of $\partial T$ with respect to $T.$

\medskip


\medskip


\section{Introduction}

Our goal in this work is to generalize the boundary regularity theory for co-dimension one area-minimizing locally rectifiable currents established in \cite{HS79} and \cite{W83}. To this end, we introduce two definitions. Suppose $T$ is an $n$-dimensional area-minimizing locally rectifiable current in $\R^{n+1}.$ First, we say $T$ has $C^{1,\alpha}$ tangentially immersed boundary with $\alpha \in (0,1]$ if $\partial T$ is locally a finite sum of orientable $(n-1)$-dimensional (embedded) submanifolds which meet only tangentially with equal orientation; see Definition \ref{immersedboundary}. Second, we say the boundary of $T$ has co-oriented mean curvature if $\partial T$ has generalized mean curvature $H_{\partial T} = h \nu_{T}$ for $h$ a real-valued function and $\nu_{T}$ the generalized outward pointing unit normal of $\partial T$ with respect to $T$; see Definition \ref{cmcboundary}, in addition to Lemma 3.1 of \cite{B77} and (2.9) of \cite{E89} for the existence of $\nu_{T}.$

\bigskip

In order to proceed, we briefly describe the now classical results of \cite{HS79} and \cite{W83}. Suppose $T$ is an $n$-dimensional area-minimizing locally rectifiable current in $\R^{n+1},$ and suppose $x$ is in the support of $\partial T.$ First, \cite{HS79} shows that if $\partial T$ near $x$ corresponds to integrating over an orientable $(n-1)$-dimensional $C^{1,\alpha}$ submanifold with multiplicity one and $\alpha \in (0,1),$ then one of the following two occurs: $T$ corresponds to integrating (with multiplicity one) over a $C^{1,\alpha}$ hypersurface-with-boundary, with boundary the support of $\partial T$; the support of $T$ near $x$ is an analytic minimal hypersurface $M$ containing the support of $\partial T$ near $x,$ and $T$ corresponds to integrating over $M$ with multiplicities $\theta,(\theta-1)$ for some positive integer $\theta$ respectively over the two $C^{1,\alpha}$ regions of $M$ determined by the support of $\partial T.$ This result is extended by \cite{W83}. If instead $\partial T$ near $x$ corresponds to integrating over an orientable $(n-1)$-dimensional $C^{1,\alpha}$ submanifold with multiplicity a positive integer $m$ and $\alpha \in (0,1),$ then one of the following two occurs:  near $x$ we have that $T$ corresponds to integrating over a disjoint union of $C^{1,\alpha}$ hypersurfaces-with-boundary each with some multiplicity, each hypersurface having boundary $\partial T$; the support of $T$ near $x$ is an analytic minimal hypersurface $M$ containing the support of $\partial T$ near $x,$ and $T$ corresponds to integrating over $M$ with multiplicities $(m+\theta),\theta$ for some positive integer $\theta$ respectively over the two $C^{1,\alpha}$ regions of $M$ determined by the support of $\partial T.$

\bigskip

We now describe our main results, the first of which is the most important:

\bigskip

{\bf Theorem \ref{hyperplanetangentconegraph}:} \emph{Suppose $T$ is an $n$-dimensional area-minimizing locally rectifiable current in $\R^{n+1}$ with $C^{1,\alpha}$ tangentially immersed boundary where $\alpha \in (0,1].$ Suppose $x$ is in the support of $\partial T,$ and that $T$ at $x$ has a tangent cone $\C$ which is a \emph{hyperplane with constant orientation but non-constant multiplicity}. Then the support of $T$ near $x$ is the graph of a smooth solution to the minimal surface equation $u$ off the hyperplane supporting the (now unique) tangent cone $\C,$ and the orientation vector of $T$ near $x$ corresponds to the upward pointing unit normal of the graph of $u.$}

\bigskip

We say that an $n$-dimensional current $\C$ in $\R^{n+1}$ is a \emph{hyperplane with constant orientation but non-constant multiplicity} if (after rotation) $\C$ corresponds to integrating over $\{ (y_{1},\ldots,y_{n}) \in \R^{n}:  y_{n}>0 \}$ with orientation $e_{n+1}$ and multiplicity $(m+\theta)$ and over $\{ (y_{1},\ldots,y_{n}) \in \R^{n}:  y_{n}<0 \}$ with orientation $e_{n+1}$ and multiplicity $\theta$ for positive integers $m,\theta$; see Definition \ref{conesdefinition}. 

\bigskip

If $\partial T$ is $C^{1,1}$ tangentially immersed with Lipschitz co-oriented mean curvature, we conclude $\partial T$ is regular:

\bigskip

{\bf Theorem \ref{hyperplanetangentconeregularity}:} \emph{Suppose $T$ is an $n$-dimensional area-minimizing locally rectifiable current in $\R^{n+1}$ with $C^{1,1}$ tangentially immersed boundary and where $\partial T$ has co-oriented mean curvature $H_{\partial T} = h \nu_{T}$ with $h$ a Lipschitz real-valued function and $\nu_{T}$ the generalized outward pointing unit normal of $\partial T$ with respect to $T.$ If $x$ is in the support of $\partial T$ and $T$ at $x$ has a tangent cone which is a hyperplane with constant orientation but non-constant multiplicity, then $\partial T$ near $x$ corresponds to integrating over an $(n-1)$-dimensional submanifold, which is $C^{2,\alpha}$ for any $\alpha \in (0,1),$ with multiplicity.}

\bigskip

In this case, \cite{HS79} and \cite{W83} imply that $T$ near $x$ is supported in an analytic hypersurface $M$ containing the support of $\partial T$ near $x,$  and $T$ corresponds to integrating over $M$ with multiplicities $(m+\theta),\theta$ for some positive integers $m,\theta$ respectively over the two $C^{2,\alpha}$ regions of $M$ determined by the support of $\partial T.$ 

\bigskip

Finally, we have the following very geometric partial regularity theorem: 

\bigskip

{\bf Theorem \ref{main}:} \emph{Suppose $T$ is an $n$-dimensional area-minimizing locally rectifiable current in $\R^{n+1}$ with $C^{1,1}$ tangentially immersed boundary and where $\partial T$ has co-oriented mean curvature $H_{\partial T} = h \nu_{T}$ with $h$ a Lipschitz real-valued function and $\nu_{T}$ the generalized outward pointing unit normal of $\partial T$ with respect to $T.$ Suppose $x$ is in the support of $\partial T,$ and that near $x$ the support of $T$ equals a finite union of $C^{1}$ hypersurfaces-with-boundary. Then the support of $T$ near $x$ is the finite union of $C^{1,1}$ hypersurfaces-with-boundary which pairwise meet only at common boundary points.}

\bigskip

The better way to understand Theorem \ref{main} is through the contrapositive: at any point $x$ in the support of $\partial T$ near which the support of $T$ does \emph{not} equal a finite union of $C^{1,1}$ hypersurfaces-with-boundary which pairwise meet only at common boundary points, the support of $T$ near $x$ must have extremely complicated structure; perhaps for example infinite topology at $x.$ 

\bigskip

Before we discuss our main results with more detail, we describe a partial boundary regularity result given by the author in \cite{R15d}. For convenience, we state this result here as Theorem \ref{partialboundaryregularity}, as it will be crucial to the proofs. We describe it loosely now:

\bigskip

{\bf Theorem \ref{partialboundaryregularity} (Theorem 2.1 of \cite{R15d}):} \emph{Suppose $T$ is an $n$-dimensional area-minimizing locally rectifiable current in $\R^{n+1}$ with $C^{1,\alpha}$ tangentially immersed boundary where $\alpha \in (0,1].$ Suppose $x$ is in the support of $\partial T,$ and that $T$ at $x$ has a tangent cone $\C$ which is a hyperplane with constant orientation but non-constant multiplicity. Then near $x$ a large portion of the support of $T$ can be written as the graph of a $C^{1,\frac{\alpha}{4n+6}}$ function defined over a large region of the hyperplane supporting $\C$; this region is large enough to conclude that $\C$ is the unique tangent cone of $T$ at $x.$}

\bigskip

We now discuss Theorems \ref{hyperplanetangentconegraph},\ref{hyperplanetangentconeregularity},\ref{main} in more detail.

\subsection{Tangentially Immersed Boundaries}

Throughout this section, let $T$ be a $n$-dimensional area-minimizing locally rectifiable current over $\R^{n+1}$ with $C^{1,\alpha}$ tangentially immersed boundary, where $\alpha \in (0,1].$ We also suppose that $0$ is in the support of $\partial T,$ and that $T$ at $0$ has a tangent cone $\C$ which is integrating over $\{ (y_{1},\ldots,y_{n}) \in \R^{n}:  y_{n}>0 \}$ with orientation $e_{n+1}$ and multiplicity $(m+\theta)$ and over $\{ (y_{1},\ldots,y_{n}) \in \R^{n}:  y_{n}<0 \}$ with orientation $e_{n+1}$ and multiplicity $\theta$ for positive integers $m,\theta.$ 

\bigskip

We now discuss Theorem \ref{hyperplanetangentconegraph}, which concludes that the support of $T$ near $0$ is the graph of $u$ a smooth solution to the minimal surface equation, and the orientation vector $\ast \vec{T}$ of $T$ near $0$ is the upward pointing unit normal of the graph of $u.$ It will be evident that Theorem \ref{hyperplanetangentconegraph} is a natural consequence of the partial boundary regularity result Theorem \ref{partialboundaryregularity} (Theorem 2.1 of \cite{R15d}). To prove Theorem \ref{hyperplanetangentconegraph}, we must first prove that having unique tangent cone a hyperplane with constant orientation but non-constant multiplicity (see Definition \ref{conesdefinition}) is an open condition along $\partial T$:

\bigskip

{\bf Lemma \ref{hyperplanetangentcone}:} \emph{For all $x$ in the support of $\partial T$ near $0,$ we have that $T$ has at $x$ a unique tangent cone which is a hyperplane with constant orientation but non-constant multiplicity. } 

\bigskip

Having stated Lemma \ref{hyperplanetangentcone}, it is then expected that Theorem \ref{hyperplanetangentconegraph} is true. Lemma \ref{hyperplanetangentcone} is the linchpin of the present work, and so we discuss it in detail. The proof of Lemma \ref{hyperplanetangentcone} is a generalization of the proof of Theorem 8 of \cite{R16}, where special cases of the present results appear in the context of two-dimensional solutions to the $c$-Plateau problem in space; see \S 1.4. In turn, the techniques used in \cite{R16} (and presently) are similar to those used by the author in \cite{R10},\cite{R11},\cite{R15b} to study the two-valued minimal surface equation, a degenerate second-order PDE first introduced in \cite{SW07} to produce examples of co-dimension one $C^{1,\alpha}$ stable branched minimal immersions.

\bigskip

The proof of Lemma \ref{hyperplanetangentcone} follows by analyzing the cross-sections of $T$ across affine planes $\{z\} \times \R^{2}$ for $z \in \R^{n-1}$ near the origin. For simplicity, we describe the proof in the case $n=2.$ First, Theorem \ref{partialboundaryregularity} implies (after rescaling) that the support of $T$ in the bored-out unit ball
$$\{ x=(x_{1},x_{2},x_{3}) \in \R^{3}: \delta < |(x_{2},x_{3})|, \ |x| < 1 \}$$
is the graph of a function defined off the region in the horizontal plane
$$\{ y=(y_{1},y_{2}) \in \R^{2}: \delta < |y_{2}|, \ |y| < 1 \},$$
for some small $\delta>0.$ We now suppose for contradiction $\tilde{x}=(\tilde{x}_{1},\tilde{x}_{2},\tilde{x}_{3})$ is a point in the support of $\partial T$ with $|\tilde{x}| < \delta$ and so that $T$ at $\tilde{x}$ does \emph{not} have unique tangent cone which is a plane with constant orientation but non-constant multiplicity. It follows (see Theorem \ref{monotonicity}) that every tangent cone $\C$ of $T$ at $\tilde{x}$ must be a \emph{sum of half-planes with constant orientation after rotation} meeting along a common line; see Definition \ref{conesdefinition}. As $T$ has $C^{1,\alpha}$ tangentially immersed boundary, the \emph{spine} of any such $\C$ is a fixed line, namely the tangent line of $\partial T$ at $x,$ which is close to the $x_{1}$-axis depending on how close $\tilde{x}$ is to the origin.

\bigskip

We now consider the cross-section of the support of $T$ along the affine plane $\{ \tilde{x}_{1}\} \times \R^{2}.$ In the region $\{  x=(\tilde{x}_{1},x_{2},x_{3}) \in \R^{3}: \delta < |(x_{2},x_{3})|, \ |x| < 1 \}$ this cross-section is by Theorem \ref{partialboundaryregularity} two curves $\Gamma_{1},\Gamma_{2}$ which are respectively $C^{1}$ close to the line segments 
$$\{ (\tilde{x}_{1},t,0) \in \R^{3}: t \in (\delta,1-|\tilde{x}_{1}|) \} \text{ and } \{ (\tilde{x}_{1},t,0) \in \R^{3}: t \in (-\delta,-1+|\tilde{x}_{1}|) \}.$$
In particular, the orientation vector $\ast \vec{T}$ of $T$ (the generalized unit normal vector field of $T$) satisfies $\ast \vec{T} \approx e_{3}$ along $\Gamma_{1},\Gamma_{2}.$ 

\bigskip

On the other hand, assuming (without loss of generality) that $\tilde{x}=(\tilde{x}_{1},\tilde{x}_{2},\tilde{x}_{3})$ is the only point in the support of $\partial T$ on the cross-section $\{\tilde{x}_{1}\} \times \R^{2},$ and that $\tilde{x}_{1}$ is a regular value of the function $f(x_{1},x_{2},x_{3})=x_{1}$ over the support of $T$ away from the support of $\partial T$ (by Sard's theorem and interior regularity of co-dimension one area-minimizing locally rectifiable currents), then we essentially argue that $\Gamma_{1},\Gamma_{2}$ are contained respectively in smooth curves $\tilde{\Gamma}_{1},\tilde{\Gamma}_{2}$ meeting only at $\tilde{x}.$ Moreover, $\ast \vec{T}$ is smooth over $\tilde{\Gamma}_{1},\tilde{\Gamma}_{2}.$ However, as every tangent cone of $T$ at $\tilde{x}$ is a sum of half-planes with constant orientation after rotation (see Definition \ref{conesdefinition}) with spine close to the $x_{1}$-axis (as noted above), then considering the continuity of $\ast \vec{T}$ along $\tilde{\Gamma}_{1}$ and $\tilde{\Gamma}_{2}$ yields a contradiction; essentially, we contradict that $\ast \vec{T}$ is continuous over $\tilde{\Gamma}_{2},$ due to orientation.

\bigskip

Having done most of the work in proving Lemma \ref{hyperplanetangentcone}, then the proof of Theorem \ref{hyperplanetangentconegraph} is a relatively short proof by induction based on the number of $(n-2)$-dimensional $C^{1,\alpha}$ submanifolds in the decomposition of $\partial T$ near $0.$

\subsection{Tangentially Immersed Boundaries with Co-Oriented Mean Curvature}

We describe in this subsection the proofs of Theorems \ref{hyperplanetangentconeregularity},\ref{main}. Throughout this subsection suppose $T$ is an $n$-dimensional area-minimizing locally rectifiable current over $\R^{n+1}$ with $C^{1,1}$ tangentially immersed boundary and so that $\partial T$ has co-oriented mean curvature $H_{\partial T} = h \nu_{T}$ where $h$ is a Lipschitz real-valued function and $\nu_{T}$ is the generalized outward pointing unit normal of $\partial T$ with respect to $T.$ Suppose as well that $0$ is in the support of $\partial T.$

\bigskip

We begin by discussing Theorem \ref{hyperplanetangentconeregularity}, which concludes that if $T$ at $0$ has tangent cone which is a hyperplane with constant orientation but non-constant multiplicity, then $\partial T$ is regular near $0.$ We prove Theorem \ref{hyperplanetangentconeregularity} using the Hopf boundary point lemma at \emph{half-regular} singular points of $\partial T.$ We roughly describe such points in the following lemma: 

\bigskip

{\bf Lemma \ref{halfregular}:} \emph{If $0$ is a singular point of $\partial T,$ then for any $\rho>0$ there exists $x$ in the singular set of $\partial T$ with $|x|< \rho$ and a non-empty open set $U \subset \R^{n+1}$ so that $x \in \partial U$ and the support of $\partial T$ in $U$ is a union of \emph{disjoint} non-empty $(n-1)$-dimensional submanifolds; we call such an $x$ a \emph{half-regular} point.}

\bigskip

The proof of Lemma \ref{halfregular} is by induction on the number of distinct submanifolds in the decomposition of $\partial T$ near $0,$ if $0$ is a singular point. A version of Lemma \ref{halfregular} appears as Lemma 1 of \cite{R16} in the context of two-dimensional solutions to the $c$-Plateau problem in space; see \S 1.4.

\bigskip

The proof of Theorem \ref{hyperplanetangentconeregularity} then follows through naturally. Suppose for contradiction that $T$ at $0$ has unique tangent cone $\C$ which is a hyperplane with constant orientation but non-constant multiplicity, but that $0$ is a singular point of $\partial T.$ Theorem \ref{hyperplanetangentconegraph} implies that $T$ near $0$ is supported in the graph of a function $u$ defined off the hyperplane supporting $\C,$ and the orientation vector of $T$ near $0$ corresponds to the upward pointing unit normal of the graph of $u.$ Since $0$ is a singular point of $\partial T,$ then we can consider a half-regular point $x$ in the singular set of $\partial T$ with $x$ sufficiently close to the origin. Then there is a nonempty open set $U \subset \R^{n+1}$ so that $x \in \partial U$ and the support of $\partial T$ in $U$ is a union of \emph{disjoint} non-empty $(n-1)$-dimensional submanifolds contained in the graph of $u.$ The proof of Theorem \ref{hyperplanetangentconeregularity} then proceeds by applying the Hopf boundary point lemma to $\partial T$ in $U$ at $x,$ using $H_{\partial T} = h \nu_{T}.$

\bigskip

We now discuss Theorem \ref{main}, which concludes that if near $0$ the support of $T$ equals a finite union of $C^{1}$ hypersurfaces-with-boundary, then near $0$ the support of $T$ equals a finite union of $C^{1,1}$ hypersurfaces-with-boundary which pairwise meet only at common boundary points. The proof of Theorem \ref{main} is very geometric, and technical, similar to the proof of Lemma \ref{hyperplanetangentcone}. 

\bigskip

Suppose now that near the $0$ the support of $T$ is a finite union of $C^{1}$ hypersurfaces-with-boundary, and consider the set
$$W = \left\{ \text{ $x$ in the support of $T:$ } 
\begin{aligned}
& \text{$x$ is not in the support of $\partial T$} \\
& \text{or $T$ at $x$ has a tangent cone} \\
& \text{which is a hyperplane with } \\
& \text{constant orientation but} \\
& \text{non-constant multiplicity} 
\end{aligned} \right\}.$$
Interior regularity for co-dimension one area-minimizing locally rectifiable currents, Theorem \ref{hyperplanetangentconeregularity}, as well as the classical boundary regularity given by \cite{HS79} and \cite{W83} imply that $W$ is a smooth embedded hypersurface. Then an analysis similar to the proof of Lemma \ref{hyperplanetangentcone}, using Sard's theorem, shows that $W$ near the origin decomposes into finitely many connected components, each of which is a $C^{1,1}$ hypersurface-with-boundary. To this end, the proof of Theorem \ref{main} closely mirrors the proof of Theorem 9 of \cite{R16}, which again is a version of Theorem \ref{main} in the context of two-dimensional solutions to the $c$-Plateau problem in space; see \S 1.4.

\subsection{Future Work}

We expect to relax the assumption of Theorems \ref{hyperplanetangentconeregularity},\ref{main} that $T$ has $C^{1,1}$ tangentially immersed boundary to $T$ having more generally $C^{1,\alpha}$ tangentially immersed boundary with $\alpha \in (0,1].$ For this, we nevertheless expect it necessary to assume that $\partial T$ has co-oriented mean curvature $H_{\partial T} = h \nu_{T}$ with $h$ a Lipschitz real-valued function. Indeed, we suspect that to generalize Theorems \ref{hyperplanetangentconeregularity},\ref{main}, we must show the following:

\bigskip

{\bf Conjecture:} \emph{Suppose $T$ is an $n$-dimensional area-minimizing locally rectifiable current in $\R^{n+1}$ with $C^{1,\alpha}$ tangentially immersed boundary where $\alpha \in (0,1),$ and that $\partial T$ has co-oriented mean curvature $H_{\partial T} = h \nu_{T}$ with $h$ a Lipschitz real-valued function. Also suppose $x$ is in the support of $\partial T,$ and that $T$ at $x$ has a tangent cone which is a hyperplane with constant orientation but non-constant multiplicity. Then $T$ near $x$ has $C^{1,1}$ tangentially immersed boundary.}

\bigskip

We suspect the proof of this conjecture is standard regularity for elliptic systems, but many details are needed; we leave this for future work. Note that more general versions of Theorems \ref{hyperplanetangentconeregularity},\ref{main} (namely, if we only assume $T$ has $C^{1,\alpha}$ tangentially immersed boundary with $\alpha \in (0,1]$) hold in case $n=2,$ for two-dimensional area-minimizing currents in space. For this, see Theorems \ref{spaceplanetangentconeregularity},\ref{spacemain}. 

\bigskip

It remains to investigate the geometric structure of co-dimension one area-minimizing locally rectifiable currents $T$ with $C^{1,\alpha}$ tangentially immersed boundary where $\alpha \in (0,1]$ near points $x$ in the support of $\partial T$ so that $T$ at $x$ has a tangent cone which is a sum of half-hyperplanes with constant orientation after rotation (see Definition \ref{conesdefinition}). To this end, we perhaps expect to at least assume $\partial T$ has co-oriented mean curvature $H_{\partial T} = h \nu_{T}$ with $h$ Lipschitz. An very optimistic conjecture is that the support of $T$ near such an $x$ is a finite union of $C^{1}$ hypersurfaces-with-boundary, so that then Theorem \ref{main} applies. It may at least be possible to show that $T$ has unique tangent cone at every $x$ in the support of $\partial T.$ Much more work is needed.

\bigskip

There is no hope to extend the present results in the case of general $n$-dimensional area-minimizing currents in $\R^{n+k}.$ It is well-known that the boundary regularity of \cite{HS79},\cite{W83} do not hold in higher co-dimensions. Furthermore, we refer the reader to \S 1.2 of \cite{R15d}, which gives examples showing that not even the partial boundary regularity result Theorem \ref{partialboundaryregularity} (Theorem 2.1 of \cite{R15d}) holds in higher co-dimensions, or even for merely \emph{stable} two-dimensional currents in space.

\subsection{Applications}

We now discuss the thread problem and the $c$-Plateau problem, both of which yield categories of area-minimizing currents for which the present results are relevant.

\bigskip

First, we describe the thread problem as defined in \cite{E89}. We say that $T,$ an $n$-dimensional locally rectifiable current in $\R^{n+1},$ is a minimizer of the thread problem with respect to $\Gamma,$ an $(n-1)$-dimensional locally rectifiable current in $\R^{n+1},$ if $\mathbf{M}_{U}(T) \leq \mathbf{M}_{U}(S)$ whenever $U \subset \R^{n+1}$ is an open set, $\mathbf{M}_{U}$ is the usual mass on currents in $U,$ and $S$ is an $n$-dimensional locally rectifiable current in $\R^{n+1}$ such that the support of $S-T$ is contained in $U$ and $\mathbf{M}_{U}(\partial S - \Gamma) = \mathbf{M}_{U}(\partial T-\Gamma).$

\bigskip

Naturally, minimizers of the thread problem are area-minimizing; see (2) of \cite{E89}. Theorem 2.3 of \cite{E89} states that if $T$ is a minimizer of the thread problem with respect to $\Gamma,$ then the \emph{free boundary} $\Sigma = \partial T - \Gamma$ has away from the support of $\Gamma$ co-oriented mean curvature $H_{\Sigma} = \frac{-1}{\lambda_{\Sigma}} \nu_{T}$ where $\lambda_{\Sigma}$ is a positive number; see in particular (2.8) of \cite{E89}. Also, Theorem 3.4 of \cite{E89} concludes in case $n=2$ that $T$ away from the support of $\Gamma$ has $C^{1,1}$ tangentially immersed boundary. The (tangentially immersed) regularity of the free boundary is unknown for general dimensions $n>2.$

\bigskip

Second, we describe the $c$-Plateau problem as defined in \cite{R15a}. We say that $T,$ an $n$-dimensional locally rectifiable current in $\R^{n+1},$ is a solution to the $c$-Plateau problem with respect to $\Gamma,$ an $(n-1)$-dimensional locally rectifiable current in $\R^{n+1},$ if $\partial T = \Sigma+\Gamma$ where $\Sigma,\Gamma$ have disjoint supports and
$$\mathbf{M}_{U}(T)+c \mathbf{M}_{U}(\partial T)^{\frac{n}{n-1}} \leq \mathbf{M}_{U}(S)+c \mathbf{M}_{U}(\partial S)^{\frac{n}{n-1}}$$
whenever $U \subset \R^{n+1}$ is an open set and $S$ is an $n$-dimensional locally rectifiable current in $\R^{n+1}$ such that $(\partial S-\Gamma),\Gamma$ have disjoint supports and the support of $S-T$ is contained in $U.$

\bigskip

Solutions to the $c$-Plateau problem are not only area-minimizing, they solve the thread problem; see Theorem 5.2 of \cite{R15a}. However, it follows directly from the definition that if $T$ is a solution to the $c$ Plateau problem with respect to $\Gamma,$ then the \emph{free boundary} $\Sigma = \partial T-\Gamma$ has co-oriented mean curvature $H_{\Sigma} = -\frac{1}{c(\frac{n}{n-1}) \mathbf{M}_{U}(\partial T)^{\frac{1}{n-1}}} \nu_{T}$; see (6.1) of \cite{R15a}. In case $n=2,$ we can use Theorem 3.4 of \cite{E89} to conclude that $T$ away from the support of $\Gamma$ has $C^{1,1}$ tangentially immersed boundary; see Theorem 4 of \cite{R16}. The (tangentially immersed) regularity of the free boundary is unknown for general dimensions $n>2.$

\subsection{Outline}

We begin in \S 2 by listing the notation we shall use. Next, in \S 3 we carefully define what it means for $T$ to have $C^{k,\alpha}$ tangentially immersed boundary for $k \geq 1$ in Definition \ref{immersedboundary}, and eventually state and prove Theorem \ref{hyperplanetangentconegraph}; to this we also include an appendix, where we prove a technical result Lemma \ref{sardslemma} needed to prove Lemma \ref{hyperplanetangentcone}. In \S 4 we define what it means for $\partial T$ to have co-oriented mean curvature in Definition \ref{cmcboundary}, and give a few basic related results. Finally, in \S 5 we state and prove Theorems \ref{hyperplanetangentconeregularity},\ref{main}.

\section{Notation}

We list basic notation and terminology we shall use throughout.

\begin{itemize}
 \item $\N,\R$ will denote the natural and real numbers respectively. We shall let $n \in \N$ with $n \geq 2.$ In this section we will let $\tilde{n} \in \{ 1,\ldots,n \}.$
 \item We shall typically write points $x = (x_{1},\ldots,x_{n+1}) \in \R^{n+1}.$ Depending on context, we shall let 
 $$\R^{\tilde{n}} = \{ (x_{1},\ldots,x_{\tilde{n}},0,\ldots,0) \in \R^{n+1} : x_{1},\ldots,x_{\tilde{n}} \in \R \}.$$
 We shall typically write $z = (z_{1},\ldots,z_{n-1}) = (z_{1},\ldots,z_{n-1},0,0,0) \in \R^{n-1}$ and $y = (y_{1},\ldots,y_{n}) =(y_{1},\ldots,y_{n},0,0) \in \R^{n}.$ For $z = (z_{1},\ldots,z_{n-1}) \in \R^{n-1}$ we will often make use of the affine plane
 $$\{z\} \times \R^{2} = \{ (z_{1},\ldots,z_{n-1},x_{n},x_{n+1}) \in \R^{n+1}: x_{n},x_{n+1} \in \R \}.$$
 We will let $0$ denote the zero vector in different dimensions, depending on context. Hence, we will also make use of $\{0\} \times S^{1} = \{ (0,\vec{x}_{n},x_{n+1}) \in \R^{n+1}: (x_{n},x_{n+1}) \in S^{1} \}.$
 \item Let $e_{1},\ldots,e_{n+1}$ be the standard basis vectors for $\R^{n+1}.$ 
 \item For $A \subseteq \R^{n+1},$ let $\clos A$ be the closure of $A.$ $U \subseteq_{o} \R^{n+1}$ shall denote $U$ is open in $\R^{n+1}.$ For $U,\tilde{U} \subseteq_{o} \R^{n+1},$ we write $\tilde{U} \subset \subset U$ if $\clos \tilde{U}$ is compact with $\clos \tilde{U} \subset U.$
 \item We shall let $B_{\rho}(x)$ be the open ball in $\R^{n+1}$ of radius $\rho>0$ centered at $x.$ For $x \in \R^{\tilde{n}},$ we write $B^{\tilde{n}}_{\rho}(x) = B_{\rho}(x) \cap \R^{\tilde{n}}.$
 \item For $x \in \R^{n+1}$ and $\lambda > 0,$ we let $\eta_{x,\lambda} : \R^{n+1} \rightarrow \R^{n+1}$ be the map $\eta_{x,\lambda}(\tilde{x}) = \frac{\tilde{x}-x}{\lambda}.$ We shall often make use of $\eta_{-x,1},$ which is translation to the right by $x.$ 
 \item For two sets $A,B \subseteq \R^{n+1},$ we denote the Hausdorff distance between $A$ and $B$ by $\dist_{\HH}(A,B).$
 \item For $U \subseteq_{o} \R^{n+1},$ we say $M \subset U$ is a hypersurface-with-boundary (embedded) in $U$ if $M$ is a hypersurface (embedded) in $U$ and $M$ attains its topological boundary $\partial M = (\clos M) \setminus M$ in the strong sense in $U$: $M$ near every point of $(\partial M) \cap U$ is diffeomorphic to the graph of a smooth function defined over a half-hyperplane. We similarly define a $C^{1},$ or more generally a $C^{k,\alpha}$ for $\alpha \in (0,1]$ and $k \in \N,$ hypersurface-with-boundary.
 \item We let $\ast : \bigwedge_{n} \R^{n+1} \rightarrow \R^{n}$ be the Hopf map
 $$\ast \left( \sum_{i=1}^{n+1} x_{i} (-1)^{i-1} e_{1} \wedge \ldots \wedge e_{i-1} \wedge e_{i+1} \wedge \ldots \wedge e_{n+1} \right) = \sum_{i=1}^{n+1} x_{i} e_{i}.$$
 Note that $\ast(e_{1} \wedge \ldots \wedge e_{n}) =(-1)^{n} e_{n+1}.$
 \item We shall let $D$ denote differentiation generally over $\R^{n+1}$ or $\R^{\tilde{n}},$ depending on context. In the proof of Lemma \ref{halfregulartangentcones} we will let $\D$ denote differentiation over $\R^{n-1},$ for emphasis. 
 \item $\HH^{\tilde{n}}$ shall denote $\tilde{n}$-dimensional Hausdorff measure in $\R^{n+1}.$ We also let $\omega_{\tilde{n}} = \HH^{\tilde{n}}(B^{\tilde{n}}_{1}(0)).$ 
 \item A smooth Jordan arc is a curve $\gamma \in C([a,b];\R^{n+1}) \cap C^{\infty}((a,b);\R^{n+1})$ which is injective over $(a,b)$; we may have $\gamma(a)=\gamma(b).$ A smooth closed Jordan curve is an injective curve $\gamma \in C^{\infty}(S^{1};\R^{n+1}).$ A continuous Jordan arc is a curve $\gamma \in C([a,b];\R^{n+1})$ injective over $(a,b),$ and a continuous closed Jordan curve is an injective curve $\gamma \in C(S^{1};\R^{n+1}).$
\end{itemize}
 
We now give notation related to currents in $\R^{n+1}.$ For a thorough introduction to currents, see \cite{F69},\cite{S83}.
 
\begin{itemize}
 \item Recall that $\mathcal{D}^{\tilde{n}}(U)$ denotes for $U \subseteq_{o} \R^{n+1}$ the smooth $\tilde{n}$-forms compactly supported in $U.$
 \item For $T$ a current in $U \subseteq_{o} \R^{n+1}$ and $f:U \rightarrow \R^{n+1},$ we denote $f_{\#} T$ the push-forward current of $T$ by $f$; we shall frequently make use of $\eta_{x,\lambda \#} T.$
 \item We say a current $\C$ is a cone if $\eta_{0,\lambda \#} \C = \C$ for every $\lambda>0.$
 \item Given an orientable $\tilde{n}$-dimensional submanifold $M \subset \R^{n+1},$ we denote $\lsem M \rsem$ the associated multiplicity one current, given an orientation. 
 \item Denote by $\E^{\tilde{n}}$ the $\tilde{n}$-dimensional current in $\R^{n+1}$ given by $\E^{\tilde{n}}(\omega) = \int_{\R^{\tilde{n}}} \langle \omega, e_{1} \wedge \ldots \wedge e_{\tilde{n}} \rangle \ d \HH^{\tilde{n}}$ for $\omega \in \mathcal{D}^{\tilde{n}}(\R^{n+1}).$
 \item For $U \subseteq_{o} \R^{n+1}$ and $T$ an $\tilde{n}$-dimensional current in $U,$ we let $\mu_{T}$ denote the associated mass measure of $T.$ This is given for $\tilde{U} \subseteq_{o} U$ by $\mu_{T}(\tilde{U}) = \sup_{\omega \in \mathcal{D}^{\tilde{n}}(\tilde{U}),|\omega| \leq 1} T(\omega).$ As usual, we set $\spt T = \spt \mu_{T}.$
 
 \bigskip
 
 For $A$ a $\mu_{T}$-measurable set, we let $T \res A$ denote the restriction current $(T \res A)(\omega) = \int_{A} <\omega,\vec{T}> \ d \mu_{T}$ for $\omega \in \mathcal{D}^{\tilde{n}}(U),$ where $\vec{T}$ is the orientation vector of $T.$ 

 \bigskip
 
Given $x \in U,$ we denote the density of $T$ at $x$ by $\Theta_{T}(x) = \lim_{\rho \searrow 0} \frac{\mu_{T}(B_{\rho}(x))}{\omega_{\tilde{n}} \rho^{\tilde{n}}}$ for $\omega_{\tilde{n}}=\HH^{\tilde{n}}(B^{\tilde{n}}_{\rho}(0)),$ wherever this limit exists.   
 \item Given $U \subseteq_{o} \R^{n+1},$ we let $\I_{\tilde{n},loc}(U)$ be the set of $\tilde{n}$-dimensional currents $T$ so that $T,\partial T$ are respectively $\tilde{n}$- and $(\tilde{n}-1)$-rectifiable integer multiplicity. 
 
 \bigskip
 
 For $T \in \I_{\tilde{n},loc}$ we let $T_{x}T$ denote the approximate tangent space of $T$ for the $\mu_{T}$-almost-every $x \in U$ such that this space exists; naturally, we let $T^{\perp}_{x}T$ denote the orthogonal complement of $T_{x}T$ in $\R^{n+1}.$ 
 \item For $T \in \I_{\tilde{n},loc}(U),$ we denote $\delta T$ to be the first variation of mass, given by 
 $$\delta T (X) = \int \dive{T}{X} \ d \mu_{T}$$ 
 for $X \in C_{c}^{1}(U;\R^{n+1}).$ 
 
 \bigskip
 
We say that $T$ has mean curvature $H_{T}:U \rightarrow \R^{n+1}$ if $H_{T}$ is $\mu_{T}$-measurable and if
$$\delta T(X) = \int X \cdot H_{T} \ d \mu_{T}$$
 for every $X \in C_{c}^{1}(U;\R^{n+1})$
 \item For $T \in \I_{\tilde{n},loc}(U),$ we let $\reg T$ denote the regular set of $T$: the set of $x \in \spt T$ so that there is a $\rho > 0$ such that $T \res B_{\rho}(x) = \theta \lsem M \rsem$ for $\theta \in \N$ and $M$ an $\tilde{n}$-dimensional orientable (embedded) $C^{1}$ submanifold of $B_{\rho}(x).$ We define the singular set $\sing T = \spt T  \setminus \reg T.$
  \item We say $T \in \I_{\tilde{n},loc}(U)$ is area-minimizing if $\mu_{T}(\tilde{U}) \leq \mu_{R}(\tilde{U})$ whenever $\tilde{U} \subset \subset U$ and $R \in \I_{\tilde{n},loc}(U)$ with $\partial R = \partial T$ and $\spt (T-R) \subset \tilde{U}.$
 \item For $T \in \I_{\tilde{n},loc}(U)$ area-minimizing there is, by Lemma 3.1 of \cite{B77} (see as well (2.10) of \cite{E89}), a $\mu_{\partial T}$-measurable vectorfield $\nu_{T}:U \rightarrow \R^{n+1}$ satisfying $|\nu_{T}| \leq 1$ for $\mu_{\partial T}$-almost-everywhere so that $$\delta T (X) = \int \nu_{T} \cdot X \ d \mu_{\partial T}$$ for every $X \in C^{1}_{c}(U;\R^{n+1}).$ We call $\nu_{T}$ the generalized outward pointing normal of $\partial T$ with respect to $T.$ Note that since $|\delta T (X)| \leq \int |X \wedge \vec{\partial T}| \ d \mu_{\partial T}$ by Lemma 3.1 of \cite{B77} (see also (2.9) of \cite{E89}) for $X \in C^{1}_{c}(U;\R^{n+1}),$ we conclude $\nu_{T}(x) \in T^{\perp}_{x} \partial T$ for $\mu_{\partial T}$-almost-every $x \in U.$ 
\end{itemize}

\section{Tangentially Immersed Boundary}

We begin by defining precisely what it means to have $C^{k,\alpha}$ tangentially immersed boundary. The goal of this section is to prove Theorem \ref{hyperplanetangentconegraph}. In the subsequent sections, \S 4,5, we will also assume our boundaries have \emph{co-oriented mean curvature}; see Definition \ref{cmcboundary}. The first part of this section will consist of giving more general lemmas as well as recalling some needed previous results.

\begin{definition} \label{immersedboundary} Let $U \subseteq_{o} \R^{n+1},$ $k \in \N,$ and $\alpha \in (0,1].$ We define $\TI^{k,\alpha}_{n,loc}(U)$ to be the set of area-minimizing $T \in \I_{n,loc}(U)$ so that $\partial T$ is locally $C^{k,\alpha}$ tangentially immersed: for every $x \in \spt \partial T$ there is $\rho>0,$ an orthogonal rotation $Q,$ and $N \in \N$ so that
$$\partial T \res B_{\rho}(x) = (-1)^{n} \sum_{\ell=1}^{N} m_{\ell} \Big[ (\eta_{-x,1} \circ Q \circ \Phi_{T,\ell})_{\#}(\E^{n-1} \res B^{n-1}_{\rho}(0)) \Big] \res B_{\rho}(x),$$
where for each $\ell =1,\ldots,N$ we have $m_{\ell} \in \N,$ and $\Phi_{T,\ell} \in C^{k,\alpha}(B^{n-1}_{\rho}(0);\R^{n+1})$ is the map
$$\Phi_{T,\ell}(z) = (z,\varphi_{T,\ell}(z),\psi_{T,\ell}(z)),$$
where $\varphi_{T,\ell},\psi_{T,\ell} \in C^{k,\alpha}(B^{n-1}_{\rho}(0))$ satisfy
$$\varphi_{T,\ell}(0)=\psi_{T,\ell}(0)=0 \text{ and } D \varphi_{T,\ell}(0) = D \psi_{T,\ell}(0)=0.$$ \end{definition}

Observe that we could define what it means for a current to have $C^{k,\alpha}$ tangentially immersed boundary in general. But we include the requirement that $T \in \TI^{k,\alpha}_{n,loc}(U)$ must be area-minimizing for future brevity. Observe that if $T \in \TI^{1,\alpha}_{n,loc}(U)$ then the approximate tangent space $T_{x}\partial T$ exists for every $x \in \spt \partial T.$

\bigskip

In order to state and prove the results of this section, we need the following lemma, leading to Definition \ref{conesdefinition}. 

\begin{lemma} \label{cones} Suppose $\C \in \I_{n,loc}(\R^{n+1})$ is an area-minimizing cone with $\partial \C = m (-1)^{n} Q_{\#} \E^{n}$ for some $m \in \N$ and an orthogonal rotation $Q.$ Then $\C$ is of one of the following two forms:
\begin{enumerate}
  \item[(1)] There is $N \in \{ 1,\ldots,m \}$ and distinct orthogonal rotations $Q_{1},\ldots,Q_{N}$ about $\R^{n-1}$ so that
  $$\C = \sum_{k=1}^{N} m_{k} (Q \circ Q_{k})_{\#} (\E^{n} \res \{ y \in \R^{n}: y_{n}>0 \}),$$ 
  where $m_{1},\ldots,m_{N}$ are positive integers with $\sum_{k=1}^{N} m_{k} = m.$
  \item[(2)] There is $\theta \in \N$ so that 
  $$\C = Q_{\#} \Big( (m+\theta) \E^{n} \res \{ y \in \R^{n}: y_{n}>0 \} + \theta \E^{n} \res \{ y \in \R^{n}: y_{n}<0 \} \Big).$$
\end{enumerate}
\end{lemma}

{\bf Proof:} By the boundary regularity for area-minimizing currents of \cite{HS79} and \cite{W83}, we get that the density $\Theta_{\C}(x)$ is constant for $x \in \R^{n-1}.$ The lemma then follows from Theorem 5.1 of \cite{B77}. $\square$

\bigskip

We thus give the following definition:

\begin{definition} \label{conesdefinition} Suppose $\C \in \I_{n,loc}(\R^{n+1})$ If $\C$ is as in (1) of Lemma \ref{cones} (with $Q,Q_{1},\ldots,Q_{k}$ orthogonal rotations, $m,N,m_{1},\ldots,m_{N} \in \N$), then we say that $\C$ is a \emph{sum of half-hyperplanes with constant orientation after rotation}. If $\C$ is as in (2) of Lemma \ref{cones} (with $Q$ an orthogonal rotation and $m,\theta \in \N$), then we say that $\C$ is a \emph{hyperplane with constant orientation but non-constant multiplicity}.  \end{definition}

The first step in proving Theorem \ref{hyperplanetangentconegraph}, the main result of this section, is to show that $T \in \TI^{1,\alpha}_{n,loc}(U)$ has tangent cones at every point $x \in \spt \partial T.$ We also categorize the tangent cones, using Lemma \ref{cones}. To do this requires proving a monotonicity formula.

\begin{theorem} \label{monotonicity} Let $U \subseteq_{o} \R^{n+1},$ $\alpha \in (0,1],$ and suppose $T \in \TI^{1,\alpha}_{n,loc}(U).$ For every $x \in \spt T$ the following hold:
\begin{enumerate}
 \item[(a)] For $0 < \sigma < \rho \leq \dist(x,\partial U)$
 $$\begin{aligned} 
 \frac{\mu_{T}(B_{\tau}(x))}{\tau^{n}} + \frac{1}{n} \int_{B_{\tau}(x)} & \left( \frac{1}{|\tilde{x}-x|^{n}}-\frac{1}{\tau^{n}} \right) (\tilde{x}-x) \cdot \nu_{T} \ d \mu_{\partial T}(\tilde{x}) \Big|^{\rho}_{\tau=\sigma} \\ & = \int_{B_{\rho}(x) \setminus B_{\sigma}(x)} \frac{|\proj{T_{\tilde{x}}^{\perp}T}{(\tilde{x}-x)}|^{2}}{|\tilde{x}-x|^{n+2}} \ d \mu_{T}(\tilde{x})
 \end{aligned}$$
where $\nu_{T}$ is the generalized outward pointing unit normal of $\partial T$ with respect to $T.$
 \item[(b)] The density $\Theta_{T}(x)$ of $T$ at $x$ exists.
 \item[(c)] There exists an area minimizing oriented tangent cone of $T$ at $x.$ Every tangent cone of $T$ at $x$ is either a sum of half-hyperplanes with constant orientation after rotation or a hyperplane with constant orientation but non-constant multiplicity; see Definition \ref{conesdefinition}.
\end{enumerate}
\end{theorem}

{\bf Proof:} We wish to apply Theorem 3.3 and Corollary 3.5 of \cite{B77}. For this, we need only to check that for $\rho \in (0,\dist(x,\partial U))$
$$\int_{B_{\rho}(x)} \frac{|\proj{T^{\perp}_{\tilde{x}}T}{(\tilde{x}-x)}|}{|\tilde{x}-x|^{n+1}} \ d \mu_{\partial T}(\tilde{x}) < \infty.$$
This readily follows from $\partial T$ being $C^{1,\alpha}$ tangentially immersed. From this we get (a),(b) and that $T$ has an area-minimizing oriented tangent cone at $x,$ with every tangent cone $\C$ of $T$ at $x$ being area-minimizing. The remainder of (c) follows from Lemma \ref{cones} and $T \in \TI^{1,\alpha}_{n,loc}(U).$ $\square$

\bigskip

Given Definition \ref{conesdefinition}, and having established the existence and structure of the tangent cones along the boundary of any $T \in \TI^{1,\alpha}_{n,loc}(U),$ it is convenient now to give Theorem 2.1 of \cite{R15d}, a partial boundary regularity result. Note that the assumption that $T \in \TI^{1,\alpha}_{n,loc}(U)$ is not needed in Theorem 2.1 of \cite{R15d}, but only that $\partial T$ consists of a finite sum of $(n-1)$-dimensional $C^{1,\alpha}$ submanifolds which meet tangentially, with same orientation, at the boundary point considered. Nevertheless, for convenience we give here the more restrictive version.

\begin{theorem} \label{partialboundaryregularity} For every pair $m,\theta \in \N$ and $\alpha \in (0,1]$ there is $\delta=\delta(n,m,\theta,\delta) \in (0,1)$ so that the following holds:

\bigskip

Let $U \subseteq_{o} \R^{n+1}$ and suppose $T \in \TI^{1,\alpha}_{n,loc}(U)$ with $0 \in \spt \partial T.$ Also suppose that $T$ at $0$ has tangent cone
$$\C = (m+\theta) \E^{n} \res \{ y \in \R^{n}: y_{n}>0 \} + \theta \E^{n} \res \{ y \in \R^{n}: y_{n}<0 \}.$$
Then, for $\rho \in (0,\dist(0,\partial U))$ sufficiently small (depending on $T$), and with $\beta = \frac{\alpha}{4n+6},$ there is a function
$$u \in \left\{ 
\begin{aligned}
& C^{\infty}(\{ y \in B^{n}_{\rho}(0) : |y_{n}| > (\delta/\rho)^{\beta} |y|^{1+\beta} \}) \\
& C^{1,\beta}(\{ y \in B^{n}_{\rho}(0) : |y_{n}| \geq (\delta/\rho)^{\beta} |y|^{1+\beta} \})
\end{aligned} \right.$$
with $u(0)=0$ and $Du(0)=0$ so that
$$\begin{aligned}
T \res & B_{\rho}(0) \cap \left( \{ y \in B^{n}_{\rho}(0): |y_{n}| \geq (\delta/\rho)^{\beta} |y|^{1+\beta} \} \times \R \right) \\
= & \Big[ (m+\theta) F_{\#} (\E^{n} \res  \{ y \in B^{n}_{\rho}(0) : y_{n} \geq (\delta/\rho)^{\beta} |y|^{1+\beta} \}) \\
& + \theta F_{\#} (\E^{n} \res  \{ y \in B^{n}_{\rho}(0) : y_{n} \leq -(\delta/\rho)^{\beta} |y|^{1+\beta} \}) \Big] \\
& \res B_{\rho}(0) \cap \left( \{ y \in B^{n}_{\rho}(0): |y_{n}| \geq (\delta/\rho)^{\beta} |y|^{1+\beta} \} \times \R \right)
\end{aligned}$$
where $F(y) = (y,u(y)).$
\end{theorem}

{\bf Proof:} This is Theorem 2.1 of \cite{R15d}; note that we have reversed the roles of $m,\theta$ found therein. $\square$

\bigskip

We can now proceed to the main result of this section, Theorem \ref{hyperplanetangentconegraph}. Recall that the brunt of proving Theorem \ref{hyperplanetangentconegraph} involves giving the following Lemma \ref{hyperplanetangentcone}. In fact, Lemma \ref{hyperplanetangentcone} will be instrumental to the subsequent results. On account of this, we give a careful proof of Lemma \ref{hyperplanetangentcone}. This proof is technical, but very geometric. 

\begin{lemma} \label{hyperplanetangentcone} Let $U \subseteq_{o} \R^{n+1},$ $\alpha \in (0,1],$ and suppose $T \in \TI^{1,\alpha}_{n,loc}(U)$ with $0 \in \spt \partial T.$ If $T$ at $0$ has a tangent cone which is a hyperplane with constant orientation but non-constant multiplicity (see Definition \ref{conesdefinition}), then there is $\rho \in (0,\dist(0,\partial U))$ so that $T$ at every $x \in \spt \partial T \cap B_{\rho}(0)$ has a unique tangent cone which is a hyperplane with constant orientation but non-constant multiplicity.

\bigskip

Furthermore suppose $T_{0} \partial T = \R^{n-1}.$ Then $\rho \in (0,\dist(0,\partial U))$ can also be chosen so that for $\HH^{n-1}$-almost-every $z \in B^{n-1}_{\rho}(0)$ we have
$$\spt T \cap \Big( \{ z \} \times \R^{2} \Big) \cap B_{\rho}(0) = \Gamma \cup L$$
where $\Gamma$ is a continuous Jordan arc, smooth away from $\spt \partial T,$ with endpoints in $\partial B_{\rho}(0),$ and $L$ is a finite disjoint collection of smooth closed Jordan curves with $L \cap (\Gamma \cup \spt \partial T) = \emptyset$.
\end{lemma}

{\bf Proof:} The proof of the first part of the lemma essentially requires us to prove the second part. To this end, we analyze the cross-sections of the support of $T$ perpendicular to $T_{0} \partial T,$ in hopes that we can apply Lemma \ref{sardslemma} found in the appendix. This will require us to first make a number of geometric observations, using Theorem \ref{partialboundaryregularity} as well as the regularity theory from Theorem 1 of \cite{SS81}. 

\bigskip

Assume (after rotation) that $T$ at $0$ has tangent cone 
$$\C = (m^{\C}+\theta^{\C}) \E^{n} \res \{ y \in \R^{n}: y_{n}>0 \} + \theta^{\C} \E^{n} \res \{ y \in \R^{n}: y_{n} < 0 \}$$
for $m^{\C},\theta^{\C} \in \N.$ Also assume $\rho \in (0,\dist(0,\partial U))$ is sufficiently small so that the following three occur:

\bigskip

First, by Definition \ref{immersedboundary} for some $N,m_{1},\ldots,m_{N} \in \N$
$$\partial T \res B_{\rho}(0) = (-1)^{n} \sum_{\ell=1}^{N} m_{\ell} \Phi_{T,\ell \#}( \E^{n-1} \res B^{n-1}_{\rho}(0)) \res B_{\rho}(0)$$
where $\Phi_{T,\ell}(z) = (z,\varphi_{T,\ell}(z),\psi_{T,\ell}(z))$ with $\varphi_{T,\ell},\psi_{T,\ell} \in C^{1,\alpha}(B^{n-1}_{\rho}(0))$ satisfy $\varphi_{T,\ell}(0) = \psi_{T,\ell}(0)=0$ and $D \varphi_{T,\ell}(0) = D \psi_{T,\ell}(0)=0.$ 

\bigskip
 
Second, for each $z \in B^{n-1}_{\rho}(0)$ and $\ell \in \{ 1,\ldots,N \}$ there is an orthogonal rotation $Q_{\ell}^{z}$ so that 
\begin{equation} \label{orthogonalrotation}
\| Q_{\ell}^{z} - I \| < \frac{1}{8} \text{ and } Q_{\ell}^{z}(\R^{n-1}) = T_{\Phi_{T,\ell}(z)} \Phi_{T,\ell}(B^{n-1}_{\rho}(0)).
\end{equation}

\bigskip
 
Third, by Theorem \ref{partialboundaryregularity} we can assume that with $\beta = \frac{\alpha}{4n+6}$
\begin{equation} \label{outergraphs}
\begin{aligned}
T \res & B_{\rho}(0) \cap \left( \{ y \in B^{n}_{\rho}(0): |y_{n}| \geq (\delta/\rho)^{\beta} |y|^{1+\beta} \} \times \R \right) = \\
= & \Big[ (m^{\C}+\theta^{\C}) F_{\#} (\E^{n} \res \{ y \in B^{n}_{\rho}(0): y_{n} \geq (\delta/\rho)^{\beta} |y|^{1+\beta} \}) \\
& +\theta^{\C} F_{\#} (\E^{n} \res \{ y \in B^{n}_{\rho}(0): y_{n} \leq - (\delta/\rho)^{\beta} |y|^{1+\beta} \}) \Big] \\
& \res B_{\rho}(0) \cap \left( \{ y \in B^{n}_{\rho}(0): |y_{n}| \geq (\delta/\rho)^{\beta} |y|^{1+\beta} \} \times \R \right)
\end{aligned}
\end{equation}
where $F(y) = (y,u(y))$ for a function
$$u \in \left\{
\begin{aligned}
& C^{\infty}(\{ y \in B^{n}_{\rho}(0): |y_{n}| > (\delta/\rho)^{\beta} |y|^{1+\beta} \}) \\
& C^{1,\beta}(\{ y \in B^{n}_{\rho}(0): |y_{n}| \geq (\delta/\rho)^{\beta} |y|^{1+\beta} \})
\end{aligned} \right.$$
satisfying $u(0)=0$ and $Du(0)=0.$

\bigskip

To prove the first part of the lemma, suppose for contradiction there is $z \in B^{n-1}_{\rho/4}(0)$ so that (after relabeling) $T$ at $\Phi_{T,1}(z)$ has a tangent cone which is of a sum of half-hyperplanes with constant orientation after rotation (as in Definition \ref{conesdefinition}). Fixing such a $z \in B^{n-1}_{\rho/4}(0),$ our goal is to contradict Lemma \ref{sardslemma}. We begin by identifying the centers of the disjoint balls we will use in applying Lemma \ref{sardslemma}. For this, choose $N_{z} \in \{ 1,\ldots,N \}$ and $N^{(1)}_{z} \in \{ 1,\ldots,N_{z} \}$ so that by Theorem \ref{monotonicity} (after relabeling) the following three occur:
\begin{itemize}
 \item $\{ \Phi_{T,\ell}(z) \}_{\ell=1}^{N_{z}}$ is the set of distinct points in $\{ \Phi_{T,\ell}(z) \}_{\ell=1}^{N}.$ 
 \item $T$ at $\Phi_{T,\ell}(z)$ for each $\ell \in \{ 1,\ldots,N^{(1)}_{z} \}$ has a tangent cone 
\begin{equation}
\label{hhcone}
\C_{\ell} = \sum_{k=1}^{N_{\C_{\ell}}} m^{\C_{\ell}}_{k} (Q^{z}_{\ell} \circ Q_{\ell,k})_{\#} (\E^{n} \res \{ y \in \R^{n}: y_{n}>0 \})
\end{equation}
where $N_{\C_{\ell}},m^{\C_{\ell}}_{1},\ldots,m^{\C_{\ell}}_{N_{\C_{\ell}}} \in \N$ and $Q_{\ell,1},\ldots,Q_{\ell,N_{\C_{\ell}}}$ are distinct orthogonal rotations about $\R^{n-1}.$ Recall that the orthogonal rotations $Q^{z}_{\ell}$ satisfy \eqref{orthogonalrotation}.
 \item $T$ at $\Phi_{T,\ell}(z)$ for each $\ell$ in the possibly empty set $\{ N^{(1)}_{z}+1,\ldots,N_{z} \}$ has tangent cone 
\begin{equation}
\label{hcone}
\begin{aligned} \C_{\ell} = (Q^{z}_{\ell} \circ Q_{\ell,1})_{\#} \Big( (m^{\C_{\ell}} & +\theta^{\C_{\ell}}) \E^{n} \res \{ y \in \R^{n}: y_{n} > 0 \} \\
& + \theta^{\C_{\ell}} \E^{n} \res \{ y \in \R^{n}: y_{n} < 0 \} \Big) \end{aligned}
\end{equation}
where $m^{\C_{\ell}},\theta^{\C_{\ell}} \in \N$ and $Q_{\ell,1}$ is an orthogonal rotation about $\R^{n-1}.$ Note that $\C_{\ell}$ is unique by Theorem \ref{partialboundaryregularity}.
\end{itemize}

\bigskip

Next, we choose the correct disjoint balls using the points $\{ \Phi_{T,\ell}(z) \}^{N_{z}}_{\ell=1},$ so that we may apply Lemma \ref{sardslemma}. First, for any $x \in \R^{n+1}$ and $\sigma>0$ denote the bored-out ball:
\begin{equation} \label{puncturedball}
\hat{B}_{\sigma}(x) = B_{\sigma}(x) \setminus \{ \tilde{x} \in \R^{n+1}: |\proj{\{0\} \times \R^{2}}{(\tilde{x}-x)}| \leq \sigma/2 \}.
\end{equation}
Choose for each $\ell = 1,\ldots,N_{z}$ a $\sigma_{\ell}>0$ sufficiently small so that $\{ B_{\sigma_{\ell}}(\Phi_{T,\ell}(z)) \}_{\ell=1}^{N_{z}}$ is a disjoint collection of balls with $B_{\sigma_{\ell}}(\Phi_{T,\ell}(z)) \subset B_{\rho/2}(0),$ and so that the following two occur:

\bigskip

First, for each $\ell \in \{ 1,\ldots,N^{(1)}_{z} \}$ by Theorem 1 of \cite{SS81} and \eqref{orthogonalrotation},\eqref{hhcone},\eqref{puncturedball}
\begin{equation} \label{innerbadgraphs}
\begin{aligned} 
\eta_{\Phi_{T,\ell}(z),1 \#} & T \res \hat{B}_{\sigma_{\ell}}(0) \\
= & \sum_{k=1}^{N_{\C_{\ell}}} \sum_{j=1}^{m^{\C_{\ell}}_{k}} (Q^{z}_{\ell} \circ Q_{\ell,k} \circ F^{\ell}_{j,k})_{\#} (\E^{n} \res \{ y \in B^{n}_{\sigma_{\ell}}(0): y_{n}> \sigma_{\ell}/4 \}) \\
& \res \hat{B}_{\sigma_{\ell}}(0)
\end{aligned}
\end{equation}
where $F^{\ell}_{j,k}(y) = (y,u^{\ell}_{j,k}(y))$ for $u^{\ell}_{j,k} \in C^{\infty}(\{ y \in B^{n}_{\sigma_{\ell}}(0): y_{n}> \sigma_{\ell}/4 \}).$ For each $k \in \{ 1,\ldots,N_{\C_{\ell}} \}$ and $j,\tilde{j} \in \{ 1,\ldots,m^{\C_{\ell}}_{k} \},$ the graphs 
$$\gph{\{ y \in B^{n}_{\sigma_{\ell}}(0): y_{n}> \sigma_{\ell}/4 \}}{u^{\ell}_{j,k}} \text{ and } \gph{\{ y \in B^{n}_{\sigma_{\ell}}(0): y_{n}> \sigma_{\ell}/4 \}}{u^{\ell}_{\tilde{j},k}}$$
are either equal or disjoint. Meanwhile, if $k,\tilde{k} \in \{ 1,\ldots,N_{\C_{\ell}} \}$ with $k \neq \tilde{k}$ then
$$Q_{\ell,k}(\gph{\{ y \in B^{n}_{\sigma_{\ell}}(0): y_{n}> \sigma_{\ell}/4 \}}{u^{\ell}_{j,k}}) \cap Q_{\ell,\tilde{k}}(\gph{\{ y \in B^{n}_{\sigma_{\ell}}(0): y_{n}> \sigma_{\ell}/4 \}}{u^{\ell}_{\tilde{j},\tilde{k}}}) = \emptyset$$
for any $j \in \{ 1,\ldots,m^{\C_{\ell}}_{k} \}$ and $\tilde{j} \in \{ 1,\ldots,m^{\C_{\ell}}_{\tilde{k}} \}.$

\bigskip

Second, for each $\ell  \in \{ N^{(1)}_{z}+1,\ldots,N_{z} \}$ and with $\beta = \frac{\alpha}{4n+6}$
\begin{equation} \label{innergoodgraphs}
\begin{aligned} 
& \eta_{\Phi_{T,\ell}(z),1 \#} T \\
& \res \{ x \in B_{\sigma_{\ell}}(0): |x_{n}| > (\delta/\sigma_{\ell})^{\beta} |\proj{\R^{n}}{x}|^{1+\beta} \} \\
= & (Q_{\ell}^{z} \circ Q_{\ell,1})_{\#} \Big[ (m^{\C_{\ell}}+\theta^{\C_{\ell}}) F^{\ell}_{\#} (\E^{n} \res \{ y \in B^{n}_{\sigma_{\ell}}(0): y_{n}> (\delta/\sigma_{\ell})^{\beta} |y|^{1+\beta} \}) \\ 
& + \theta^{\C_{\ell}} F^{\ell}_{\#} (\E^{n} \res \{ y \in B^{n}_{\sigma_{\ell}}(0): y_{n} < - (\delta/\sigma_{\ell})^{\beta} |y|^{1+\beta} \}) \Big] \\
& \res \{ x \in B_{\sigma_{\ell}}(0): |x_{n}| > (\delta/\sigma_{\ell})^{\beta} |\proj{\R^{n}}{x}|^{1+\beta} \}
\end{aligned}
\end{equation}
where $F^{\ell}(y) = (y,u^{\ell}(y))$ for 
$$u^{\ell} \in \left\{
\begin{aligned}
& C^{\infty}( \{ y \in B^{n}_{\sigma_{\ell}}(0): |y_{n}| > (\delta/\sigma_{\ell})^{\beta} |y|^{1+\beta} \}) \\
& C^{1,\beta}( \{ y \in B^{n}_{\sigma_{\ell}}(0): |y_{n}| \geq (\delta/\sigma_{\ell})^{\beta} |y|^{1+\beta} \}).
\end{aligned} \right.$$

\bigskip

In order to apply Lemma \ref{sardslemma}, and obtain a contradiction, we must consider the support of $T$ over a cross-section $\{ \tilde{z} \} \times \R^{2}$ with suitably chosen $\tilde{z} \in \R^{n-1}$ near $z.$ To this end, define the set 
\begin{equation}
\label{set}
\mathcal{T} = \left\{ \tilde{z} \in B^{n-1}_{\rho}(0): 
\begin{aligned} 
& (\sing T \setminus \spt \partial T) \cap (\clos B_{\rho}(0)) \cap (\{\tilde{z}\} \times \R^{2}) = \emptyset \\
& \text{and } \ast \vec{T}(x) \notin \R^{n-1} \text{ for all } \\
& x \in (\spt T \setminus \spt \partial T) \cap (\clos B_{\rho}(0)) \cap (\{\tilde{z}\} \times \R^{2})
\end{aligned} \right\}
\end{equation}
Observe that $\ast \vec{T}(x) \notin \R^{n-1}$ for $x \in \reg T$ implies that $T_{x} T \cap (\{0\} \times \R^{2})$ is a one-dimensional subspace in $\R^{n+1}.$ Therefore $\HH^{n-1}(\mathcal{T})=0$ by Sard's theorem as well as interior regularity for co-dimension one area-minimizing currents. Fix $\tilde{z} \in B^{n-1}_{\min \{ \sigma_{1},\ldots,\sigma_{N_{z}} \}}(z) \cap \mathcal{T}.$ 

\bigskip

We conclude $\spt T \cap (\clos B_{\rho}(0)) \cap (\{\tilde{z}\} \times \R^{2})$ is a disjoint union of smooth Jordan arcs with endpoints only at $\partial B_{\rho}(0)$ or $\spt \partial T,$ together with a disjoint collection of smooth closed Jordan curves (we shall soon be more precise). This will allow us to apply Lemma \ref{sardslemma}. In order to do so, let
$$V \in C \Big( (\spt T \setminus \spt \partial T) \cap (\clos B_{\rho}(0)) \cap (\{\tilde{z}\} \times \R^{2}); \{ 0 \} \times S^{1} \Big)$$
be defined for $x \in (\spt T \setminus \spt \partial T) \cap (\clos B_{\rho}(0)) \cap (\{\tilde{z}\} \times \R^{2})$ by
\begin{equation} \label{sardsvector}
V(x) = \frac{\proj{\{0\} \times \R^{2}}{\ast \vec{T}(x)}}{|\proj{\{0\} \times \R^{2}}{\ast \vec{T}(x)}|}.
\end{equation}
Then by \eqref{outergraphs},\eqref{innerbadgraphs},\eqref{innergoodgraphs},\eqref{set} we have
$$\begin{aligned}
\spt T & \cap (\clos B_{\rho}(0)) \cap (\{\tilde{z}\} \times \R^{2}) \setminus K_{z} \\
& = \left[ \Gamma_{1} \cup \Gamma_{2} \cup \left( \bigcup_{\ell=1}^{N^{(1)}_{z}} \bigcup_{k=1}^{N_{\C_{\ell}}} \bigcup_{j=1}^{m^{\C_{\ell}}_{k}} \gamma^{\ell}_{j,k} \right) \cup \left( \bigcup_{\ell=N^{(1)}_{z}+1}^{N_{z}} (G^{\ell} \cup g^{\ell}) \right) \cup L \right] \setminus K_{z}
\end{aligned}$$
where 
$$K_{z} = \bigcup_{\ell=1}^{N_{z}} \{ x \in B_{\sigma_{\ell}}(\Phi_{T,\ell}(z)): |\proj{\{0\} \times \R^{2}} (x-\Phi_{T,\ell}(z))| < \sigma_{\ell}/2 \}$$
and the following hold:
\begin{enumerate}
 \item[(1)] $\Gamma_{1},\Gamma_{2}$ are Jordan arcs so that by \eqref{outergraphs} (for $\delta>0$ sufficiently small) with $\hat{B}_{\rho}(0)$ as in \eqref{puncturedball},
 $$\spt T \cap (\clos \hat{B}_{\rho}(0)) \cap (\{ \tilde{z} \} \times \R^{2}) = (\Gamma_{1} \cup \Gamma_{2}) \cap (\clos \hat{B}_{\rho}(0)).$$

\bigskip

We can parameterize $\Gamma_{1},\Gamma_{2}$ by arc-length so that
$$\Gamma_{1} \in \left\{ \begin{aligned} & C([0,\HH^{1}(\Gamma_{1})];(\clos B_{\rho}(0)) \cap (\{ \tilde{z} \} \times \R^{2})) \\ & C^{\infty}((0,\HH^{1}(\Gamma_{1})); B_{\rho}(0) \cap (\{ \tilde{z} \} \times \R^{2})) \end{aligned} \right.$$
and 
$$\Gamma_{2} \in \left\{ \begin{aligned} & C([0,\HH^{1}(\Gamma_{2})];(\clos B_{\rho}(0)) \cap (\{ \tilde{z} \} \times \R^{2})) \\ & C^{\infty}((0,\HH^{1}(\Gamma_{2})); B_{\rho}(0) \cap (\{ \tilde{z} \} \times \R^{2})), \end{aligned} \right.$$
$\Gamma_{1}(0),\Gamma_{2}(0) \in \partial B_{\rho}(0),$ 
$$\Gamma_{1}([0,\HH^{1}(\Gamma_{1}))) \cap \Gamma_{2}([0,\HH^{1}(\Gamma_{2}))) = \emptyset$$
(for this, see \eqref{densityouter} below), while $\Gamma_{1}(\HH^{1}(\Gamma_{1})) = \Phi_{T,\ell_{1}}(z_{1})$ and $\Gamma_{2}(\HH^{1}(\Gamma_{2})) = \Phi_{T,\ell_{2}}(z_{1})$ for some $\ell_{1},\ell_{2} \in \{ 1,\ldots,N \}.$ Thus, $\Gamma_{1}(\HH^{1}(\Gamma_{1})),\Gamma_{2}(\HH^{1}(\Gamma_{2})) \in K_{z}$ by \eqref{orthogonalrotation}.

\bigskip

With this parameterization, we have in the sense of Definition \ref{arcorientation} that $\Gamma_{1},V$ are positively oriented while $\Gamma_{2},V$ are negatively oriented, by \eqref{outergraphs}. The density of $T$ also satisfies by \eqref{outergraphs}
\begin{equation} \label{densityouter}
\begin{aligned}
 \Theta_{T}(x) & = m^{\C}+\theta^{\C} \text{ for } x \in \Gamma_{1}([0,\HH^{1}(\Gamma_{1}))), \\
 \Theta_{T}(x) & = \theta^{\C} \text{ for } x \in \Gamma_{2}([0,\HH^{1}(\Gamma_{2}))).
\end{aligned}
\end{equation}

 \item[(2)] For each $\ell = 1,\ldots,N^{(1)}_{z}$ we have that $\{ \gamma^{\ell}_{j,k} \}_{j=1,k=1}^{m^{\C_{\ell}}_{k},N_{\C_{\ell}}}$ is a collection of Jordan arcs so that by \eqref{innerbadgraphs}
$$\spt T \cap \hat{B}_{\sigma_{\ell}}(\Phi_{T,\ell}(z)) \cap (\{ \tilde{z} \} \times \R^{2}) = \bigcup_{k=1}^{N_{\C_{\ell}}} \bigcup_{j=1}^{m^{\C_{\ell}}_{k}} \gamma^{\ell}_{j,k} \cap \hat{B}_{\sigma_{\ell}}(\Phi_{T,\ell}(z)).$$

\bigskip

We can parameterize each arc
$$\gamma^{\ell}_{j,k} \in \left\{ \begin{aligned} & C([0,\HH^{1}(\gamma^{\ell}_{j,k})];(\clos B_{\rho}(0)) \cap (\{ \tilde{z} \} \times \R^{2})) \\ & C^{\infty}((0,\HH^{1}(\gamma^{\ell}_{j,k})); B_{\rho}(0) \cap (\{ \tilde{z} \} \times \R^{2})) \end{aligned} \right.$$
by arc-length so that  
$$\gamma^{\ell}_{j,k}(\HH^{1}(\gamma^{\ell}_{j,k})) \in \{ x \in B_{\sigma_{\ell}}(\Phi_{T,\ell}(z)): |\proj{\{0\} \times \R^{2}}{(x-\Phi_{T,\ell}(z))}| < \sigma_{\ell}/2 \}.$$
Meanwhile, $\gamma^{\ell}_{j,k}(0) \in (\partial B_{1}(0)) \cup K_{z}$ with $\gamma^{\ell}_{j,k} \cap (\partial B_{\sigma_{\ell}}(\Phi_{T,\ell}(z))) \cap (\{ \tilde{z} \} \times \R^{2}) \neq \emptyset.$ For each $k \in \{ 1,\ldots,N_{\C_{\ell}} \}$ the images $\gamma^{\ell}_{j,k}((0,\HH^{1}(\gamma^{\ell}_{j,k})))$ and $\gamma^{\ell}_{\tilde{j},k}((0,\HH^{1}(\gamma^{\ell}_{\tilde{j},k})))$ are either equal or disjoint for each $j,\tilde{j} \in \{ 1,\ldots,m^{\C_{\ell}}_{k} \},$ by \eqref{innerbadgraphs}. On the other hand, the images $\gamma^{\ell}_{j,k}((0,\HH^{1}(\gamma^{\ell}_{j,k})))$ and $\gamma^{\ell}_{\tilde{j},\tilde{k}}((0,\HH^{1}(\gamma^{\ell}_{\tilde{j},\tilde{k}})))$ are disjoint if $k \neq \tilde{k}.$

\bigskip

With this parameterization we have that $\gamma^{\ell}_{j,k},V$ are positively oriented (as in Definition \ref{arcorientation}), by \eqref{innerbadgraphs}. 

 \item[(3)] For each $\ell$ in the (possibly empty) set $\{ N^{(1)}_{z}+1,\ldots,N_{z} \}$ we have that $G^{\ell},g^{\ell}$ are Jordan arcs so that by \eqref{innergoodgraphs}
$$\spt T \cap \hat{B}_{\sigma_{\ell}}(\Phi_{T,\ell}(z)) \cap (\{ \tilde{z} \} \times \R^{2}) = (G^{\ell} \cup g^{\ell}) \cap \hat{B}_{\sigma_{\ell}}(\Phi_{T,\ell}(z)).$$

\bigskip

We can parameterize each arc
$$G^{\ell} \in \left\{ \begin{aligned} & C([0,\HH^{1}(G^{\ell})];(\clos B_{\rho}(0)) \cap (\{ \tilde{z} \} \times \R^{2})) \\ & C^{\infty}((0,\HH^{1}(G^{\ell})); B_{\rho}(0) \cap (\{ \tilde{z} \} \times \R^{2})) \end{aligned} \right.$$
and
$$g^{\ell} \in \left\{ \begin{aligned} & C([0,\HH^{1}(g^{\ell})];(\clos B_{\rho}(0)) \cap (\{ \tilde{z} \} \times \R^{2})) \\ & C^{\infty}((0,\HH^{1}(g^{\ell})); B_{\rho}(0) \cap (\{ \tilde{z} \} \times \R^{2})) \end{aligned} \right.$$
by arc-length so that 
$$G^{\ell}(\HH^{1}(G^{\ell})),g^{\ell}(\HH^{1}(g^{\ell})) \in \{ x \in B_{\sigma_{\ell}}(\Phi_{T,\ell}(z)): |\proj{\{0\} \times \R^{2}}{(x-\Phi_{T,\ell}(z))}| < \sigma_{\ell}/2 \}.$$
Meanwhile, $G^{\ell}(0),g^{\ell}(0) \in (\partial B_{1}(0)) \cup K_{z},$ 
$$\begin{aligned}
& G_{\ell} \cap (\partial B_{\sigma_{\ell}}(\Phi_{T,\ell}(z))) \cap (\{ \tilde{z} \} \times \R^{2}) \neq \emptyset \\
\text{and } & g_{\ell} \cap (\partial B_{\sigma_{\ell}}(\Phi_{T,\ell}(z))) \cap (\{ \tilde{z} \} \times \R^{2}) \neq \emptyset.
\end{aligned}$$
Moreover, the images $G^{\ell}((0,\HH^{1}(G^{\ell})))$ and $g^{\ell}((0,\HH^{1}(g^{\ell})))$ are disjoint (see \eqref{densityinner} below).

\bigskip

With this parameterization $G^{\ell},V$ are positively oriented while $g^{\ell},V$ are negatively oriented, by \eqref{innergoodgraphs}. The density of $T$ also satisfies
\begin{equation} \label{densityinner}
\begin{aligned}
 \Theta_{T}(x) & = m^{\C_{\ell}}+\theta^{\C_{\ell}} \text{ for } x \in G^{\ell}((0,\HH^{1}(G^{\ell}))), \\
 \Theta_{T}(x) & = \theta^{\C_{\ell}} \text{ for } x \in g^{\ell}((0,\HH^{1}(g^{\ell}))).
\end{aligned}
\end{equation}

 \item[(4)] $L$ is a disjoint union of smooth closed Jordan curves with
$$L \subset \{ x \in B_{\rho}(0): |\proj{\{0\} \times \R^{2}}{x}| < \rho/2 \} \cap (\{ \tilde{z} \} \times \R^{2}) \setminus \left( \bigcup_{\ell=1}^{N_{z}} B_{\sigma_{\ell}}(\Phi_{T,\ell}(z)) \right)$$
and so that 
$$L \cap \left[ \Gamma_{1} \cup \Gamma_{2} \cup \left( \bigcup_{\ell=1}^{N^{(1)}_{z}} \bigcup_{k=1}^{N_{\C_{\ell}}} \bigcup_{j=1}^{m^{\C_{\ell}}_{k}} \gamma^{\ell}_{j,k} \right) \cup \left( \bigcup_{\ell=N^{(1)}_{z}+1}^{N_{z}} (G^{\ell} \cup g^{\ell}) \right) \right] = \emptyset.$$
\end{enumerate}
Setting
$$\mathcal{G} = \Gamma_{1} \cup \Gamma_{2} \cup \left( \bigcup_{\ell=1}^{N^{(1)}_{z}} \bigcup_{k=1}^{N_{\C_{\ell}}} \bigcup_{j=1}^{m^{\C_{\ell}}} \gamma^{\ell}_{j,k} \right) \cup \left( \bigcup_{\ell=N^{(1)}_{z}+1}^{N_{z}} (G^{\ell} \cup g^{\ell}) \right)$$
means that we can apply Lemma \ref{sardslemma}$(a)$ with $V$ as in \eqref{sardsvector} in order to contradict that $\Gamma_{1},V$ are positively oriented while $\Gamma_{2},V$ are negatively oriented. This shows the first part of the theorem.

\bigskip

For the second part of the theorem, our aim is to apply Lemma \ref{sardslemma}(b). Choose any $z \in B^{n-1}_{\rho/4}(0) \cap \mathcal{T}$ with $\mathcal{T}$ as in \eqref{set}. Take again $N_{z} \in \{ 1,\ldots,N \}$ so that (after relabeling) $\{ \Phi_{T,\ell}(z) \}_{\ell=1}^{N_{z}}$ is the collection of distinct points in $\{ \Phi_{T,\ell}(z) \}_{\ell=1}^{N}.$  Again, we choose a disjoint collection of balls $\{ B_{\sigma_{\ell}}(\Phi_{T,\ell}(z)) \}_{\ell=1}^{N_{z}}$ with $B_{\sigma_{\ell}}(\Phi_{T,\ell}(z)) \subset B_{\rho/2}(z)$ and $\sigma_{\ell}>0$ sufficiently small so that \eqref{innergoodgraphs} holds now for each $\ell \in \{ 1,\ldots,N_{z} \}.$

\bigskip

In this case we conclude
$$\spt T \cap (\clos B_{\rho}(0)) \cap (\{ z \} \times \R^{2}) = \Gamma_{1} \cup \Gamma_{2} \cup \left( \bigcup_{\ell=1}^{N_{z}} (G^{\ell} \cup g^{\ell}) \right) \cup L$$
where the following hold:
\begin{itemize}
 \item $\Gamma_{1},\Gamma_{2}$ are as in (1) above, in the proof of the first part of the theorem, but of course with $\tilde{z}=z.$ 
 \item For each $\ell \in \{ 1,\ldots,N_{z} \},$ the arcs $G^{\ell},g^{\ell}$ are as in (3) above, but of course with $\tilde{z}=z.$ In this case we further conclude by \eqref{innergoodgraphs} that $G^{\ell}(\HH^{1}(G^{\ell}))=g^{\ell}(\HH^{1}(g^{\ell})) = \Phi_{T,\ell}(z).$
 \item $L$ is a disjoint collection of smooth closed Jordan curves as in (4) above, disjoint from the arcs $\Gamma_{1},\Gamma_{2},\{ G^{\ell},g^{\ell} \}_{\ell=1}^{N_{z}}.$
\end{itemize}
With $V$ as in \eqref{sardsvector} (but of course again with $\tilde{z}=z$) and
$$\mathcal{G} = \Gamma_{1} \cup \Gamma_{2} \cup \bigcup_{\ell=1}^{N_{z}} (G^{\ell} \cup g^{\ell}),$$
then we may apply Lemma \ref{sardslemma}(b), since $\Gamma_{1},V$ are positively oriented while $\Gamma_{2},V$ are negatively oriented. 

\bigskip

Assume for contradiction that, in applying Lemma \ref{sardslemma}(b), we conclude
$$\mathcal{G} = \Gamma \cup \bigcup_{\ell=1}^{N^{(2)}_{z}} L_{\ell}$$
where $\Gamma \subset (\clos B_{\rho}(0)) \cap (\{ z \} \times \R^{2})$ is a continuous Jordan arc, smooth away from the points $\{ \Phi_{T,\ell}(z) \}_{\ell=1}^{N_{z}}$ and with endpoints in $\partial B_{\rho}(0).$ Meanwhile, for each $\ell \in \{ 1,\ldots,N^{(2)}_{z} \}$ we have that $L_{\ell} \subset B_{\rho/2}(z) \cap (\{ z \} \times \R^{2})$ is a continuous closed Jordan curve, smooth away from the points $\{ \Phi_{T,\ell}(z) \}_{\ell=1}^{N_{z}}.$ Furthermore, $\Gamma,L_{1},\ldots,L_{N^{(2)}_{z}}$ have pairwise disjoint images.

\bigskip

Consider $L_{1},$ then $L_{1} \cap \{ \Phi_{T,\ell}(z) \}_{\ell=1}^{N_{z}} \neq \emptyset$ by Lemma \ref{sardslemma}(b). We can suppose (after relabeling) that $\Phi_{T,1}(z) \in L_{1}.$ This exactly means $G^{1},g^{1} \subset L_{1}.$ On the other hand, $G^{1}((0,\HH^{1}(G^{1})))$ and $g^{1}((0,\HH^{1}(g^{1})))$ are disjoint by (3) above. This implies that (after relabeling) $\Phi_{T,2}(z) \in L_{1}$ with $G^{1} = g^{2}$ (see the proof of Lemma \ref{sardslemma}, particularly in concluding \eqref{sards1},\eqref{sards2}). We therefore have $G^{2} \subset L_{1}$ as well. 

\bigskip

Arguing iteratively, we conclude there is $N^{L_{1}}_{z} \in \{ 2,\ldots,N_{z} \}$ so that (after relabeling)
$$\begin{aligned}
G^{1} & = g^{2} \\
& \vdots \\
G^{N^{L_{1}}_{z}-1} & = g^{N^{L_{1}}_{z}} \\
G^{N^{L_{1}}_{z}} & = g^{1}.
\end{aligned}$$
But then \eqref{densityinner} implies
$$\begin{aligned}
m^{1}+\theta^{1} & = \theta^{2} \\
& \vdots \\
m^{N^{L_{1}}_{z}-1}+\theta^{N^{L_{1}}_{z}-1} & = \theta^{N^{L_{1}}_{z}} \\
m^{N^{L_{1}}_{z}}+\theta^{N^{L_{1}}_{z}} & = \theta^{1},
\end{aligned}$$
which gives
$$\theta^{1} < m^{1}+\theta^{1} = \theta^{2} \leq \ldots \leq \theta^{N^{L_{1}}_{z}} < m^{N^{L_{1}}_{z}}+\theta^{N^{L_{1}}_{z}} = \theta^{1},$$
a contradiction. We conclude the second part of the theorem. $\square$

\bigskip

We are now ready to state and prove our first main result.

\begin{theorem} \label{hyperplanetangentconegraph} Let $U \subseteq_{o} \R^{n+1},$ $\alpha \in (0,1],$ and $T \in \TI^{1,\alpha}_{n,loc}(U).$ Suppose $x \in \spt \partial T$ and that $T$ at $x$ has tangent cone which is a hyperplane with constant orientation but non-constant multiplicity (as in Definition \ref{conesdefinition}). Then there is a $\rho \in \dist (0,\dist(x,\partial U))$ and a solution to the minimal surface equation $u \in C^{\infty}(B^{n}_{\rho}(0))$ with $u(0)=0$ and $Du(0)=0$ such that
$$\spt T \cap B_{\rho}(x) = \eta_{-x,1}(Q(\gph{B^{n}_{\rho}(0)}{u})) \cap B_{\rho}(x)$$
for an orthogonal rotation $Q.$ The orientation vector for $T$ if $\tilde{x} \in \spt T \cap B_{\rho}(x),$ is given by
$$\ast \vec{T}(\tilde{x}) = Q \left( \left. \left( \frac{-Du}{\sqrt{1+|Du|^{2}}},\frac{1}{\sqrt{1+|Du|^{2}}} \right) \right|_{\proj{\R^{n}}{\tilde{x}}} \right) .$$ \end{theorem}

{\bf Proof:} Suppose (after translation) that $0 \in \spt \partial T,$ and that $T$ at $0$ has a tangent cone which is a hyperplane with constant orientation but non-constant multiplicity. Choose $\rho \in (0,\dist(0,\partial U))$ so that Lemma \ref{hyperplanetangentcone} holds, and so that (after rotation)
$$\partial T \res B_{\rho}(x) = (-1)^{n} \sum_{\ell=1}^{N} m_{\ell} \Phi_{T,\ell \#}(\E^{n-1} \res B^{n-1}_{\rho}(0)) \res B_{\rho}(0)$$
for $N,m_{1},\ldots,m_{N} \in \N$ where $\Phi_{T,\ell}(z) = (z,\varphi_{T,\ell}(z),\psi_{T,\ell}(z))$ with $\varphi_{T,\ell},\psi_{T,\ell} \in C^{1,\alpha}(B^{n-1}_{\rho}(0))$ satisfying $\varphi_{\ell}(0) = \psi_{\ell}(0)=0$ and $D\varphi_{\ell}(0) = D\psi_{\ell}(0)=0$ for each $\ell \in \{ 1,\ldots,N \}.$ We proceed by induction on $N.$

\bigskip

If $N=1,$ then $\spt \partial T \cap B_{\rho}(0) \subset \reg \partial T.$ We conclude the theorem in this case by \cite{HS79} and \cite{W83}, since $T$ at $0$ has tangent cone which is a hyperplane with constant orientation but non-constant multiplicity.

\bigskip

Now suppose $N \geq 2.$ We can of course assume $0 \in \sing \partial T.$ Define the set
$$M = \left( (\spt T \setminus \sing T) \cup (\spt \partial T \setminus \bigcap_{\ell=1}^{N} \{ \Phi_{T,\ell}(z): z \in B^{n-1}_{\rho}(0) \}) \right) \cap B_{\rho}(0).$$
We claim $M$ is a smooth orientable hypersurface, smoothly oriented by $\ast \vec{T}.$ For this, consider $x \in M,$ then there are three possibilities. First, $x \in \reg T \cap B_{\rho}(0),$ in which case $M$ near $x$ is a smooth orientable hypersurface, smoothly oriented by $\vec{T}.$ Second, $x \in \reg \partial T \cap B_{\rho}(0),$ in which case by the first part of Lemma \ref{hyperplanetangentcone} and the boundary regularity given by \cite{HS79},\cite{W83} there is $\sigma>0$ such that $M \cap B_{\sigma}(x)$ is a smooth orientable hypersurface, smoothly oriented by $\ast \vec{T} \in C^{\infty}(M \cap B_{\sigma}(x);\R^{n+1}).$ Third,  
$$x \in \sing \partial T \cap B_{\rho}(0) \setminus \bigcap_{\ell=1}^{N} \{ \Phi_{T,\ell}(z): z \in B^{n-1}_{\rho}(0) \},$$ 
in which case by induction there is again $\sigma>0$ such that $M \cap B_{\sigma}(x)$ is a smooth orientable hypersurface, smoothly oriented by $\ast \vec{T} \in C^{\infty}(M \cap B_{\sigma}(x);\R^{n+1}).$ We conclude $M$ is a smooth orientable hypersurface, smoothly oriented by $\ast \vec{T} \in C^{\infty}(M;\R^{n+1}).$

\bigskip

Consider the current $\lsem M \rsem,$ where $M$ is oriented by $\ast \vec{T}|_{M} \in C^{\infty}(M;\R^{n+1}).$ Note that $\spt \lsem M \rsem \cap B_{\rho}(0) = \spt T \cap B_{\rho}(0)$ (by, for example, Lemma \ref{hyperplanetangentcone} and Theorem \ref{partialboundaryregularity}). We claim $\partial \lsem M \rsem \res B_{\rho}(0) = 0.$ To see this, observe that
$$\spt \partial \lsem M \rsem \cap B_{\rho}(0) \subseteq \bigcap_{\ell=1}^{N} \{ \Phi_{T,\ell}(z): z \in B^{n-1}_{\rho}(0) \} \cap B_{\rho}(0)$$
by interior regularity for area-minimizing currents, \cite{HS79}, \cite{W83}, and induction. By the constancy theorem (see 4.1.31 of \cite{F69}) we conclude that for each $\ell \in \{ 1,\ldots,N \}$
$$\partial \lsem M \rsem \res B_{\rho}(0) = m^{M}_{\ell}  \Phi_{T,\ell \#} (\E^{n-1} \res B^{n-1}_{\rho}(0)) \res B_{\rho}(0)$$
for some integer $m^{M}_{\ell}.$ 

\bigskip

Recall that we have assumed $0 \in \sing \partial T.$ Thus, there is a $z \in B^{n-1}_{\rho}(0)$ so that (after relabeling) $\Phi_{T,1}(z) \in B_{\rho}(0) \setminus \bigcap_{\ell=1}^{N} \Phi_{T,\ell}(B^{n-1}_{\rho}(0)).$ This implies we must have $m^{M}_{1}=0,$ so that $\partial \lsem M \rsem \res B_{\rho}(0)=0$ as claimed.

\bigskip

Now, $\lsem M \rsem$ is area-minimizing by Lemma 33.4 of \cite{S83}. Recall as well that $\spt \lsem M \rsem \cap B_{\rho}(0) = \spt T \cap B_{\rho}(0).$ We conclude (using, for example, Theorem 1 of \cite{SS81} along with the fact that $T$ has tangent cone at $0$ which is a hyperplane with constant orientation but non-constant multiplicity) that $0 \in \reg \lsem M \rsem.$ We conclude the theorem, by induction. $\square$

\bigskip

In the remainder of this section we give two lemmas that we shall need in the subsequent sections. We give these lemmas now, as they regard $T \in \TI^{1,\alpha}_{n,loc}(U).$ First, Lemma \ref{halfregular} shows the existence of \emph{half-regular} points along the boundary of $T \in \TI^{1,\alpha}_{n,loc}$ near any singular point of $\partial T.$ Half-regular points were used in studying the $c$-Plateau problem (see \S 1.4), and in that context appeared in Lemma 1 of \cite{R16}; see as well Lemma \ref{spacehalfregular} forthcoming. We shall make use of half-regular points in proving Theorem \ref{hyperplanetangentconeregularity}. 

\begin{lemma} \label{halfregular} Let $U \subseteq_{o} \R^{n+1},$ $\alpha \in (0,1],$ and $T \in \TI^{1,\alpha}_{n,loc}(U).$ Suppose $x \in \sing \partial T$ and that $\rho \in (0,\dist(x,\partial U))$ is as in Definition \ref{immersedboundary}, so that
$$\partial T \res B_{\rho}(x) = (-1)^{n} \sum_{\ell=1}^{N} m_{\ell} \Big[ (\eta_{-x,1} \circ Q \circ \Phi_{T,\ell})_{\#}(\E^{n-1} \res B^{n-1}_{\rho}(0)) \Big] \res B_{\rho}(x)$$
for $N,m_{1},\ldots,m_{N} \in \N,$ an orthogonal rotation $Q,$ and $\Phi_{T,\ell} \in C^{1,\alpha}(B^{n-1}_{\rho}(0);\R^{n+1})$ for each $\ell = 1,\ldots,N.$ Then there is $\z \in B^{n-1}_{\rho}(0),$ a radius $\sigma \in (0,\rho-|\z|],$ and a non-empty set $\mathcal{N} \subseteq \{ 1,\ldots,N \}$ so that:
\begin{itemize}
 \item $\Phi_{T,\ell}(B^{n-1}_{\sigma}(\z)) \subset \reg \partial T$ for each $\ell \in \mathcal{N},$
 \item $\Phi_{T,\ell}(B^{n-1}_{\sigma}(\z)) \cap \Phi_{T,\tilde{\ell}}(B^{n-1}_{\sigma}(\z)) = \emptyset$ for some $\ell,\tilde{\ell} \in \mathcal{N},$
 \item $\bigcap_{\ell \in \mathcal{N}} \Phi_{T,\ell}(\partial B^{n-1}_{\sigma}(\z)) \neq \emptyset.$
\end{itemize}
With this in mind, we say that any point
$$\x \in \bigcap_{\ell \in \mathcal{N}} (\eta_{-x,1} \circ Q \circ \Phi_{T,\ell})(\partial B^{n-1}_{\sigma}(\z)))$$
is \emph{half-regular}.
\end{lemma}

{\bf Proof:} Suppose (after translation) $0 \in \sing \partial T.$ Also assume $\rho \in (0,\dist(0,\partial U))$ is such that (after rotation)
$$\partial T \res B_{\rho}(0) =  (-1)^{n} \sum_{\ell=1}^{N} m_{\ell} \Phi_{T,\ell \#}(\E^{n-1} \res B^{n-1}_{\rho}(0)) \res B_{\rho}(0).$$
We now argue by induction on $N \geq 2.$

\bigskip

First, suppose $N=2.$ Since $0 \in \sing \partial T,$ then there is $\z \in B^{n-1}_{\rho/2}(0) \setminus \{0\}$ so that $\Phi_{T,1}(\z) \neq \Phi_{T,2}(\z).$ We can thus find $\sigma \in (0,|\z|]$ so that $\Phi_{T,1}(z) \neq \Phi_{T,2}(z)$ for each $z \in B^{n-1}_{\sigma}(\z),$ but there is $\x \in \Phi_{T,1}(\partial B^{n-1}_{\sigma}(\z)) \cap \Phi_{T,2}(\partial B^{n-1}_{\sigma}(\z)).$ This proves the case $N=2.$

\bigskip

Second, suppose $N \geq 3.$ Since $0 \in \sing \partial T$ we can suppose (after relabeling) there is $z^{1} \in B^{n-1}_{\rho/2}(0)$ and $\sigma_{1} \in (0,|z^{1}|]$ so that $\Phi_{T,1}(z) \neq \Phi_{T,2}(z)$ for each $z \in B^{n-1}_{\sigma_{1}}(z^{1})$ while $\Phi_{T,1}(\partial B^{n-1}_{\sigma_{1}}(z^{1})) \cap \Phi_{T,2}(\partial B^{n-1}_{\sigma_{1}}(z^{1})) \neq \emptyset.$ We now consider the two cases: $\Phi_{T,\ell}(z^{1}) \in \reg \partial T$ for each $\ell = 1,\ldots,N$; $\Phi_{T,\ell_{z^{1}}}(z^{1}) \in \sing \partial T$ for some $\ell_{z^{1}} \in \{1,\ldots,N\}.$

\bigskip

\begin{itemize} 
 \item Suppose $\Phi_{T,\ell}(z^{1}) \in \reg \partial T$ for each $\ell = 1,\ldots,N.$ 
 
 \bigskip
 
 Then we can find $\sigma_{2} \in (0,\sigma_{1}]$ such that $\Phi_{T,\ell}(B^{n-1}_{\sigma_{2}}(z^{1})) \subset \reg \partial T$ for each $\ell=1,\ldots,N$ but so that there is $z^{2} \in \partial B^{n-1}_{\sigma_{2}}(z^{1})$ with $\Phi_{T,\ell_{z^{2}}}(z^{2}) \in \sing \partial T$ for some $\ell_{z^{2}} \in \{ 1,\ldots,N \}.$ This leads to two sub-cases:
 \begin{itemize}
  \item Suppose $\Phi_{T,\ell_{z^{2}}}(z^{2}) \in \bigcap_{\ell=1}^{N} \Phi_{T,\ell}( \partial B^{n-1}_{\sigma_{2}}(z^{1})).$
  
   \bigskip
   
   We thus let $\z=z^{1},$ $\sigma = \sigma_{2},$ and $\x = \Phi_{T,1}(z^{2}) = \ldots = \Phi_{T,N}(z^{2}).$ Then we have $\Phi_{T,\ell}(B^{n-1}_{\sigma}(\z)) \subset \reg \partial T$ for each $\ell \in \{1,\ldots,N\}.$ Furthermore, since $\Phi_{T,1}(\z) = \Phi_{T,1}(z_{1}) \neq \Phi_{T,2}(z_{1}) = \Phi_{T,2}(\z),$ then $\Phi_{T,1}(B^{n-1}_{\sigma}(\z)) \cap \Phi_{T,2}(B^{n-1}_{\sigma}(\z)) = \emptyset.$ Since $\x \in \bigcap_{\ell=1}^{N} \Phi_{T,\ell}(\partial B^{n-1}_{\sigma}(\z)),$ then we conclude the lemma in this case with $\mathcal{N} = \{ 1,\ldots,N \}.$
  \item Suppose $\Phi_{T,\ell_{z^{2}}}(z^{2}) \notin \bigcap_{\ell=1}^{N} \Phi_{T,\ell}( \partial B^{n-1}_{\sigma_{2}}(z^{1})).$
  
  \bigskip
  
  Recall that $\Phi_{T,\ell} \in C^{1,\alpha}(B^{n-1}_{\rho}(0);\R^{n+1})$ for each $\ell = 1,\ldots,N.$ Then $\Phi_{T,\ell_{z^{2}}}(z^{2}) \in \sing \partial T$ and $\Phi_{T,\ell_{z^{2}}}(z^{2}) \notin \bigcap_{\ell=1}^{N} \Phi_{T,\ell}( \partial B^{n-1}_{\sigma_{2}}(z^{1}))$ imply that, by induction, we can find $\z \in B^{n-1}_{\rho}(0)$ (in fact, with $\z$ close to $z^{2}$), a radius $\sigma \in (0,\rho-|\z|],$ and a non-empty set $\mathcal{N} \subseteq \{ 1,\ldots,N \}$ with $\ell_{z^{2}} \in \mathcal{N}$ such that there is a half-regular $\x \in \bigcap_{\ell \in \mathcal{N}} \Phi_{T,\ell}(\partial B^{n-1}_{\sigma}(\z))$ as required. 
 \end{itemize}
 \item Suppose $\Phi_{T,\ell_{z^{1}}}(z^{1}) \in \sing \partial T$ for some $\ell_{z^{1}} \in \{ 1,\ldots,N \}.$
 
 \bigskip
 
 Since $\Phi_{T,1}(z^{1}) \neq \Phi_{T,2}(z^{1}),$ then either $\Phi_{T,\ell_{z^{1}}}(z^{1}) \neq \Phi_{T,1}(z^{1})$ or $\Phi_{T,\ell_{z^{1}}}(z^{1}) \neq \Phi_{T,2}(z^{1}).$ Since $\Phi_{T,\ell_{z^{1}}}(z^{1}) \in \sing \partial T,$ then this case also holds by induction.
\end{itemize}

We thus conclude the lemma. $\square$

\bigskip

The next lemma will be convenient for the proof of Theorem \ref{main}.

\begin{lemma} \label{submain} Let $U \subseteq_{o} \R^{n+1},$ $\alpha \in (0,1],$ and suppose $T \in \I^{1,\alpha}_{n,loc}(U).$ Suppose $x \in \spt \partial T$ and there exists $\rho \in (0,\dist(x,\partial U))$ so that
$$\spt T \cap B_{\rho}(x) = \bigcup_{a=1}^{A} (\clos M_{a}) \cap B_{\rho}(x)$$
for $C^{1}$ hypersurfaces-with-boundary $M_{1},\ldots,M_{a}$ in $B_{\rho}(x).$ Then there is $\sigma \in (0,\rho)$ so that  

$$\spt T \cap B_{\sigma}(x) = \bigcup_{ \{ a \in \{ 1,\ldots,A \}: x \in \clos M_{a} \}} (\clos M_{a}) \cap B_{\sigma}(x),$$
   
and $T_{x} \partial T \subset T_{x} M_{a}$ for each $a \in \{ 1,\ldots,A\}$ with $x \in \clos M_{a}.$
\end{lemma}

{\bf Proof:} Choose $\sigma \in (0,\rho)$ so that $B_{\sigma}(x) \cap \clos M_{a} = \emptyset$ for each $a \in \{ 1,\ldots,A\}$ with $x \notin \clos M_{a},$ then $\spt T \cap B_{\sigma}(x) = \bigcup_{ \{ a \in \{ 1,\ldots,A \}: x \in \clos M_{a} \}} (\clos M_{a}) \cap B_{\sigma}(x).$ 

\bigskip

Suppose for contradiction $x \in \clos M_{a}$ but $T_{x} \partial T \not \subset T_{x}M_{a}.$ Let $\C$ be any tangent cone of $T$ at $x,$ then Theorem \ref{monotonicity} implies that either $\C$ is a sum of half-hyperplanes, each containing $T_{x} \partial T,$ with constant orientation after rotation, or $\C$ is a hyperplane, containing $T_{x} \partial T,$ with constant orientation but non-constant multiplicity. Then $T_{x} \partial T \not \subset T_{x} M_{a}$ implies we can find $\tilde{x} \in T_{x} M_{a} \setminus \{0\}$ so that $\{ t \tilde{x}: t \in \R \} \cap \spt \C = \{0\}.$ However, the fact that $M_{a}$ is a $C^{1}$ hypersurface-with-boundary with $0 \in \clos M_{a}$ and Theorem 5.4.2 of \cite{F69} imply that either $\{ t \tilde{x}: t \geq 0 \} \subset \spt \C,$ $\{ t \tilde{x}: t \leq 0 \} \subset \spt \C,$ or $\{ t \tilde{x}: t \in \R \} \subset \spt \C$ (the last occurs if, for example, if $x \in M_{a}$), giving a contradiction. $\square$

\section{Boundaries With Co-Oriented Mean Curvature}

In this section we define what it means for a co-dimension one area-minimizing locally rectifiable current $T$ to have boundary with co-oriented mean curvature. We also give two basic results that we will need in \S 6, but are also of independent interest; these results will follow straight from well-known theory.

\begin{definition} \label{cmcboundary} Let $U \subseteq_{o} \R^{n+1},$ and suppose $T \in \I_{n,loc}(U)$ is area-minimizing. Then we say $\partial T$ has \emph{co-oriented mean curvature} if $\partial T$ has mean curvature $H_{\partial T} = h \nu_{T}$ for $h:U \rightarrow \R$ a $\mu_{\partial T}$-locally integrable function, where $\nu_{T}:U \rightarrow \R^{n+1}$ is the generalized outward pointing normal of $\partial T$ with respect to $T$; this means that
$$\int \dive{\partial T}{X} \ d \mu_{\partial T} = \int X \cdot (h \nu_{T}) \ d \mu_{\partial T}$$
for all $X \in C^{1}_{c}(U;\R^{n+1}).$ \end{definition}

The assumption that $T$ is area-minimizing in Definition \ref{cmcboundary} is merely to guarantee the existence of the generalized outward pointing unit normal $\nu_{T}$ of $\partial T$ with respect to $T$; see Lemma 3.1 of \cite{B77} and (2.10) of \cite{E89}. 

\bigskip

We now give the following boundary regularity result, which follows directly from \cite{HS79},\cite{W83}.

\begin{theorem} \label{boundaryregularity} Let $U \subseteq_{o} \R^{n+1}$ and suppose $T \in \I_{n,loc}(U)$ is area-minimizing, and that $\partial T$ has locally bounded co-oriented mean curvature; meaning $H_{\partial T} = h \nu_{T}$ with $h \in L^{\infty}_{loc}(\mu_{\partial T}).$ If $x \in \reg \partial T,$ then there is $\rho>0$ such that $(\spt \partial T) \cap B_{\rho}(0)$ is an $(n-1)$-dimensional $C^{1,\alpha}$ submanifold for any $\alpha \in (0,1).$ Furthermore, one of the following holds:
\begin{enumerate}
 \item[(1)] $T$ at $x$ has unique tangent cone which is a sum of half-hyperplanes with constant orientation after rotation; see Definition \ref{conesdefinition}.
 
 \bigskip

Moreover, there are disjoint orientable $C^{1,\alpha}$ hypersurfaces-with-boundary $M_{1},\ldots,M_{N}$ in $B_{\rho}(x),$ with $(\partial M_{\ell}) \cap B_{\rho}(x) = (\spt \partial T) \cap B_{\rho}(x)$ and $m_{1},\ldots,m_{N} \in \N$ so that
$$T \res B_{\rho}(x) = \sum_{\ell=1}^{N} m_{\ell} \lsem M_{\ell} \rsem;$$
each $M_{\ell}$ is oriented so that $\partial \lsem M_{\ell} \rsem \res B_{\rho}(x) = \lsem \spt \partial T \rsem \res B_{\rho}(x)$ for each $\ell=1,\ldots,N,$ where $\lsem \spt \partial T \rsem \res B_{\rho}(x)$ has orientation $\vec{\partial T}.$ Furthermore, no two $M_{\ell}$ can meet tangentially at any point of $(\spt \partial T) \cap B_{\rho}(x).$ Thus, for any $\tilde{x} \in (\spt \partial T) \cap B_{\rho}(x)$ we have
$$\nu_{T}(\tilde{x}) = \frac{\sum_{\ell=1}^{N} m_{\ell} \nu_{M_{\ell}}(\tilde{x})}{\sum_{\ell=1}^{N} m_{\ell}}$$
where $\nu_{M_{\ell}}$ is the outward pointing unit normal of $(\spt \partial T) \cap B_{\rho}(x)$ with respect to $M_{\ell}$ for each $\ell = 1,\ldots,N.$

\bigskip

If in addition $h|_{(\spt \partial T) \cap B_{\rho}(x)} \in C^{k,\alpha}((\spt \partial T) \cap B_{\rho}(x),\R^{n+1})$ with $k \in \{0\} \cup \N$ (or $C^{\infty},$ analytic), then $(\spt \partial T) \cap B_{\rho}(x)$ is $C^{2+k,\alpha}$ (respectively $C^{\infty},$ analytic) and each $M_{\ell}$ is $C^{2+k,\alpha}$ (respectively $C^{\infty},$ analytic).

 \item[(2)] $T$ at $x$ has unique tangent cone which is a hyperplane with constant orientation but non-constant multiplicity; see Definition \ref{conesdefinition}.
 
 \bigskip
 
 Furthermore, there is an orientable analytic minimal hypersurface $M$ in $B_{\rho}(x)$ containing $(\spt \partial T) \cap B_{\rho}(x)$ and $m,\theta \in \N$ so that
$$\begin{aligned}
T \res B_{\rho}(x) = & (m+\theta) \lsem M^{+} \rsem + \theta \lsem M^{-} \rsem \\
= & m \lsem M^{+} \rsem + \theta \lsem M \rsem,
\end{aligned}$$
where we have the following: each $M^{\pm} \subset M$ is an orientable $C^{1,\alpha}$ hypersurface-with-boundary in $B_{\rho}(x),$ with $(\partial M^{\pm}) \cap B_{\rho}(x) = (\spt \partial T) \cap B_{\rho}(x)$; each of $\lsem M \rsem,\lsem M^{\pm} \rsem$ is oriented by $\vec{T}$; $M \cap B_{\rho}(x) = (M^{+} \cup \spt \partial T \cup M^{-}) \cap B_{\rho}(x)$; $M^{\pm},\spt \partial T$ are pairwise disjoint in $B_{\rho}(x).$ In this case, $\nu_{T}(\tilde{x})$ for $\tilde{x} \in (\spt \partial T) \cap B_{\rho}(x)$ is the outward pointing unit normal of $\spt \partial T$ with respect to $M^{+}.$

\bigskip
 
If in addition $h|_{(\spt \partial T) \cap B_{\rho}(x)} \in C^{k,\alpha}((\spt \partial T) \cap B_{\rho}(x),\R^{n+1})$ for $k \in \{0\} \cup \N$ (or $C^{\infty},$ analytic), then $(\spt \partial T) \cap B_{\rho}(x)$ is $C^{2+k,\alpha}$ (respectively $C^{\infty},$ analytic), and each $M^{\pm}$ is a $C^{2+k,\alpha}$ (respectively $C^{\infty},$ analytic) hypersurface-with-boundary.
\end{enumerate}
\end{theorem}

{\bf Proof:} Since $\partial T$ has locally bounded mean curvature, then standard regularity theory for $C^{1}$ solutions to the mean curvature system (see for example \S 6.8 of \cite{M66}) implies that for $x \in \reg \partial T$ there is $\rho \in (0,\dist(x,\partial U))$ so that $(\spt \partial T) \cap B_{\rho}(x)$ is an $(n-1)$-dimensional $C^{1,\alpha}$ submanifold for any $\alpha \in (0,1).$ If in addition $h|_{(\spt \partial T) \cap B_{\rho}(x)} \in C^{k,\alpha}((\spt \partial T) \cap B_{\rho}(x),\R^{n+1})$ for $k \in \{0\} \cup \N$ (or $C^{\infty},$ analytic), then $(\spt \partial T) \cap B_{\rho}(x)$ is $C^{k+2}$ (respectively $C^{\infty},$ analytic; see again \S 6.8 of \cite{M66}). 
\bigskip

The remainder of the theorem then follows from the boundary regularity theory for co-dimension one area-minimizing currents, given by \cite{HS79},\cite{W83}. $\square$

\bigskip

We end this section by giving the following monotonicity formula.

\begin{theorem} \label{boundarymonotonicity} Let $U \subseteq_{o} \R^{n+1}$ and suppose $T \in \I_{n,loc}(U)$ is area-minimizing. Also suppose $\partial T$ has co-oriented mean curvature $H_{\partial T} = h \nu_{T}.$ If $x \in U,$ $R \in (0,\dist(x,\partial U)),$ and for some $\alpha \in (0,1]$ and $\Lambda \geq 0$ we have
$$\frac{1}{\alpha} \int_{B_{\rho}(x)} |h| \ d \mu_{\partial T} \leq \Lambda (\rho/R)^{\alpha-1} \mu_{\partial T}(B_{\rho}(x)) \text{ for all } \rho \in (0,R),$$
then
$$\begin{aligned}
e^{\Lambda R^{1-\alpha} \rho^{\alpha}} \frac{\mu_{\partial T}(B_{\rho}(x))}{\rho^{n-1}} - & e^{\Lambda R^{1-\alpha} \sigma^{\alpha}} \frac{\mu_{\partial T}(B_{\sigma}(x))}{\sigma^{n-1}} \\
& \geq \int_{B_{\rho}(x) \setminus B_{\sigma}(x)} \frac{|\proj{T^{\perp}_{\tilde{x}}\partial T}{(\tilde{x}-x)}|^{2}}{|\tilde{x}-x|^{n+1}} \ d \mu_{\partial T}(\tilde{x})
\end{aligned}$$
whenever $0 < \sigma < \rho \leq R.$
\end{theorem}

{\bf Proof:} This is the usual monotonicity formula applied to $\partial T$; see for example Theorem 17.6 of \cite{S83}. $\square$

\section{Tangentially Immersed Boundaries with Co-Oriented Mean Curvature.}

We now study co-dimension one area-minimizing currents with boundary being both $C^{1,1}$ tangentially immersed and having co-oriented Lipschitz mean curvature. The first main result we wish to show is that the boundary $\partial T$ of any such current $T$ is regular near any point $x \in \spt \partial T$ such that $T$ at $x$ has tangent cone which is a hyperplane with constant orientation but non-constant multiplicity; this is Theorem \ref{hyperplanetangentconeregularity}. For this, we must prove Lemma \ref{halfregulartangentcones}, which shows that if $x \in \spt \partial T$ is half-regular (see Lemma \ref{halfregular}), then every tangent cone of $T$ at $x$ must be a sum of half-hyperplanes with constant orientation after rotation. Theorem \ref{hyperplanetangentconeregularity} then follows from Theorem \ref{hyperplanetangentconegraph}. We then give Theorem \ref{main}, which states that near any $x \in \sing \partial T,$ either $T$ near $x$ exhibits a reasonable amount of regularity or $\spt T$ near $x$ must be extremely irregular.

\bigskip

As mentioned in \S 1.3, we wish in the future to prove the results of this section for $T \in \TI^{1,\alpha}_{n,loc}(U)$ with boundary having co-oriented Lipschitz mean curvature, but with $\alpha \in (0,1]$ more generally. 

\begin{lemma} \label{halfregulartangentcones} Let $U \subseteq_{o} \R^{n+1},$ and suppose $T \in \TI^{1,1}_{n,loc}(U)$ where $\partial T$ has co-oriented mean curvature $H_{\partial T} = h \nu_{T}$ with $h:U \rightarrow \R$ Lipschitz. For any $x \in \sing \partial T,$ there is $\rho \in (0,\dist(x,\partial U))$ so that for any half-regular $\x \in \sing \partial T \cap B_{\rho}(0)$ (see Lemma \ref{halfregular}), every tangent cone of $T$ at $\x$ is the sum of half-hyperplanes with constant orientation after rotation (see Definition \ref{conesdefinition}).
\end{lemma}

{\bf Proof:} Suppose (after translation) $0 \in \sing \partial T,$ and choose $\rho \in (0,\dist(0,\partial U))$ so that (after rotation)
$$\partial T \res B_{\rho}(0) = (-1)^{n} \sum_{\ell=1}^{N} m_{\ell} \Phi_{T,\ell \#}(\E^{n-1} \res B^{n-1}_{\rho}(0)) \res B_{\rho}(0)$$
for $N,m_{1},\ldots,m_{N} \in \N$ and $\Phi_{T,\ell} \in C^{1,1}(B^{n-1}_{\rho}(0);\R^{n+1})$ is the map $\Phi_{T,\ell}(z)=(z,\varphi_{T,\ell}(z),\psi_{T,\ell}(z))$ where $\varphi_{T,\ell},\psi_{T,\ell} \in C^{1,1}(B^{n-1}_{\rho}(0))$ satisfy $\varphi_{T,\ell}(0)=\psi_{T,\ell}(0)=0$ and $D \varphi_{T,\ell}(0)=D \psi_{T,\ell}(0)=0$ for each $\ell=1,\ldots,N.$ We also choose $\rho \in (0,\dist(x,\partial U))$ sufficiently small depending on $\epsilon=\epsilon(n,\partial T)>0,$ to be chosen later, so that 
\begin{equation} \label{hrtc1}
\| D \varphi_{T,\ell} \|_{C(B^{n-1}_{\rho}(0))},\| D \psi_{T,\ell} \|_{C(B^{n-1}_{\rho}(0))} < \epsilon.\end{equation}

\bigskip

Suppose $\x \in \spt \partial T \cap B_{\rho/3}(0)$ is a half-regular point. Thus, by Lemma \ref{halfregular} there is $\z \in B^{n-1}_{\rho/3}(0)$ and $\sigma \in (0,\rho/3-|\z|]$ is such that (after relabeling)
\begin{equation} \label{hrtc2}
\begin{aligned}
\x \in & \Phi_{T,1}(\partial B^{n-1}_{\sigma}(\z)) \cap \Phi_{T,2}(\partial B^{n-1}_{\sigma}(\z)), \\
& \Phi_{T,1}(B^{n-1}_{\sigma}(\z)),\Phi_{T,2}(B^{n-1}_{\sigma}(\z)) \subset \reg \partial T, \text{ and} \\ 
& \Phi_{T,1}(B^{n-1}_{\sigma}(\z)) \cap \Phi_{T,2}(B^{n-1}_{\sigma}(\z)) = \emptyset.
\end{aligned}
\end{equation}

\bigskip

Suppose for contradiction that $T$ at $\x$ has tangent cone $\C$ which is a hyperplane with constant orientation but non-constant multiplicity. In fact, \eqref{hrtc1} implies 
\begin{equation} \label{hrtc3}
\C = Q^{T_{\x} \partial T}_{\#} Q_{\#} \big( (m+\theta) \E^{n} \res \{ y \in \R^{n} : y_{n}>0 \} + \theta \E^{n} \res \{ y \in \R^{n} : y_{n}<0 \} \big)
\end{equation}
for some $m,\theta \in N,$ a rotation $Q$ about $\R^{n-1},$ and $Q^{T_{\x} \partial T}$ an orthogonal rotation with $\|Q^{T_{\x} \partial T}-I\| < c \epsilon$ for some $c=c(n)>0$ (in fact, $Q^{T_{\x} \partial T}(\R^{n-1}) = T_{\x} \partial T$).

\bigskip

Taking $\z$ closer to $\x$ in that direction, we can also assume $\z \in B^{n-1}_{\rho/3}(0)$ and $\sigma \in (0,\rho/3-|\z|]$ are such that \eqref{hrtc2} continues to hold while
\begin{equation} \label{hrtc4}
\begin{aligned}
& \spt ((Q^{-1} \circ \eta_{\x,1})_{\#} T) \cap B_{3\sigma}(0) = (\gph{B^{n}_{3\sigma}(0)}{u}) \cap B_{3 \sigma}(0) \text{ and} \\
& \ast ((Q^{-1} \circ \eta_{\x,1})_{\#}T) (y,u(y)) = \frac{(-Du(y),1)}{\sqrt{1+|Du(y)|^{2}}} \text{ for } (y,u(y)) \in B_{3 \sigma}(0)
\end{aligned}
\end{equation}
where $u \in C^{\infty}(B^{n}_{3 \sigma}(0))$ with $u(0)=0$; this follows by Theorem \ref{hyperplanetangentconegraph} and \eqref{hrtc3}, in particular using $\|Q^{T_{\x} \partial T}-I\| < c \epsilon$ for some $c=c(n)>0$ where we choose $\epsilon=\epsilon(n)>0$ sufficiently small in \eqref{hrtc1}. Note that while we may have $Du(0) \neq 0,$ we do have by \eqref{hrtc3} that $|Du(0)| < c \epsilon$ for some $c=c(n)>0.$

\bigskip

For $\ell=1,2$ define for $z \in B^{n-1}_{3 \sigma}(0)$
\begin{equation}
\label{hrtc5}
s^{(\ell)}(z) = (Q^{-1} \circ \eta_{\x,1} \circ \Phi_{T,\ell})(z+\proj{\R^{n-1}}{\x}) \cdot e_{n}.
\end{equation}
Observe that $s^{(1)},s^{(2)} \in C^{1,1}(B^{n-1}_{\rho}(0)),$ as $\Phi_{T,1},\Phi_{T,2} \in C^{1,1}(B^{n-1}_{\rho}(0);\R^{n+1}).$ Since $Q$ is a rotation about $\R^{n-1},$ then by \eqref{hrtc4} (if $\epsilon=\epsilon(n)>0$ is sufficiently small in \eqref{hrtc1}) we have for each $\ell=1,2$
$$\begin{aligned}
(Q^{-1} \circ \eta_{\x,1} \circ \Phi_{T,\ell}) & (B^{n-1}_{2 \sigma}(\proj{\R^{n-1}}{\x})) \\
& = \{ (z,s^{(\ell)}(z),u(z,s^{(\ell)}(z))): z \in B^{n-1}_{2\sigma}(0) \} \subset B_{3 \sigma}(0).
\end{aligned}$$

\bigskip

We as well define for each $\ell=1,2$ and $z \in B^{n-1}_{2 \sigma}(0)$
\begin{equation} \label{hrtc6}
\begin{aligned}
h^{(\ell)}(z) & = h \Big((\eta_{-\x,1} \circ Q) \big( z,s^{(\ell)}(z),u(z,s^{(\ell)}(z)) \big) \Big), \\
Du^{(\ell)}(z) & = (Du)(z,s^{(\ell)}(z)), \\
u^{(\ell)}_{i}(z) & = (D_{i}u)(z,s^{(\ell)}(z)) \text{ for } i \in \{ 1,\ldots,n \}, \\
D^{2}u^{(\ell)}(z) & = (D^{2}u)(z,s^{(\ell)}(z)), \\ 
s^{(\ell)}_{j}(z) & = \D_{j} s^{(\ell)}(z) = \D s^{(\ell)}(z) \cdot e_{j} \text{ and } \\
\partial^{(\ell)}_{j}(z) & = (e_{j},s^{(\ell)}_{j}(z),u^{(\ell)}_{j}(z)+u^{(\ell)}_{n}(z) s^{(\ell)}_{j}(z)) \text{ for } j \in \{ 1,\ldots,n-1 \},
\end{aligned}
\end{equation}
where for emphasis we let $\D$ be the gradient over $\R^{n-1}.$ Also consider for $\ell=1,2$ the downward pointing unit normal of the graph of $s^{(\ell)}$ within the graph of $u$; this is given by
\begin{equation} \label{hrtc7}
\nu^{(\ell)} = \frac{ \left( ( \D s^{(\ell)},-1,0 ) + \left( \frac{ ( \D s^{(\ell)},-1 ) \cdot Du^{(\ell)}}{1+|Du^{(\ell)}|^{2}} \right) ( -Du^{(\ell)},1 ) \right)}{\sqrt{1+|\D s^{(\ell)}|^{2} - \frac{( ( \D s^{(\ell)},-1 ) \cdot Du^{(\ell)})^{2}}{1+|Du^{(\ell)}|^{2}}}}.
\end{equation}

\bigskip

Let $g \in C^{\infty}(\R^{n-1} \times \R^{n})$ be the $(n-1) \times (n-1)$ matrix
\begin{equation} \label{hrtc8}
\begin{aligned}
g(p,q) = I & +(1+q_{n}^{2}) pp^{T}+q_{n} \big( (\proj{\R^{n-1}}{q}) p^{T}+ p(\proj{\R^{n-1}} {q})^{T} \big) \\
& + (\proj{\R^{n-1}}{q})(\proj{\R^{n-1}}{q})^{T}
\end{aligned}
\end{equation}
for $p \in \R^{n-1}$ and $q \in \R^{n},$ where $I$ is the $(n-1) \times (n-1)$ identity matrix. Also let $g^{ij} \in C^{\infty}(\R^{n-1} \times \R^{n})$ for $i,j \in \{ 1,\ldots,n-1 \}$ be so that $g^{ij}(p,q)$ is the $ij$-th entry of $g(p,q)^{-1}$ for $p \in \R^{n-1}$ and $q \in \R^{n}$ ($g$ is generally invertible, but by \eqref{hrtc1},\eqref{hrtc3},\eqref{hrtc4},\eqref{hrtc5} we may if we like restrict $g$ to $|p|,|q|$ small). 

\bigskip

For $\ell = 1,2$ let $S^{(\ell)} = \{ (z,s^{(\ell)}(z),u(z,s^{(\ell)}(z))): z \in B^{n-1}_{2 \sigma}(0) \}.$ Since $H_{\partial T} = h \nu_{T},$ then \eqref{hrtc1},\eqref{hrtc2},\eqref{hrtc3},\eqref{hrtc4},\eqref{hrtc5} imply that for $\ell = 1,2$
$$\sum_{i,j=1}^{n-1} g^{ij}(\D s^{(\ell)},Du^{(\ell)}) \proj{T^{\perp}_{(z,s^{(\ell)}(z),u(z,s^{\ell}(z)))} S^{(\ell)}}{D_{\partial^{(\ell)}_{i}} \partial^{(\ell)}_{j}} = - h^{(\ell)} \nu^{(\ell)}$$
for $z \in B^{n-1}_{\sigma}(\z-\proj{\R^{n-1}}{\x}),$ where we have used the notation \eqref{hrtc6},\eqref{hrtc7},\eqref{hrtc8}. This implies, using that $v^{(\ell)} \perp \partial^{(\ell)}_{j}$ for each $j \in \{ 1,\ldots,n-1 \},$ 
\begin{equation}
\label{hrtc9}
\sum_{i,j=1}^{n-1} g^{ij}(\D s^{(\ell)},Du^{(\ell)}) \left( \partial^{(\ell)}_{j} \cdot D_{\partial^{(\ell)}_{i}} \nu^{(\ell)} \right) = h^{(\ell)}.
\end{equation}

\bigskip

Next, we compute by \eqref{hrtc6},\eqref{hrtc7}
$$\partial^{(\ell)}_{j} \cdot D_{\partial^{(\ell)}_{i}} \nu^{(\ell)} = \partial^{(\ell)}_{j} \cdot \frac{\partial}{\partial z_{i}} \left( \frac{ \left( ( \D s^{(\ell)},-1,0 ) + \left( \frac{ ( \D s^{(\ell)},-1 ) \cdot Du^{(\ell)}}{1+|Du^{(\ell)}|^{2}} \right) ( -Du^{(\ell)},1 ) \right)}{\sqrt{1+|\D s^{(\ell)}|^{2} - \frac{( ( \D s^{(\ell)},-1 ) \cdot Du^{(\ell)})^{2}}{1+|Du^{(\ell)}|^{2}}}} \right).$$
Since $\partial^{(\ell)}_{j} \perp ( \D s^{(\ell)},-1,0 ),( -Du^{(\ell)},1 )$
$$\partial^{(\ell)}_{j} \cdot D_{\partial^{(\ell)}_{i}} \nu^{(\ell)} = \frac{\partial^{(\ell)}_{j} \cdot \left( \left( ( \D s^{(\ell)}_{i},0,0 ) + \left( \frac{ ( \D s^{(\ell)},-1 ) \cdot Du^{(\ell)}}{1+|Du^{(\ell)}|^{2}} \right) ( -D^{2}u^{(\ell)}(e_{i},s^{(\ell)}_{i}),0 ) \right) \right)}{\sqrt{1+|\D s^{(\ell)}|^{2} - \frac{( ( \D s^{(\ell)},-1 ) \cdot Du^{(\ell)})^{2}}{1+|Du^{(\ell)}|^{2}}}},$$
using the notation from \eqref{hrtc6}. We simplify this to get
$$\partial^{(\ell)}_{j} \cdot D_{\partial^{(\ell)}_{i}} \nu^{(\ell)} = \frac{s^{(\ell)}_{ij}-\left( \frac{ ( \D s^{(\ell)},-1 ) \cdot Du^{(\ell)}}{1+|Du^{(\ell)}|^{2}} \right) \big( (e_{j},s^{(\ell)}_{j}) D^{2}u^{(\ell)}(e_{i},s^{(\ell)}_{i}) \big)}{\sqrt{1+|\D s^{(\ell)}|^{2} - \frac{( ( \D s^{(\ell)},-1 ) \cdot Du^{(\ell)})^{2}}{1+|Du^{(\ell)}|^{2}}}}.$$
Using this in \eqref{hrtc9}, we conclude for $\ell=1,2$ and $z \in B^{n-1}_{\sigma}(\z-\proj{\R^{n-1}}{\x})$
\begin{equation}
\label{hrtc10}
\sum_{i,j=1}^{n-1} a_{ij}(z,s^{(\ell)},\D s^{(\ell)}) s^{(\ell)}_{ij}+b(z,s^{(\ell)},\D s^{(\ell)})=0
\end{equation}
where we define the functions $a_{ij} \in C^{\infty}(B^{n}_{2 \sigma}(0) \times \R^{n-1})$ and the Lipschitz function $b:B^{n}_{2 \sigma}(0) \times \R^{n-1} \rightarrow \R$ given by
$$\begin{aligned}
a_{ij}(y,p) = & g^{ij}(p,Du(y)) \\
b(y,p) = & - \sum_{i,j=1}^{n-1} g^{ij}(p,Du(y)) \left( \frac{ ( p,-1 ) \cdot Du(y)}{1+|Du(y)|^{2}} \right) \big( (e_{j},p_{j}) D^{2}u(y) (e_{i},p_{i}) \big) \\
& - h((\eta_{-\x,1} \circ Q)(y,u(y))) \sqrt{1+|p|^{2} - \frac{( ( p,-1 ) \cdot Du(y))^{2}}{1+|Du(y)|^{2}}}
\end{aligned}$$
for $y \in B^{n}_{2 \sigma}(0)$ and $p \in \R^{n-1}.$ 

\bigskip

Subtracting \eqref{hrtc10} across $\ell=1,2$ we get for $z \in B^{n-1}_{\sigma}(\z-\proj{\R^{n-1}}{\x})$
$$\begin{aligned}
\sum_{i,j=1}^{n-1} & \int^{1}_{0} \frac{d}{dt} \left[ a_{ij}(z,t s^{(2)}+(1-t) s^{(1)},t \D s^{(2)}+(1-t) \D s^{(1)})(t s^{(2)}_{ij}+(1-t) s^{(1)}_{ij}) \right] \ dt \\
& + \int^{1}_{0} \frac{d}{dt} b(z,t s^{(2)}+ (1-t) s^{(1)}, t \D s^{(2)} + (1-t) \D s^{(1)})=0,
\end{aligned}$$
using that $b$ is Lipschitz. This gives that $s = s^{(2)}-s^{(1)}$ is a solution to the equation over $B^{n-1}_{\sigma}(\z-\proj{\R^{n-1}}{\x})$ given by
$$\sum_{i,j=1}^{n-1} \mathcal{A}_{ij} \D_{i} \D_{j} s + \mathcal{B} \cdot \D s + \mathcal{C} s = 0$$
where for $z  \in B^{n-1}_{2 \sigma}(0)$
$$\begin{aligned}
\mathcal{A}_{ij}(z) = & \sum_{i,j=1}^{n-1} \int^{1}_{0} a_{ij}(z,t s^{(2)}+(1-t) s^{(1)},t \D s^{(2)}+(1-t) \D s^{(1)}) \ dt \\
\mathcal{B}(z) = & \sum_{i,j=1}^{n-1} \int^{1}_{0} \frac{\partial a_{ij}}{\partial p} \Big|_{(z,t s^{(2)}+(1-t) s^{(1)},t \D s^{(2)}+(1-t) \D s^{(1)})} (ts^{(2)}_{ij}+(1-t)s^{(1)}_{ij}) \ dt \\
& + \int^{1}_{0} \frac{\partial b}{\partial p}|_{(z,t s^{(2)}+ (1-t) s^{(1)}, t \D s^{(2)} + (1-t) \D s^{(1)})} \ dt \\
\mathcal{C}(z) = & \sum_{i,j=1}^{n-1} \int^{1}_{0} \frac{\partial a_{ij}}{\partial s} \Big|_{(z,t s^{(2)}+(1-t) s^{(1)},t \D s^{(2)}+(1-t) \D s^{(1)})} (ts^{(2)}_{ij}+(1-t)s^{(1)}_{ij}) \ dt \\
& + \int^{1}_{0} \frac{\partial b}{\partial s}|_{(z,t s^{(2)}+ (1-t) s^{(1)}, t \D s^{(2)} + (1-t) \D s^{(1)})} \ dt.
\end{aligned}$$

\bigskip

By definition (see \eqref{hrtc8} as well as the definitions after \eqref{hrtc10}), the coefficients $\mathcal{A}_{ij}$ are uniformly elliptic (again, we may use \eqref{hrtc2} with $\epsilon = \epsilon(n,\partial T)$ sufficiently small). Since $u \in C^{\infty}(B^{n}_{3 \sigma}(0)),$ the functions $s^{(2)},s^{(1)} \in C^{1,1}(B^{n-1}_{2 \sigma}(0)),$ and $h$ is Lipschitz, then the coefficients $\mathcal{B},\mathcal{C}$ are bounded in $B^{n-1}_{2 \sigma}(0).$ On the other hand, $T \in \TI^{1,1}_{n,loc}(U)$ and \eqref{hrtc5} give $s(0)=0$ and $\D s(0)=0.$ But by \eqref{hrtc2} (after relabeling, if necessary) $s^{(1)}(z) < s^{(2)}(z)$ so that $s(z)>0$ for each $z \in B^{n-1}_{\sigma}(\z-\proj{\R^{n-1}}{\x}).$ This contradicts the Hopf boundary point lemma (see for example Lemma 3.4 of \cite{GT83}). $\square$

\bigskip

Having shown Lemma \ref{halfregulartangentcones}, then the proof of the first main result of this section is relatively short.

\begin{theorem} \label{hyperplanetangentconeregularity} Let $U \subseteq_{o} \R^{n+1}$ and suppose $T \in \TI^{1,1}_{n,loc}(U)$ where $\partial T$ has co-oriented mean curvature $H_{\partial T} = h \nu_{T}$ with $h:U \rightarrow \R$ Lipschitz. If $x \in \spt \partial T$ and $T$ at $x$ has tangent cone which is a hyperplane with constant orientation but non-constant multiplicity (as in Definition \ref{conesdefinition}), then $x \in \reg \partial T.$ \end{theorem}

{\bf Proof:} Suppose for contradiction $x \in \sing \partial T$ and that $T$ at $x$ has tangent cone which is a hyperplane with constant orientation but non-constant multiplicity. By Theorem \ref{hyperplanetangentconegraph} there is $\rho \in (0,\dist(x,\partial U))$ so that $T$ at every $\tilde{x} \in B_{\rho}(x) \cap \spt \partial T$ has unique tangent cone which is a hyperplane with constant orientation but non-constant multiplicity. However, by Lemma \ref{halfregular} we can find a half-regular $\x \in B_{\rho}(x) \cap \sing \partial T.$ This contradicts Lemma \ref{halfregulartangentcones}. $\square$

\bigskip

We are not ready to prove the final main result of this section. 

\begin{theorem} \label{main} Let $U \subseteq_{o} \R^{n+1}$ and suppose $T \in \TI^{1,1}_{n,loc}(U)$ where $\partial T$ has co-oriented mean curvature $H_{\partial T} = h \nu_{T}$ with $h:U \rightarrow (0,\infty)$ Lipschitz. Suppose $x \in \spt \partial T$ and that there exists $\rho \in (0,\dist(x,\partial U))$ and $C^{1}$ hypersurfaces-with-boundary $M_{1},\ldots,M_{A}$ in $B_{\rho}(x)$ so that either:
\begin{enumerate}
 \item[(1)] $\spt T \cap B_{\rho}(x) = \bigcup_{a=1}^{A} (\clos M_{a}) \cap B_{\rho}(x),$ or
 \item[(2)] $\spt T \cap B_{\rho}(x) \subseteq \bigcup_{a=1}^{A} (\clos M_{a}) \cap B_{\rho}(x)$ and $T^{\perp}_{x} \partial T \not \subset T_{x} M_{a}$ for each $a \in \{ 1,\ldots,A \}$ such that $x \in \clos M_{a}.$
\end{enumerate}
Then there is $\sigma \in (0,\rho)$ and $\mathcal{B} \in \{ 1,\ldots,2 \Theta_{T}(x) \}$ so that
$$\spt T \cap B_{\sigma}(x) = \bigcup_{b=1}^{\mathcal{B}} (\clos W_{b}) \cap B_{\sigma}(x)$$
for orientable $C^{1,1}$ hypersurfaces-with-boundary $W_{1},\ldots,W_{\mathcal{B}}$ in $B_{\sigma}(x).$ For each $b \in \{ 1,\ldots,\mathcal{B} \}$ we have $x \in \partial W_{b}$ and $W_{b} \cap \spt \partial T \subset \reg \partial T.$ Furthermore, for each $b,\tilde{b} \in \{ 1,\ldots,\mathcal{B} \}$ we have 
$$(\clos W_{b}) \cap (\clos W_{\tilde{b}}) \cap B_{\sigma}(x) \subseteq (\partial W_{b}) \cap (\partial W_{\tilde{b}}) \cap B_{\sigma}(x).$$
\end{theorem}

Theorem \ref{main} essentially appears as Theorem 9 of \cite{R16} in the context of two-dimensional solutions to the $c$-Plateau problem; see \S 1.4. The proof of Theorem \ref{main} follows closely the proof of Theorem 9 of \cite{R16}, with both some simplifications (primarly due to Theorem \ref{hyperplanetangentconegraph}) to the argument and some more subtle analysis due to the more general present setting.

\bigskip

{\bf Proof:} Observe that by Lemma \ref{submain}, it suffices to show (2). This we do in what follows.

\bigskip

Suppose (after translation) that $0 \in \spt \partial T$ and (after rotation) $T_{0} \partial T = \R^{n-1}.$ Also suppose (by scaling) that $U=B_{1}(0)$ and
$$\spt T \cap B_{1}(0) \subseteq \bigcup_{a=1}^{A} (\clos M_{a}) \cap B_{1}(0)$$ 
with $0 \in \clos M_{a}$ for each $a = 1,\ldots,A.$ Observe that for any $\epsilon>0$ there is $\rho=\rho(\epsilon,\{M_{a}\}_{a=1}^{A}) \in (0,1)$ so that for every $a=1,\ldots,A$
\begin{equation}
\label{main1}
\dist_{\HH}(T_{x} M_{a},T_{0} M_{a}) < \epsilon \text{ for each } x \in (\clos M_{a}) \cap B_{\rho}(0).
\end{equation}
Recall that we are assuming (1) holds, so that $\{0\} \times \R^{2} = T^{\perp}_{0} \partial T \not \subset T_{0} M_{a}$ for each $a=1,\ldots,A.$ Hence, if $\epsilon=\epsilon(n,\{T_{0}M_{a}\}_{a=1}^{A})>0$ is sufficiently small, then \eqref{main1} with $\rho=\rho(\epsilon,\{M_{a}\}_{a=1}^{A}) \in (0,1)$ implies $\{0\} \times \R^{2} \not \subset T_{x} M_{a}$ for each $x \in (\clos M_{a}) \cap B_{\rho}(0)$ and $a = 1,\ldots,A.$

\bigskip

By definition of $T \in \TI^{1,1}_{n,loc}(B_{1}(0)),$ for any $\epsilon>0$ there is $\rho = \rho(\epsilon,\partial T) \in (0,1)$ so that
\begin{equation} \label{main2}
\partial T \res B_{\rho}(0) = (-1)^{n} \sum_{\ell=1}^{N} m_{\ell} \Phi_{T,\ell \#}(\E^{n-1} \res B^{n-1}_{\rho}(0)) \res B_{\rho}(0)
\end{equation}
where $\Phi_{T,\ell}(z) = (z,\varphi_{T,\ell}(z),\psi_{T,\ell}(z))$ with $\varphi_{T,\ell},\psi_{T,\ell} \in C^{1,1}(B^{n-1}_{\rho}(0))$ satisfying $\phi_{T,\ell}(0) = \psi_{T,\ell}(0)=0$ and $D \phi_{T,\ell}(0) = D \psi_{T,\ell}(0)=0$ as well as
$$\| D\varphi_{T,\ell} \|_{C(B^{n-1}_{\rho}(0))},\| D\psi_{T,\ell} \|_{C(B^{n-1}_{\rho}(0))} < \epsilon$$
for each $\ell \in \{ 1,\ldots,N \}$; using \eqref{main1} we can also take $\rho = \rho(\epsilon,\partial T) \in (0,1)$ so that for each $x \in \spt \partial T \cap B_{\rho}(0)$ we have $T_{x}^{\perp} \partial T \not \subset T_{x} M_{a}$ for each $a \in \{1,\ldots,A\}$ with $x \in \clos M_{a}.$

\bigskip

We shall proceed by induction on $N \geq 1$ in \eqref{main2}. If $N=1,$ then $0 \in \reg \partial T$ and the theorem follows by Theorem \ref{boundaryregularity}. So we assume $N \geq 2,$ and again by Theorem \ref{boundaryregularity} assume that $0 \in \sing \partial T.$ Our first goal is to define a smooth hypersurface $W$ in $B_{1}(0).$ We will then show that $W$ near $0$ decomposes into finitely many connected components $\{ W_{b} \}_{b=1}^{\mathcal{B}}$ with each $W_{b}$ a $C^{1,1}$ hypersurface-with-boundary as required. Before we define $W,$ we make three sets of observations using Theorem \ref{hyperplanetangentconeregularity}.

\bigskip

First, $0 \in \sing \partial T$ and Theorem \ref{hyperplanetangentconeregularity} imply that every tangent cone of $T$ at $0$ is a sum of half-hyperplanes with constant orientation after rotation. In fact, since $T_{0} \partial T = \R^{n-1}$ then every tangent cone $\C$ of $T$ at $0$ is of the form
\begin{equation} \label{main3}
\C = \sum_{k=1}^{N^{\C}} m^{\C}_{k} Q_{k \#} (\E^{n} \res \{ y \in \R^{n}: y_{n}>0 \})
\end{equation}
for distinct orthogonal rotations $Q_{1},\ldots,Q_{N^{\C}}$ about $\R^{n-1},$ and where $N^{\C},m^{\C}_{1},\ldots,m^{\C}_{N^{\C}} \in \N$ satisfy $\sum_{k=1}^{N^{\C}} m^{\C}_{k} = \sum_{\ell=1}^{N} m_{\ell}.$ 

\bigskip

Second, Theorem 1 of \cite{SS81} and \eqref{main3} imply there exists $\delta_{0}=\delta_{0}(n) \in (0,1)$ with the following property: for any $\delta \in (0,\delta_{0})$ and $\C$ a tangent cone of $T$ at $0,$ there is $\rho = \rho(\delta,T,\C) \in (0,1)$ so that
\begin{equation} \label{main4}
\begin{aligned}
T \res \{ x \in & B_{\rho}(0): |\proj{\{0\} \times \R^{2}}{x}| \geq \delta \rho \} \\
= & \sum_{k=1}^{N^{\C}} \sum_{j=1}^{N^{\C}_{k}} m^{\C}_{j,k} (Q_{k} \circ F_{j,k})_{\#} \left( \E^{n} \res \{ y \in B^{n}_{\rho}(0): y_{n}> \delta \rho/2 \} \right) \\
& \res \{ x \in B_{\rho}(0): |\proj{\{0\} \times \R^{2}}{x}| \geq \delta \rho \},
\end{aligned}
\end{equation}
where for each $k \in \{ 1,\ldots, N^{\C} \}$ we have that $N^{\C}_{k},m^{\C}_{1,k},\ldots,m^{\C}_{N^{\C}_{k},k} \in \N$ satisfy $\sum_{j=1}^{N^{\C}_{k}} m^{\C}_{j,k} = m^{\C}_{k}$; for each $j \in \{ 1,\ldots,N^{\C}_{k} \}$ the map $F_{j,k}$ is given by $F_{j,k}(y) = (y,u_{j,k}(y))$ for $y \in B^{n}_{\rho}(0)$ with $y_{n} > \delta \rho/2$ where 
$$u_{j,k} \in C^{\infty}(\{ y \in B^{n}_{\rho}(0): y_{n}> \delta \rho/2 \}) \text{ with } \| Du_{j,k} \|_{C(\{ y \in B^{n}_{\rho}(0): y_{n}> \delta \rho/2 \})} < c \delta$$ 
for $c=c(n)$; for each $j \in \{ 1,\ldots,N^{\C}_{k}-1\}$ we have $u_{j,k}(y) < u_{j+1,k}(y)$ for each $y \in B^{n}_{\rho}(0)$ with $y_{n}> \delta \rho/2.$

\bigskip

Third, by Theorem 5.4.2 of \cite{F69} we have that
\begin{equation}
\label{main5}
\spt \C \subseteq \bigcup_{\{ a \in \{ 1,\ldots,A \}: x \in \clos M_{a} \}} T_{x} M_{a}
\end{equation}
for any tangent cone $\C$ of $T$ at (any) $x \in \spt T \cap B_{1}(0).$ In case $x \in \reg T \cap B_{1}(0),$ then $T_{x}T = T_{x} M_{a}$ for some $a \in \{ 1,\ldots,A\},$ and \eqref{main1} further implies $\dist_{\HH}(T_{x}T,T_{0}M_{a}) < \epsilon$ for some $a \in \{ 1,\ldots,A \}.$ If $\C$ is a tangent cone of $T$ at $0,$ then $\spt \C \subseteq \bigcup_{a=1}^{A} T_{0}M_{a},$ as $0 \in \clos M_{a}$ for each $a \in \{1,\ldots,A\}$ (see the paragraph containing \eqref{main1}).

\bigskip

Now define the set
$$\begin{aligned}
W = (\spt T & \cap B_{1}(0) \setminus \spt \partial T) \\
& \bigcup \left\{ x \in \spt \partial T \cap B_{1}(0): \begin{aligned} & \text{$T$ has tangent cone at $x$} \\ & \text{which is a hyperplane} \\ & \text{with constant orientation} \\ & \text{but non-constant multiplicity} \end{aligned} \right\}.
\end{aligned}$$
By Theorem \ref{hyperplanetangentconegraph}, the (topological) boundary of $W$ in $B_{1}(0)$ satisfies
\begin{equation} \label{main6}
\partial W \cap B_{1}(0) = \left\{ x \in \spt \partial T \cap B_{1}(0): \begin{aligned} & \text{$T$ has tangent cone at $x$} \\ & \text{which is a sum of} \\ & \text{half-hyperplanes with constant} \\ & \text{orientation after rotation} \end{aligned} \right\}.
\end{equation}
We claim that $W$ is a smooth hypersurface. We prove this in what follows.

\bigskip

First, if $x \in \spt \partial T \cap B_{1}(0)$ and $T$ at $x$ has a tangent cone consisting of a hyperplane with constant orientation but non-constant multiplicity, then $W$ near $x$ is a smooth hypersurface by Theorem \ref{hyperplanetangentconegraph} (or Theorem \ref{hyperplanetangentconeregularity}). Second, if $x \in \spt T \cap B_{1}(0) \setminus \spt \partial T,$ then standard interior regularity for area-minimizing currents shows $x \in \reg T.$ To see this more clearly, let $\C$ be a tangent cone of $T$ at $x.$ Then $\C$ is area-minimizing with $\partial \C = 0,$ and thus the singular set of $\C$ must satisfy Federer's alternatives; see \S37 of \cite{S83}, in particular see Theorem 37.7 of \cite{S83}. We thus conclude $\C$ must be a hyperplane with multiplicity by \eqref{main5} (for this, we can use a proof by induction on the number of distinct planes $T_{x}M_{a}$ with $x \in \clos M_{a}$). It follows that $x \in \reg T$ (using, for example, Theorem 1 of \cite{SS81}), and thus $W$ near $x$ is a smooth hypersurface.

\bigskip

Our goal now is to show that $W$ near $0$ decomposes into finitely many connected components, each of which is a $C^{1,1}$ hypersurface-with-boundary containing the origin. To do this, we fix in what follows $\epsilon=\epsilon(n,\{T_{0}M_{a}\}_{a=1}^{A})>0$ and $\delta = \delta(n,\epsilon) \in (0,\delta_{0})$ to be chosen later. Fix as well any tangent cone $\C$ of $T$ at $0,$ and consider $\rho=\rho(\epsilon,\{M_{a}\}_{a=1}^{A},\partial T,\delta,T,\C)$ so that \eqref{main1}-\eqref{main4} hold.

\bigskip

Consider any $x \in W \cap B_{\rho/2}(0),$ and let $z = \proj{\R^{n-1}}{x}.$ We claim $x \in \gamma((0,1))$ for a Jordan arc
\begin{equation} \label{main7}
\begin{aligned}
& \gamma \in \left\{
\begin{aligned} 
& C([0,1];(\clos W) \cap (\{ z \} \times \R^{2}) \cap \clos B_{\rho}(0)) \\
& C^{\infty}((0,1); W \cap (\{z\} \times \R^{2}) \cap B_{\rho}(0))
\end{aligned} \right. \\
& \text{with } \gamma(0) \in (\partial W) \cap B_{\rho}(0) \text{ while } \gamma(1) \in \partial B_{\rho}(0).
\end{aligned}
\end{equation}
To see this, observe that \eqref{main5} implies that for each $\tilde{x} \in W$ we have $T_{\tilde{x}} W = T_{\tilde{x}} M_{a}$ for some $a \in \{ 1,\ldots,A \}.$ Thus, $\{0\} \times \R^{2} \not \subset T_{\tilde{x}}W$ for each $\tilde{x} \in W \cap B_{\rho}(0)$ by \eqref{main1} (if $\epsilon=\epsilon(n,\{ T_{0}M_{a} \}_{a=1}^{A})>0$ is sufficiently small, and $\rho \in (0,1)$ is chosen depending on $\epsilon,\{M_{a}\}_{a=1}^{A}$). Sard's theorem thus implies $x \in \gamma((0,1))$ for a Jordan arc as claimed in \eqref{main7}, although it remains to show  $\gamma(0) \in (\partial W) \cap B_{\rho}(0)$ while $\gamma(1) \in \partial B_{\rho}(0)$; we show this in what follows.

\bigskip

First, suppose for contradiction $\gamma(0),\gamma(1) \in (\partial W) \cap B_{\rho}(0).$ Let $\gamma(0)=x^{1}$ and $\gamma(1)=x^{2}.$ Thus we conclude by \eqref{main6} that $T$ has at $x^{d}$ for $d=1,2$ a tangent cone
$$\C_{x^{d}} = \sum_{k=1}^{N^{\C_{x^{d}}}} m^{\C_{x^{d}}}_{k} Q^{x^{d}}_{\#} Q^{x^{d}}_{k \#} (\E^{n} \res \{ y \in \R^{n}: y_{n}>0 \})$$
for $N^{\C_{x^{d}}},m^{\C_{x^{d}}}_{1},\ldots,m^{\C_{x^{d}}}_{N^{\C_{x^{d}}}} \in \N,$ an orthogonal rotation $Q^{x^{d}}$ with $\| Q^{x^{d}}-I \| < c \epsilon$ for some $c=c(n)>0$ by \eqref{main2}, and distinct orthogonal rotations $Q^{x^{d}}_{1},\ldots,Q^{x^{d}}_{N^{\C_{x^{d}}}}$ about $\R^{n-1}.$ Choose $\sigma \in (0,\rho)$ so that
$$B_{\sigma}(x^{1}) \cap B_{\sigma}(x^{2}) = \emptyset \text{ and } B_{\sigma}(x^{1}),B_{\sigma}(x^{2}) \subset B_{\rho}(0).$$
By Theorem 1 of \cite{SS81} (as in \eqref{main4}) we can also choose $\sigma \in (0,\rho)$ sufficiently small so that for each $d=1,2$
\begin{equation} \label{main8}
\begin{aligned}
(( & Q^{x^{d}})^{-1} \circ \eta_{x^{d},1})_{\#} T \res \{ \tilde{x} \in B_{\sigma}(0): |\proj{\{0\} \times \R^{2}}{\tilde{x}}| \geq \sigma/2 \} \\
= & \sum_{k=1}^{N^{\C_{x^{d}}}} \sum_{j=1}^{N^{\C_{x^{d}}}_{k}} m^{\C_{x^{d}}}_{j,k} (Q^{x^{d}}_{k} \circ F_{j,k})_{\#} \left( \E^{n} \res \{ y \in B^{n}_{\sigma}(0): y_{n}> \sigma/4 \} \right) \\
& \res \{ \tilde{x} \in B_{\sigma}(0): |\proj{\{0\} \times \R^{2}}{\tilde{x}}| \geq \sigma/2 \}
\end{aligned}
\end{equation}
for $\{ N^{\C_{x^{d}}}_{1},\ldots,N^{\C_{x^{d}}}_{N^{\C_{x^{d}}}}\},\{ m^{\C_{x^{d}}}_{j,k} \}_{j=1,k=1}^{N^{\C_{x^{d}}}_{k},N^{\C_{x^{d}}}} \subset \N$ and $F_{j,k}(y) = (y,u_{j,k}(y))$ for $y \in B^{n}_{\sigma}(0)$ with $y_{n} > \sigma/4$ where $u_{j,k} \in C^{\infty}(\{ y \in B^{n}_{\sigma}(0): y_{n}> \sigma/4 \}).$ Furthermore, for each $k \in \{ 1,\ldots,N^{\C_{x^{d}}} \}$ and $j \in \{ 1,\ldots,N^{\C_{x^{d}}}_{k}-1 \}$ we have that $u_{j,k}(y) < u_{j+1,k}(y)$ for each $y \in B^{n}_{\sigma}(0)$ with $y_{n}> \sigma/4.$ 

\bigskip

Since $\proj{\R^{n-1}}{x^{d}}=z$ and $\| Q^{x^{d}}-I \|< \epsilon,$ then we conclude (as $\epsilon=\epsilon(n,\{T_{0}M_{a}\}_{a=1}^{A})>0$ can be chosen small depending on $n$) that for each $d=1,2$
$$(\eta_{-x^{d},1} \circ Q^{x^{d}})(\{ \tilde{x} \in B_{\sigma}(0): |\proj{\{0\} \times \R^{2}}{\tilde{x}}| \geq \sigma/2 \}) \cap \gamma((0,1)) \neq \emptyset,$$
since $B_{\sigma}(x^{1}) \cap B_{\sigma}(x^{2}) = \emptyset$ and $\gamma(0)=x^{1},\gamma(1)=x^{2}.$ Thus, we can find a Jordan arc
$$g \in C^{\infty}([0,\HH^{1}(g)];\gamma((0,1)))$$
parameterized by arc-length so that
$$\begin{aligned}
g(0) & \in (\eta_{-x_{1},1} \circ Q^{x^{1}})(\{ \tilde{x} \in B_{\sigma}(0): |\proj{\{0\} \times \R^{2}}{\tilde{x}}| < \sigma/2 \}) \text{ while } \\ 
g(\HH^{1}(g)) & \in (\eta_{-x_{2},1} \circ Q^{x^{2}})(\{ \tilde{x} \in B_{\sigma}(0): |\proj{\{0\} \times \R^{2}}{\tilde{x}}| < \sigma/2 \}).
\end{aligned}$$
Define by \eqref{main1},\eqref{main5} the unit-vector field $V \in C^{\infty}(g([0,\HH^{1}(g)]);\{0\} \times S^{1})$ given for $\tilde{x} \in g([0,\HH^{1}(g)])$ by
$$V(\tilde{x}) = \frac{\proj{\{0\} \times \R^{2}}{\ast \vec{T}(\tilde{x})}}{| \proj{\{0\} \times \R^{2}}{\ast \vec{T}(\tilde{x})}|}.$$
Then \eqref{main8} with $d=1$ implies $g,V$ are negatively oriented as in Definition \ref{arcorientation}, while \eqref{main8} with $d=2$ implies $g,V$ are positively oriented, giving a contradiction. We conclude at least one $x^{d} \in \partial B_{\rho}(0).$ 

\bigskip

Second, suppose for contradiction both $x^{1},x^{2} \in \partial B_{\rho}(0).$ Now parameterize
$$\gamma:[0,\HH^{1}(\gamma)] \rightarrow (\clos W) \cap (\{ z \} \times \R^{2}) \cap \clos B_{\rho}(0)$$
by arc-length so that $\gamma(0)=x^{1}$ and $\gamma(\HH^{1}(\gamma))=x^{2}.$ Consider by \eqref{main1},\eqref{main5} 
$$V \in C(\gamma([0,\HH^{1}(\gamma)]);\{0\} \times S^{1}) \cap C^{\infty}(\gamma((0,\HH^{1}(\gamma)));\{0\} \times S^{1})$$
given for $\tilde{x} \in \gamma([0,\HH^{1}(\gamma)])$ by
$$V(\tilde{x}) = \frac{\proj{\{0\} \times \R^{2}}{\ast \vec{T}(\tilde{x})}}{|\proj{\{0\} \times \R^{2}}{\ast \vec{T}(\tilde{x})}|}.$$
Then $\gamma(0)=x^{1} \in \partial B_{\rho}(0),$ \eqref{main4}, and $z = \proj{\R^{n-1}}{x}$ with $x \in B_{\rho/2}(0)$ imply that $\gamma,V$ are positively oriented (again as in Definition \ref{arcorientation}). Contrarily $\gamma(\HH^{1}(\gamma))=x^{2} \in \partial B_{\rho}(0)$ and \eqref{main4} (along with $z = \proj{\R^{n-1}}{x}$ and $x \in B_{\rho/2}(0)$) imply that $\gamma,V$ are negatively oriented. This gives a contradiction, and so we conclude \eqref{main7}.

\bigskip

We now show that $W$ near $0$ decomposes into finitely many connected components, using \eqref{main4}. For this, consider any connected component $W_{1}$ of $W \cap (B^{n-1}_{\rho/4}(0) \times (-\rho/8,\rho/8)^{2}).$ Applying \eqref{main7} to any $x \in W_{1},$ we can assume by \eqref{main4}
$$Q_{1}\left( \gph{\{ y \in B^{n}_{\rho}(0): y_{n} \geq \delta \rho \}}{u_{1,1}} \right) \cap \big( B^{n-1}_{\rho/4}(0) \times (-\rho/8,\rho/8)^{2} \big) \subset W_{1}.$$
Suppose (after rotation about $\R^{n-1}$) that $Q_{1} = I.$ We claim this is the only graph from \eqref{main4} which is contained in $W_{1}.$ We show this in what follows.

\bigskip

First, suppose for contradiction and without loss of generality that
$$Q_{2} \left( \gph{\{ y \in B^{n}_{\rho}(0): y_{n} \geq \delta \rho \}}{u_{1,2}} \right) \cap \big( B^{n-1}_{\rho/4}(0) \times (-\rho/8,\rho/8)^{2} \big) \subset W_{1}$$
Since $W_{1}$ is connected, then there is a curve $g \in C([0,1];W_{1})$ with
$$g(0) \in \gph{\{ y \in B^{n}_{\rho}(0): y_{n} \geq \delta \rho \}}{u_{1,1}} \text{ while } g(1) \in Q_{2}(\gph{\{ y \in B^{n}_{\rho}(0): y_{n} \geq \delta \rho \}}{u_{1,2}}).$$
By \eqref{main1} and \eqref{main5}, for each $t \in [0,1]$ there is $a_{t} \in \{ 1,\ldots,A \}$ so that $\dist_{\HH}(T_{g(t)} W,T_{0}M_{a_{t}}) < \epsilon.$ However, if we choose $\epsilon=\epsilon(n,\{T_{0}M_{a}\}_{a=1}^{A})>0$ sufficiently small and subsequently $\rho \in (0,1)$ sufficiently small (so far depending on $\epsilon,\{M_{a}\}_{a=1}^{A},\partial T$) in \eqref{main1} and \eqref{main2}, then we conclude since $W$ is a smooth hypersurface that $T_{0}M_{a_{t}} = T_{0}M_{a_{0}}$ for each $t \in [0,1].$ Choosing $\delta=\delta(n,\epsilon) \in (0,\delta_{0})$ sufficiently small in \eqref{main4}, and now $\rho=\rho(\epsilon,\{M_{a}\}_{a=1}^{A},\partial T,\delta,T,\C) \in (0,1),$ we conclude $T_{0} M_{a_{0}} = \R^{n}.$ This implies that $Q_{2}$ must be either the identity or the rotation about $\R^{n-1}$ by $\pi.$ Since $Q_{1},\ldots,Q_{N^{\C}}$ are distinct rotations about $\R^{n-1}$ by \eqref{main3}, then $Q_{2}$ must be the rotation about $\R^{n-1}$ by angle $\pi.$ We thus conclude by \eqref{main4} that $\ast \vec{T}(g(0)) \cdot e_{n+1} > 0$ while $\ast \vec{T}(g(1)) \cdot e_{n+1} < 0.$ On the other hand, $\ast \vec{T}(g(t)) \cdot e_{n+1}$ is continuous for $t \in [0,1],$ and $\dist_{\HH}(T_{g(t)}W,\R^{n})<\epsilon$ implies $\ast \vec{T}(g(t)) \cdot e_{n+1} \neq 0$ for each $t \in [0,1].$ This gives a contradiction.

\bigskip

Second, suppose for contradiction 
$$\left( \gph{\{ y \in B^{n}_{\rho}(0): y_{n} \geq \delta \rho \}}{u_{1,1}} \right) \cap \big( B^{n-1}_{\rho/4}(0) \times (-\rho/8,\rho/8)^{2} \big) \subset W_{1}$$
(as assumed above, with $Q_{1}=I$) while
$$\left( \gph{\{ y \in B^{n}_{\rho}(0): y_{n} \geq \delta \rho \}}{u_{2,1}} \right) \cap \big( B^{n-1}_{\rho/4}(0) \times (-\rho/8,\rho/8)^{2} \big) \subset W_{1}$$
as well. Take $g \in C([0,1];W_{1})$ a curve with 
$$g(0) \in \gph{\{ y \in B^{n}_{\rho}(0): y_{n} \geq \delta \rho \}}{u_{1,1}} \text{ while } g(1) \in \gph{\{ y \in B^{n}_{\rho}(0): y_{n} \geq \delta \rho \}}{u_{2,1}}.$$ 
For any $t \in [0,1]$ there is by \eqref{main7} a Jordan arc 
$$\gamma^{t} \in \left\{ 
\begin{aligned}
& C([0,1]; (\clos W) \cap (\{ \proj{\R^{n-1}}{g(t)} \} \times \R^{2}) \cap \clos B_{\rho}(0)) \\
& C^{\infty}((0,1); W \cap (\{ \proj{\R^{n-1}}{g(t)} \} \times \R^{2}) \cap B_{\rho}(0)) 
\end{aligned} \right.$$
with $g(t) \in \gamma^{t}((0,1)),$ and where $\gamma^{t}(0) \in (\partial W) \cap B_{\rho}(0)$ while $\gamma^{t}(1) \in \partial B_{\rho}(0).$ If $\epsilon=\epsilon(n,\{T_{0}M_{a}\}_{a=1}^{A})>0,$ $\delta=\delta(n,\epsilon) \in (0,\delta_{0}),$ and subsequently $\rho=\rho(\epsilon,\{M_{a}\}_{a=1}^{A},\partial T,\delta,T,\C) \in (0,1)$ are chosen sufficiently small, then \eqref{main1}-\eqref{main5} again imply $\dist_{\HH}(T_{g(t)}W,\R^{n}) < \epsilon$ for each $t \in [0,1].$ We as well conclude
$$\begin{aligned}
\gamma^{t}((0,1)) \cap & \{ x \in B_{\rho}(0): |\proj{\{0\} \times \R^{2}}{x}| \geq \delta \rho \} \\
& = \left( \gph{\{ y \in B^{n}_{\rho}(0): y_{n} \geq \delta \rho \}}{u_{j^{t},1}} \right) \cap (\{ \proj{\R^{n-1}}{g(t)} \} \times \R^{2})
\end{aligned}$$
for some $j^{t} \in \{ 1,\ldots,N^{\C}_{1} \}.$ Since $W$ is a smooth hypersurface and $j \in \{ 1,\ldots,N^{\C}_{1}-1 \}$ we have that $u_{j,1}(y)<u_{j+1,1}(y)$ for each $y \in B^{n}_{\rho}(0)$ with $y_{n} \geq \delta \rho,$ then the choice of $j^{t}$ is continuous in $t.$ Finally, $j^{0}=1$ while $j^{1}=2$ gives a contradiction.

\bigskip

Fixing our choice of $\rho \in (0,1),$ we conclude by \eqref{main4} and \eqref{main7} that $W \cap \big( B^{n-1}_{\rho/4}(0) \times (-\rho/8,\rho/8)^{2} \big)$ decomposes precisely into connected components $\{ W_{j,k} \}_{j=1,k=1}^{N^{\C}_{k},N^{\C}},$ where for each $k = 1,\ldots,N^{\C}$ and $j = 1,\ldots,N^{\C}_{k}$
\begin{equation} \label{main9}
Q_{k} \left( \gph{\{ y \in B^{n}_{\rho}(0): y_{n} \geq \delta \rho \}}{u_{j,k}} \right) \cap \big( B^{n-1}_{\rho/4}(0) \times (-\rho/8,\rho/8)^{2} \big) \subset W_{j,k},
\end{equation}
while on the other hand 
\begin{equation} \label{main10}
Q_{\tilde{k}} \left( \gph{\{ y \in B^{n}_{\rho}(0): y_{n} \geq \delta \rho \}}{u_{\tilde{j},\tilde{k}}} \right) \cap W_{j,k} = \emptyset
\end{equation}
if $\tilde{k} \in \{ 1,\ldots,N^{\C} \}$ with $\tilde{k} \neq k,$ or if $\tilde{k}=k$ but $\tilde{j} \in \{ 1,\ldots,N^{\C}_{k} \}$ with $\tilde{j} \neq j.$ 

\bigskip

Next, we show each $W_{j,k}$ is a $C^{1,1}$ hypersurface-with-boundary, with $0 \in \partial W_{j,k}.$ Consider without loss of generality $W_{1,1}.$ We claim 
$$(\partial W_{1,1}) \cap \big( B^{n-1}_{\rho/4}(0) \times (-\rho/8,\rho/8)^{2} \big) = \{ (z,w_{1,1}(z)): z \in B^{n-1}_{\rho/4}(0) \}$$
where $w_{1,1} \in C^{1,1}(B^{n-1}_{\rho/4}(0)).$ For this, suppose for contradiction there are $x^{1},x^{2} \in (\partial W_{1,1}) \cap \big( B^{n-1}_{\rho/4}(0) \times (-\rho/8,\rho/8)^{2} \big)$ with $x^{1} \neq x^{2}$ but so that $\proj{\R^{n-1}}{x^{1}} = \proj{\R^{n-1}}{x^{2}}=z.$ We conclude by \eqref{main7},\eqref{main9},\eqref{main10} that there is a Jordan arc

$$\gamma \in \left\{ \begin{aligned}
& C([0,1]; (\clos W_{1,1}) \cap (\{z\} \times \R^{2})) \\
& C^{\infty}((0,1); W_{1,1} \cap (\{z\} \times \R^{2}))
\end{aligned} \right.$$
with $\gamma((0,1)) = W_{1,1} \cap (\{z\} \times \R^{2})$ and
$$\gamma(0) \in \partial W_{1,1} \cap \big( B^{n-1}_{\rho/4}(0) \times (-\rho/8,\rho/8)^{2} \big) \text{ while } \gamma(1) \in \partial \big( B^{n-1}_{\rho/4}(0) \times (-\rho/8,\rho/8)^{2} \big).$$
Since $x^{1} \neq x^{2},$ then either $\gamma(0) \neq x^{1}$ or $\gamma(0) \neq x^{2}.$ However, if for example $\gamma(0) \neq x^{1},$ then we can choose by \eqref{main6} and Theorem 1 of \cite{SS81} a $\sigma \in (0,\dist(x^{1},\gamma))$ so that we get as in \eqref{main8} with $d=1.$ However, then \eqref{main8} and $\sigma < \dist(x^{1},\gamma)$ imply $\gamma((0,1)) \neq W_{1,1} \cap (\{z\} \times \R^{2}),$ which is a contradiction.

\bigskip

We conclude that $x^{1},x^{2} \in (\partial W_{1,1}) \cap \big(B^{n-1}_{\rho/4}(0) \times (-\rho/8,\rho/8)^{2}\big)$ with $\proj{\R^{n-1}}{x^{1}} = \proj{\R^{n-1}}{x^{2}}$ must satisfy $x^{1}=x^{2}.$ By \eqref{main7} and \eqref{main9} we further conclude
$$(\partial W_{1,1}) \cap \big( B^{n-1}_{\rho/4}(0) \times (-\rho/8,\rho/8)^{2} \big) = \{ (z,w_{1,1}(z)): z \in B^{n-1}_{\rho/4}(0) \}$$
for some function $w_{1,1}:B^{n-1}_{\rho/4}(0) \rightarrow \R^{2}.$ By \eqref{main2},\eqref{main6} there is for each $z \in B^{n-1}_{\rho/4}(0)$ an $\ell_{z} \in \{ 1,\ldots,N \}$ such that $(z,w_{1,1}(z)) = \Phi_{T,\ell_{z}}(z).$ We now show that $w_{1,1} \in C^{1,1}(B^{n-1}_{\rho/4}(0);\R^{2}).$ Fix $z \in B^{n-1}_{\rho/4}(0),$ we argue two cases. 

\bigskip

First, suppose there is $\ell \in \{ 1,\ldots,N\}$ so that $\Phi_{T,\ell}(z) \neq (z,w_{1,1}(z)).$ Then there is $\sigma \in (0,1-|(z,w_{1,1}(z)|)$ so that by \eqref{main2}, after relabeling,
$$\partial T \res B_{\sigma}((z,w_{1,1}(z))) = \sum_{\ell=1}^{N_{z}} m_{\ell} \Phi_{T,\ell \#}(\E^{n-1} \res B^{n-1}_{\sigma}(z)) \res B_{\sigma}((z,w_{1,1}(z)))$$
where $N_{z} \in \{ 1,\ldots,N-1 \}.$ Furthermore, by \eqref{main2} we have $T^{\perp}_{(z,w_{1,1}(z))} \partial T \not \subset T_{(z,w_{1,1}(z))} M_{a}$ for each $a \in \{1,\ldots,A\}$ so that $(z,w_{1,1}(z)) \in \clos M_{a}.$ Recall that the theorem is true by induction for $N_{z} \in \{ 1,\ldots,N-1 \}$ (see the paragraph after \eqref{main2}). If we consider $T \res B_{\sigma}((z,w_{1,1}(z))),$ then we conclude $w_{1,1}|_{B^{n-1}_{\sigma}(z)} \in C^{1,1}(B^{n-1}_{\sigma}(z);\R^{2}).$ 

\bigskip

Second, suppose $\Phi_{T,\ell}(z)=(z,w_{1,1}(z))$ for each $\ell \in \{1,\ldots,N\}.$ Since for each $\tilde{z} \in B^{n-1}_{\rho/4}(0)$ we have $(\tilde{z},w_{1,1}(\tilde{z})) = \Phi_{T,\ell_{\tilde{z}}}(\tilde{z})$ for some $\ell_{\tilde{z}} \in \{1,\ldots,N\},$ then by Definition \ref{immersedboundary} it follows that $Dw_{1,1}(z)$ exists and 
$$Dw_{1,1}(z) = (D \varphi_{T,1}(z),D \psi_{T,1}(z)) = \ldots = (D\varphi_{T,N}(z),D \psi_{T,N}(z))$$ 
with $\varphi_{T,1},\ldots,\varphi_{T,N}$ and $\psi_{T,1},\ldots,\psi_{T,N}$ as in \eqref{main2}. 

\bigskip

The two cases together with Definition \ref{immersedboundary} imply $w_{1,1} \in C^{1,1}(B^{n-1}_{\rho/4}(0);\R^{2}).$ Next, we show $W_{1,1}$ is a $C^{1,1}$ hypersurface-with-boundary, using \cite{A75}. For this, we claim 
\begin{equation} \label{main11}
\frac{\HH^{n}(W_{1,1} \cap B_{\rho/8}(0))}{\omega_{n} (\rho/8)^{n}} < \frac{1+c \delta}{2}
\end{equation}
where $\omega_{n}=\HH^{n}(B_{1}(0))$ and $c=c(n)>0.$ For this, first observe that by \eqref{main4},\eqref{main9},\eqref{main10} we can compute
$$\frac{\HH^{n}(W_{1,1} \cap \{ x \in B_{\rho/8}(0): |\proj{\{0\} \times \R^{2}}{x}| \geq \delta \rho  \})}{\omega_{n} (\rho/8)^{n}} \leq \frac{\sqrt{1+c_{1} \delta^{2}}}{2}$$
for some $c_{1}=c_{1}(n).$ Meanwhile, using Theorem \ref{monotonicity} we can compute
$$\frac{\HH^{n}(W_{1,1} \cap \{ x \in B_{\rho/8}(0): |\proj{\{0\} \times \R^{2}}{x}| \leq \delta \rho  \})}{\omega_{n} \rho^{n}} \leq c_{2} \delta$$
for some $c_{2}=c_{2}(n).$ These two calculations together show \eqref{main11}.

\bigskip

Consider the varifold $|W_{1,1}|$ associated to $W_{1,1}.$ By \eqref{main7} and \eqref{main9} (with $z=0$) we have $0 \in \spt W_{1,1}.$ Observe that $W_{1,1}$ is stationary in $B_{\rho/8}(0) \setminus \{ (z,w_{1,1}(z)) : z \in B^{n-1}_{\rho/8}(0) \},$ where $w_{1,1} \in C^{1,1}(B^{n-1}_{\rho/4}(0);\R^{2}).$ Choosing $\delta= \delta(n,\epsilon) \in (0,\delta_{0})$ sufficiently small depending on $n,$ then \eqref{main11} and \cite{A75} imply that $W_{1,1}$ is a $C^{1,1}$ hypersurface-with-boundary in $B_{\rho/16}(0).$ 

\bigskip

Finally, $W_{1,1} \cap \spt \partial T \subseteq \reg \partial T$ follows by Theorem \ref{hyperplanetangentconeregularity}. $\square$

\section{Tangentially Immersed Boundaries in Space.}

In this section we discuss two-dimensional area-minimizing currents in $\R^{3}$ with $C^{1,\alpha}$ tangentially immersed boundaries having co-oriented Lipschitz mean curvature. In this setting, the main results of the previous section can be sharpened. In particular, Theorems \ref{hyperplanetangentconeregularity},\ref{main} hold for $T \in \TI^{1,\alpha}_{2,loc}(U)$ for $\alpha \in (0,1]$ more generally; see respectively Theorems \ref{spaceplanetangentconeregularity},\ref{spacemain}.

\bigskip

To begin, we give the following lemma which will, in Theorem \ref{spacemain}, allow us to give a more general version of Theorem \ref{main}.

\begin{lemma} \label{spacesubmain} Let $U \subseteq_{o} \R^{3},$ $\alpha \in (0,1],$ and suppose $T \in \I^{1,\alpha}_{2,loc}(U).$ Suppose $x \in \spt \partial T$ and there exists $\rho \in (0,\dist(x,\partial U))$ so that
$$\spt T \cap B_{\rho}(x) \subseteq \bigcup_{a=1}^{A} (\clos M_{a}) \cap B_{\rho}(x)$$
for $C^{1}$ hypersurfaces-with-boundary $M_{1},\ldots,M_{a}$ in $B_{\rho}(x).$ Then there is $\sigma \in (0,\rho)$ and $\mathcal{A} \subseteq \{ 1,\ldots,A\}$ so that
$$\spt T \cap B_{\sigma}(x) \subseteq \bigcup_{a \in \mathcal{A}} (\clos M_{a}) \cap B_{\sigma}(x),$$
and $x \in \clos M_{a}$ with $T_{x} \partial T \subset T_{x} M_{a}$ for each $a \in \mathcal{A}.$
\end{lemma}

{\bf Proof:} Suppose (after translation) $0 \in \spt \partial T$ and that (after rotation) $T_{0} \partial T = \R.$ We can also suppose $\rho \in (0,\dist(0,\partial U))$ is such that (after relabeling) $\spt T \cap B_{\rho}(0) \subseteq \bigcup_{a=1}^{A} (\clos M_{a}) \cap B_{\rho}(0)$ with $0 \in \clos M_{a}$ for each $a \in \{1,\ldots,A\}.$ Define
$$\mathcal{A} = \{ a \in \{ 1,\ldots,A \}: \R \subset T_{0} M_{a} \}.$$
We claim there is $\sigma \in (0,\rho)$ so that
$$\spt T \cap B_{\sigma}(0) \subseteq \bigcup_{a \in \mathcal{A}} (\clos M_{a}) \cap B_{\sigma}(0).$$
It suffices to discount the following two scenarios:

\bigskip

First, suppose for contradiction there exists a sequence 
$$\{ x^{d} \}_{d=1}^{\infty} \subset \spt \partial T \cap B_{\rho}(0) \setminus \bigcup_{a \in \mathcal{A}} \clos M_{a}$$
so that $x^{d} \rightarrow 0.$ For each $d \in \N$ let $\C_{x^{d}}$ be a tangent cone of $T$ at $x^{d}.$ By Theorem 5.4.2 of \cite{F69} and Theorem \ref{monotonicity} we have
$$T_{x^{d}} \partial T \subseteq \spt \C_{x^{d}} \subseteq \bigcup_{\{ a \in \{ 1,\ldots,A \} \setminus \mathcal{A}: x^{d} \in \clos M_{a} \}} T_{x^{d}} M_{a}.$$
Recall that $0 \in \clos M_{a}$ for each $a \in \{ 1,\ldots,A \}.$ Since for each $a \in \{ 1,\ldots,A \} \setminus \mathcal{A}$ we have that $M_{a}$ is a $C^{1}$ hypersurface-with-boundary with $\R \not \subset T_{0} M_{a},$ then we conclude there is $\delta>0$ small so that for all sufficiently large $d$ we have $\dist_{\HH}(\R,T_{x^{d}} \partial T) \geq \delta.$ This contradicts $T \in \TI^{1,\alpha}_{2,loc}(U).$

\bigskip

Second, suppose there exists a sequence
$$\{ x^{d} \}_{d=1}^{\infty} \subset \spt T \cap B_{\rho}(0) \setminus \left( \spt \partial T \cup \bigcup_{a \in \mathcal{A}} \clos M_{a} \right);$$
we can suppose $\{ x^{d} \}_{d=1}^{\infty} \subset \clos M_{a_{0}}$ for some fixed $a_{0} \in \{1,\ldots,A\} \setminus \mathcal{A}.$ By Theorem 5.4.2 of \cite{F69} and interior regularity of area-minimizing currents
\begin{equation}
\label{ssmain1}
T_{x^{d}} T = T_{x^{d}} M_{a} \text{ for some } a \in \{ 1,\ldots,A \} \setminus \mathcal{A} \text{ with } x^{d} \in \clos M_{a}.
\end{equation}
For each $d \geq 1$ let $\rho_{d} = 2|x^{d}|,$ so $\rho_{d}>0.$ 

\bigskip

By taking a subsequence and relabeling, we can suppose that $\C$ is a tangent cone of $T$ at $0$ with $\eta_{0,\rho_{d} \#} T \rightarrow \C.$ By Theorem 5.4.2 of \cite{F69} and Theorem \ref{monotonicity} we again have $\spt \C \subseteq \bigcup_{a=1}^{A} T_{0} M_{a}.$ In fact,
\begin{equation}
\label{ssmain2}
\spt \C \subseteq \bigcup_{a \in \mathcal{A}} T_{0} M_{a}.
\end{equation}  
To see this, suppose for contradiction $\tilde{x} \in \spt \C \setminus \bigcup_{a \in \mathcal{A}} T_{0} M_{a}.$ Then $\spt \C \subseteq \bigcup_{a=1}^{A} T_{0} M_{a},$ the constancy theorem (see Theorem 26.27 of \cite{S83}), and $\spt \partial \C = \R$ imply there is $\tilde{\rho}>0$ and $a \in \{1,\ldots,A\} \setminus \mathcal{A}$ so that $\spt \C \cap B_{\tilde{\rho}}(\tilde{x}) = T_{0}M_{a} \cap B_{\tilde{\rho}}(\tilde{x}).$ However, $\R \not \subset T_{0}M_{a}$ contradicts Theorem \ref{monotonicity}, and so we conclude \eqref{ssmain2}.

\bigskip

Fix $\delta = \delta(\{ T_{0}M_{a}\}_{a=1}^{A}) \in (0,\delta_{0})$ to be chosen later, with $\delta_{0}>0$ as in Theorem 1 of \cite{SS81}. Recall that $M_{1},\ldots,M_{A}$ are $C^{1}$ hypersurfaces-with-boundary, and $\R \subset T_{0} M_{a}$ for each $a \in \mathcal{A}$ while $\R \not \subset T_{0} M_{a}$ for each $a \in \{1,\ldots,A\} \setminus \mathcal{A}.$ If we choose $\delta = \delta(\{ T_{0}M_{a}\}_{a=1}^{A}) \in (0,\delta_{0})$ sufficiently small, then \eqref{ssmain1} implies that for all sufficiently large $d \in \N$
$$\dist_{\HH}(T_{x^{d}}T,T_{0}M_{a}) \geq \delta \text{ for each } a \in \mathcal{A}.$$
On the other hand, Theorem 1 of \cite{SS81}, $\eta_{0,\rho_{d} \#} T \rightarrow \C,$ and \eqref{ssmain2} imply that for all sufficiently large integers $d$
$$\dist_{\HH}(T_{x}T,T_{0} M_{a}) < \delta \text{ for some } a \in \mathcal{A}$$
for each $x \in B_{\rho_{d}}(0)$ with $|\proj{\{0\} \times \R^{2}}{x}| > \delta \rho_{d}$; see for example \eqref{main3} in the proof of Theorem \ref{main}. Since $\rho_{d} = 2|x^{d}|,$ then we conclude 
\begin{equation} \label{ssmain3}
|\proj{\{0\} \times \R^{2}}{x_{d}}| \leq \delta \rho_{d}
\end{equation}
for all sufficiently large $d \in \N.$

\bigskip

Now consider $x^{d} = \proj{T_{0}M_{a_{0}}}{x^{d}}+\proj{T^{\perp}_{0}M_{a_{0}}}{x^{d}}.$ Since $\R \not \subset T_{0}M_{a_{0}},$ then there exists $\delta_{a_{0}}>0$ so that for each $x \in T_{0}M_{a_{0}}$ we have $|\proj{\{0\} \times \R^{2}}{x}| \geq \delta_{a_{0}} |x|$; thus
$$|\proj{\{0\} \times \R^{2}}{\proj{T_{0} M_{a_{0}}}{x^{d}}}| \geq \delta_{a_{0}} |\proj{T_{0} M_{a_{0}}}{x^{d}}|.$$
Since $M_{a_{0}}$ is a $C^{1}$ hypersurface-with-boundary and $x^{d} \in \clos M_{a_{0}}$ for each $d \in \N,$ then for all sufficiently large $d \in \N$
$$|\proj{T^{\perp}_{0} M_{a_{0}}}{x^{d}}| \leq \frac{\delta_{a_{0}}}{2} |\proj{T_{0} M_{a_{0}}}{x^{d}}|.$$
Since $\rho_{d} = 2|x^{d}| \leq 2 |\proj{T_{0}M_{a_{0}}}{x^{d}}| +  2 |\proj{T^{\perp}_{0}M_{a_{0}}}{x^{d}}|$ then
$$|\proj{T_{0}M_{a_{0}}}{x^{d}}| \geq \frac{\rho_{d}}{2+\delta_{a_{0}}}.$$
Using these three computations together, we conclude for all sufficiently large $d \in \N$
$$\begin{aligned}
|\proj{ \{0\} \times \R^{2}}{x_{d}}| & \geq |\proj{\{0\} \times \R^{2}}{\proj{T_{0} M_{a_{0}}}{x^{d}}}| - |\proj{\{0\} \times \R^{2}}{\proj{T^{\perp}_{0} M_{a_{0}}}{x^{d}}}| \\
& \geq \delta_{a_{0}} |\proj{T_{0} M_{a_{0}}}{x^{d}}| - \frac{\delta_{a_{0}}}{2} |\proj{T_{0} M_{a_{0}}}{x^{d}}| \geq \frac{\delta_{a_{0}} \rho_{d}}{4+2\delta_{a_{0}}}.
\end{aligned}$$
This contradicts \eqref{ssmain3}, choosing $\delta = \delta(\{ T_{0}M_{a}\}_{a=1}^{A}) \in (0,\delta_{0})$ small. $\square$

\bigskip

We now give a different definition of half-regular points, differing slightly from that given by Lemma \ref{halfregular}. This definition of half-regular points appears in Lemma 1 of \cite{R16}, in the setting of two-dimensional solutions to the $c$-Plateau problem in space; see \S 1.4.

\begin{lemma} \label{spacehalfregular} Let $U \subseteq_{o} \R^{3},$ $\alpha \in (0,1],$ and $T \in \TI^{1,\alpha}_{2,loc}(U).$ For any $x \in \sing \partial T$ and $\rho \in (0,\dist(x,\partial U)),$ there is $\x \in \sing \partial T$ and $\sigma \in (0,\rho-|x-\x|)$ so that
$$\partial T \res B_{\sigma}(\x) = \sum_{\ell=1}^{N} m_{\ell} \Big[ (\eta_{-\x,1} \circ Q \circ \Phi_{T,\ell})_{\#}(\E^{1} \res (-\sigma,\sigma)) \Big] \res B_{\sigma}(\x);$$
for $N,m_{1},\ldots,m_{N} \in \N,$ an orthogonal rotation $Q,$ and $\Phi_{T,\ell} \in C^{1,\alpha}((-\sigma,\sigma);\R^{3})$ for each $\ell = 1,\ldots,N$ is the map $\Phi_{T,\ell}(z) = (z,\varphi_{T,\ell}(z),\psi_{T,\ell}(z))$ where $\phi_{T,\ell},\psi_{T,\ell} \in C^{1,\alpha}((-\sigma,\sigma))$ satisfy
$$\varphi_{T,\ell}(0) = \psi_{T,\ell}(0) = 0 \text{ and } D \varphi_{T,\ell}(0) = D \psi_{T,\ell}(0)=0.$$
Moreover, one of the two following occurs:

\begin{itemize}
 \item For each $\ell,\tilde{\ell} \in \{ 1,\ldots,N \}$ the images
$$\Phi_{T,\ell}((0,\sigma)) \text{ and } \Phi_{T,\tilde{\ell}}((0,\sigma))$$
are either equal or disjoint. Moreover $\Phi_{T,\ell}((0,\sigma)) \cap \Phi_{T,\tilde{\ell}}((0,\sigma)) = \emptyset$ for some $\ell,\tilde{\ell} \in \{ 1,\ldots,N \}.$
 \item For each $\ell,\tilde{\ell} \in \{ 1,\ldots,N \}$ the images
$$\Phi_{T,\ell}((-\sigma,0)) \text{ and } \Phi_{T,\tilde{\ell}}((-\sigma,0))$$
are either equal or disjoint. Moreover $\Phi_{T,\ell}((-\sigma,0)) \cap \Phi_{T,\tilde{\ell}}((-\sigma,0)) = \emptyset$ for some $\ell,\tilde{\ell} \in \{ 1,\ldots,N \}.$
\end{itemize}

\bigskip

We say that such a point $\x \in \sing \partial T$ is half-regular. \end{lemma}

{\bf Proof:} We can either modify the proof of Lemma \ref{halfregular}, or more directly merely apply the proof of Lemma 1 of \cite{R16}. $\square$

\bigskip

Our aim is now to show that Theorem \ref{hyperplanetangentconeregularity} holds for $T \in \TI^{1,\alpha}_{2,loc}(U)$ with $\alpha \in (0,1]$ more generally if $U \subseteq_{o} \R^{3}.$ The first step is to give the more general version of Lemma \ref{halfregulartangentcones}.

\begin{lemma} \label{spacehalfregulartangentcones} Let $U \subseteq_{o} \R^{3},$ $\alpha \in (0,1],$ and suppose $T \in \TI^{1,\alpha}_{2,loc}(U),$ where $\partial T$ has co-oriented mean curvature $H_{\partial T} = h \nu_{T}$ with $h:U \rightarrow \R$ Lipschitz. If $\x \in \sing \partial T$ is half-regular (as in Lemma \ref{spacehalfregular}), then every tangent cone of $T$ at $\x$ is a sum of half-planes with constant orientation after rotation (as in Definition \ref{conesdefinition} with $n=2$).
\end{lemma}

{\bf Proof:} Suppose for contradiction that (after translation) $0 \in \sing \partial T$  is half-regular and that $T$ at $0$ has tangent cone 
\begin{equation} \label{shrtc1}
\C = (m+\theta) \E^{2} \res \{ y \in \R^{2} : y_{2}>0 \} + \theta \E^{2} \res \{ y \in \R^{2}: y_{2}<0 \}
\end{equation}
where $m,\theta \in \N.$ Suppose as well that $\rho \in (0,\dist(0,\partial U))$ is such that
$$\partial T \res B_{\rho}(0) = \sum_{\ell=1}^{N} m_{\ell} \Phi_{T,\ell \#}(\E^{1} \res (-\rho,\rho)) \res B_{\rho}(0)$$
as in Lemma \ref{spacehalfregular}, and where (without loss of generality, and after relabeling)
\begin{equation} \label{shrtc2}
\begin{aligned}
& \text{for each $\ell,\tilde{\ell} \in \{ 1,\ldots,N \}$ the images $\Phi_{T,\ell}((0,\rho))$ and $\Phi_{T,\tilde{\ell}}((0,\rho))$} \\
& \text{are either equal or disjoint, while $\Phi_{T,1}((0,\rho)) \cap \Phi_{T,2}((0,\rho)) = \emptyset.$}
\end{aligned} 
\end{equation}
We can also suppose by Theorem \ref{hyperplanetangentconegraph} there exists a solution to the minimal surface equation $u \in C^{\infty}(B^{2}_{\rho}(0))$ with $u(0)=0$ and $Du(0)=0$ so that
$$\spt T \cap B_{\rho}(0) = \left( \gph{B^{2}_{\rho}(0)}{u} \right) \cap B_{\rho}(0),$$
where the orientation vector for $T$ if $x \in \spt T \cap B_{\rho}(0)$ is given by
\begin{equation} \label{shrtc3}
\ast \vec{T}(x) = \left. \left( \frac{-Du}{\sqrt{1+|Du|^{2}}},\frac{1}{\sqrt{1+|Du|^{2}}} \right) \right|_{\proj{\R^{2}}{x}}.
\end{equation}

\bigskip

For $\ell = 1,2$ let $\sigma_{\ell}$ be the curve parameterized by arc-length with image $\Phi_{T,\ell}([0,\rho)).$ Let $P:\R^{3} \rightarrow \R^{3}$ be the map $P(x) = (x_{2},-x_{1},0).$ We conclude by \eqref{shrtc1},\eqref{shrtc2},\eqref{shrtc3} that there is an $\epsilon>0$ so that for each $\ell=1,2$
$$\sigma_{\ell}''(t) = -h(\sigma_{\ell}(t)) \frac{P(\sigma_{\ell}'(t))- \left( P(\sigma_{\ell}'(t)) \cdot \ast \vec{T}(\sigma_{\ell}) \right) \ast \vec{T}(\sigma_{\ell})}{\| P(\sigma_{\ell}'(t))- \left( P(\sigma_{\ell}'(t)) \cdot \ast \vec{T}(\sigma_{\ell}) \right) \ast \vec{T}(\sigma_{\ell}) \|}$$
for each $t \in (0,\epsilon).$ Since $h$ is Lipschitz while 
$$\sigma_{1}(0)=\sigma_{2}(0)=0 \text{ and } \sigma_{1}'(0)=\sigma_{2}'(0)=e_{1},$$
then we conclude by uniqueness of ODEs that $\sigma_{1}=\sigma_{2}.$ This contradicts $\Phi_{T,1}((0,\rho)) \cap \Phi_{T,2}((0,\rho)) = \emptyset.$ $\square$

\bigskip

As in the the higher dimensional case, Lemma \ref{spacehalfregulartangentcones} readily implies the following regularity theorem:

\begin{theorem} \label{spaceplanetangentconeregularity} Let $U \subseteq_{o} \R^{3},$ $\alpha \in (0,1],$ and suppose $T \in \TI^{1,\alpha}_{n,loc}(U)$ where $\partial T$ has co-oriented mean curvature $H_{\partial T} = h \nu_{T}$ with $h:U \rightarrow \R$ Lipschitz. If $T$ at $x \in \spt \partial T$ has tangent cone a plane with constant orientation but non-constant multiplicity (as in Definition \ref{conesdefinition} with $n=2$), then $x \in \reg \partial T.$ \end{theorem}

{\bf Proof:} Follows exactly as the proof of Theorem \ref{hyperplanetangentconeregularity}, where we now use Lemma \ref{spacehalfregulartangentcones}. $\square$

\bigskip

Finally, we give a more general version of Theorem \ref{main} in case $n=2$:

\begin{theorem} \label{spacemain} Let $U \subseteq_{o} \R^{3},$ $\alpha \in (0,1],$ and suppose $T \in \TI^{1,\alpha}_{2,loc}(U)$ where $\partial T$ has co-oriented mean curvature $H_{\partial T} = h \nu_{T}$ with $h:U \rightarrow \R$ Lipschitz. Suppose $x \in \spt \partial T$ and that there exists $\rho \in (0,\dist(x,\partial U))$ and $C^{1}$ hypersurfaces-with-boundary $M_{1},\ldots,M_{A}$ in $B_{\rho}(x)$ so that 
$$\spt T \cap B_{\rho}(x) \subseteq \bigcup_{a=1}^{A} (\clos M_{a}) \cap B_{\rho}(x).$$
Then there is $\sigma \in (0,\rho)$ and $\mathcal{B} \in \{ 1,\ldots,2 \Theta_{T}(x) \}$ so that 
$$\spt T \cap B_{\sigma}(x) = \bigcup_{b=1}^{\mathcal{B}} (\clos W_{b}) \cap B_{\sigma}(x),$$
for orientable $C^{1,\alpha}$ hypersurfaces-with-boundary $W_{1},\ldots,W_{\mathcal{B}}$ in $B_{\sigma}(x).$ For each $b \in \{ 1,\ldots,\mathcal{B} \}$ we have $x \in \partial W_{b}$ and $W_{b} \cap \spt \partial T \subset \reg \partial T.$ Furthermore, for each $b,\tilde{b} \in \{ 1,\ldots,\mathcal{B} \}$ we have
$$(\clos W_{b}) \cap (\clos W_{\tilde{b}}) \cap B_{\sigma}(x) \subseteq (\partial W_{b}) \cap (\partial W_{\tilde{b}}) \cap B_{\sigma}(x).$$
\end{theorem}

{\bf Proof:} We follow the proof of Theorem \ref{main}, where we use Lemmas \ref{spacesubmain},\ref{spaceplanetangentconeregularity} in place of Lemmas \ref{submain},\ref{hyperplanetangentconeregularity}. The only other significant difference is in the argument following \eqref{main11}; we describe this difference. At that analogous point in the argument, we show there is $\sigma \in (0,\rho)$ sufficiently small so that $\spt T \cap B_{\sigma}(x) = \bigcup_{b=1}^{\mathcal{B}} (\clos W_{b}) \cap B_{\sigma}(x)$ where $W_{1},\ldots,W_{\mathcal{B}}$ are pairwise disjoint smooth surfaces with \emph{topological} boundary $(\partial W_{b}) \cap B_{\rho}(x)$ a $C^{1,\alpha}$ Jordan arc through $x$ for each $b \in \{ 1,\ldots, \mathcal{B} \}.$ For each $b \in \{ 1,\ldots,\mathcal{B} \}$ we have that \eqref{main11} holds as well, more specifically that for each $b \in \{ 1,\ldots,\mathcal{B} \}$ we have $\frac{\HH^{2}(W_{b} \cap B_{\sigma}(x))}{\pi \sigma^{2}} < \frac{1+\delta}{2},$ and that $|W_{b}|$ the varifold associated to $W_{b}$ is stationary in $B_{\sigma}(x) \setminus \partial W_{b}.$ We conclude in this case by \cite{B15} that each $W_{1},\ldots,W_{\mathcal{B}}$ is a $C^{1,\alpha}$ surface-with-boundary. $\square$

\appendix
\section{Appendix}

In this section we give the technical Lemma \ref{sardslemma}, which we need for the proof of Theorem \ref{hyperplanetangentconegraph}. We work in $\R^{2},$ so that for $x=(x_{1},x_{2}) \in \R^{2}$ we have that $B_{\rho}(x)$ is the \emph{disk} of radius $\rho>0$ centered at $x.$ First, we make the following useful definition:

\begin{definition} \label{arcorientation} Suppose $\gamma \in C([0,\HH^{1}(\gamma)];\R^{2}) \cap C^{\infty}((0,\HH^{1}(\gamma));\R^{2})$ is a Jordan arc parameterized by arc length; write $\gamma(t)=(\gamma_{1}(t),\gamma_{2}(t)).$ Also suppose $V \in C^{\infty}(\gamma((0,\HH^{1}(\gamma)));S^{1}).$ Then we say \emph{$\gamma,V$ are positively oriented} if $V(\gamma(t))=(-\gamma_{2}'(t),\gamma_{1}'(t))$ for each $t \in (0,\HH^{1}(\gamma)).$ Contrarily, we say \emph{$\gamma,V$ are negatively oriented} if $V(\gamma(t))=(\gamma_{2}'(t),-\gamma_{1}'(t))$ for each $t \in (0,\HH^{1}(\gamma)).$ \end{definition}

The following general lemma about curves then is proved geometrically:

\begin{lemma} \label{sardslemma} Suppose $\{ B_{\sigma_{\ell}}(x_{\ell}) \}_{\ell=1}^{N}$ is a collection of pairwise disjoint balls with $B_{\sigma_{\ell}}(x_{\ell}) \subset B_{1/2}(0) \subset \R^{2}$ for each $\ell \in \{ 1,\ldots,N \}.$ Then (a) and (b) hold:
\begin{enumerate} 
 \item[(a)] With $N^{(1)} \in \{ 1,\ldots,N \}$ and $N_{1},\ldots,N_{N^{(1)}} \in \N,$ suppose we have a collection of Jordan arcs
$$\Gamma_{1},\Gamma_{2}, \{ \{ \gamma^{\ell}_{k} \}_{k=1}^{N_{\ell}} \}_{\ell=1}^{N^{(1)}}, \{ G^{\ell},g^{\ell} \}_{\ell=N^{(1)}+1}^{N},$$
with images that are pairwise either equal or meet only at mutual endpoints, and let $\mathcal{E}$ be the collection of endpoints of these arcs; we allow for the collection $\{ G^{\ell},g^{\ell} \}_{\ell=N^{(1)}+1}^{N}$ to be empty. With
$$\mathcal{G} = \Gamma_{1} \cup \Gamma_{2} \cup \left( \bigcup_{\ell=1}^{N^{(1)}} \bigcup_{k=1}^{N_{\ell}} \gamma^{\ell}_{k} \right) \cup \left( \bigcup_{\ell=N^{(1)}+1}^{N} G^{\ell} \cup g^{\ell} \right)$$
suppose there is $V \in C^{\infty}(\mathcal{G} \setminus \mathcal{E}; S^{1})$ so that the following hold:
\begin{enumerate}
 \item[(1)] The curves $\Gamma_{1},\Gamma_{2}$ satisfy
 $$\mathcal{G} \cap (\clos B_{1}(0)) \setminus B_{1/2}(0) = (\Gamma_{1} \cup \Gamma_{2}) \cap (\clos B_{1}(0)) \setminus B_{1/2}(0).$$

\bigskip

 We parameterize 
 $$\Gamma_{1} \in C([0,\HH^{1}(\Gamma_{1})];\clos B_{1}(0)) \cap C^{\infty}((0,\HH^{1}(\Gamma_{1}));B_{1}(0))$$ 
 and 
 $$\Gamma_{2} \in C([0,\HH^{1}(\Gamma_{2})];\clos B_{1}(0)) \cap C^{\infty}((0,\HH^{1}(\Gamma_{2}));B_{1}(0))$$ 
 by arc-length and so that $\Gamma_{1}(0),\Gamma_{2}(0) \in \partial B_{1}(0)$ while 
$$\Gamma_{1}(\HH^{1}(\Gamma_{1})),\Gamma_{2}(\HH^{1}(\Gamma_{2})) \in \bigcup_{\ell=1}^{N} B_{\sigma_{\ell}/2}(x_{\ell}).$$
The images $\Gamma_{1}([0,\HH^{1}(\Gamma_{1})))$ and $\Gamma_{2}([0,\HH^{1}(\Gamma_{2})))$ are disjoint.
 
 \item[(2)] For each $\ell \in \{ 1,\ldots,N^{(1)} \}$
 $$\mathcal{G} \cap (\clos B_{\sigma_{\ell}}(x_{\ell})) \setminus B_{\sigma_{\ell}/2}(x_{\ell}) = \bigcup_{k=1}^{N_{\ell}} \gamma^{\ell}_{k} \cap (\clos B_{\sigma_{\ell}}(x_{\ell})) \setminus B_{\sigma_{\ell}/2}(x_{\ell}).$$

 \bigskip
 
 We parameterize for each $k \in \{ 1,\ldots,N_{\ell} \}$
 $$\gamma^{\ell}_{k} \in C([0,\HH^{1}(\gamma^{\ell}_{k})];\clos B_{1}(0)) \cap C^{\infty}((0,\HH^{1}(\gamma^{\ell}_{k}));B_{1}(0))$$
 by arc-length so that $\gamma^{\ell}_{k}(\HH^{1}(\gamma^{\ell}_{k})) \in B_{\sigma_{\ell}/2}(x_{\ell}).$ Meanwhile, $\gamma^{\ell}_{k}(0) \in (\partial B_{1}(0)) \cup \bigcup_{\tilde{\ell}=1}^{N} B_{\sigma_{\tilde{\ell}}}(x_{\tilde{\ell}})$ and $\gamma^{\ell}_{k} \cap \partial B_{\sigma_{\ell}}(x_{\ell}) \neq \emptyset.$
 
\bigskip

With this parameterization $\gamma^{\ell}_{k},V$ are positively oriented.

 \item[(c)] For each $\ell$ in the (possibly empty) set $\{ N^{(1)}+1,\ldots,N \}$
  $$\mathcal{G} \cap (\clos B_{\sigma_{\ell}}(x_{\ell})) \setminus B_{\sigma_{\ell}/2}(x_{\ell}) = \left( G^{\ell} \cup g^{\ell} \right) \cap (\clos B_{\sigma_{\ell}}(x_{\ell})) \setminus B_{\sigma_{\ell}/2}(x_{\ell}).$$

 \bigskip
       
 We parameterize 
 $$G^{\ell} \in C([0,\HH^{1}(G^{\ell})];\clos B_{1}(0)) \cap C^{\infty}((0,\HH^{1}(G^{\ell}));B_{1}(0))$$
 and 
 $$g^{\ell} \in C([0,\HH^{1}(g^{\ell})];\clos B_{1}(0)) \cap C^{\infty}((0,\HH^{1}(g^{\ell}));B_{1}(0))$$
 by arc-length so that $G^{\ell}(\HH^{1}(G^{\ell})) ,g^{\ell}(\HH^{1}(g^{\ell})) \in B_{\sigma_{\ell}/2}(x_{\ell}).$ Meanwhile, $G^{\ell}(0) ,g^{\ell}(0) \in (\partial B_{1}(0)) \cup \bigcup_{\tilde{\ell}=1}^{N} B_{\sigma_{\tilde{\ell}}/2}(x_{\tilde{\ell}})$ along with $G^{\ell} \cap \partial B_{\sigma_{\ell}}(x_{\ell}) \neq \emptyset$ and $g^{\ell} \cap \partial B_{\sigma_{\ell}}(x_{\ell}) \neq \emptyset.$ The images $G^{\ell}((0,\HH^{1}(G^{\ell})))$ and $g^{\ell}((0,\HH^{1}(g^{\ell})))$ are disjoint.

\bigskip

With this parameterization $G^{\ell},V$ are positively oriented, while $g^{\ell},V$ are negatively oriented.
\end{enumerate}
Then $\Gamma_{1},V$ and $\Gamma_{2},V$ are either both positively oriented or both negatively oriented.

\bigskip

\item[(b)] Consider now a collection of Jordan arcs
$$\Gamma_{1},\Gamma_{2}, \{ G^{\ell},g^{\ell} \}_{\ell=1}^{N}$$
with images that are pairwise either equal or meet only at mutual endpoints, and let $\mathcal{E}$ be the collection of endpoints of these arcs. Let
$$\mathcal{G} = \Gamma_{1} \cup \Gamma_{2} \cup \bigcup_{\ell=1}^{N} G^{\ell} \cup g^{\ell}$$
and suppose $V \in C^{\infty}( \mathcal{G} \setminus \mathcal{E}; S^{1})$ is such that the following hold:

\begin{enumerate}
 \item[(1)] $\Gamma_{1}$ and $\Gamma_{2}$ satisfy (a)(1) above. We also assume $\Gamma_{1},V$ are positively oriented and $\Gamma_{2},V$ are negatively oriented.
 \item[(2)] For each $\ell \in \{ 1,\ldots,N \}$
  $$\mathcal{G} \cap (\clos B_{\sigma_{\ell}}(x_{\ell})) = \left( G^{\ell} \cup g^{\ell} \right) \cap (\clos B_{\sigma_{\ell}}(x_{\ell})).$$
 We parameterize 
 $$G^{\ell} \in C([0,\HH^{1}(G^{\ell})];\clos B_{1}(0)) \cap C^{\infty}((0,\HH^{1}(G^{\ell}));B_{1}(0))$$
 and 
 $$g^{\ell} \in C([0,\HH^{1}(g^{\ell})];\clos B_{1}(0)) \cap C^{\infty}((0,\HH^{1}(g^{\ell}));B_{1}(0))$$
 by arc-length so that $G^{\ell}(\HH^{1}(G^{\ell})) = g^{\ell}(\HH^{1}(g^{\ell})) = x_{\ell}.$ Meanwhile, $G^{\ell}(0) ,g^{\ell}(0) \in (\partial B_{1}(0)) \cup \{ x_{1},\ldots,x_{N} \}.$ The images $G^{\ell}((0,\HH^{1}(G^{\ell})))$ and $g^{\ell}((0,\HH^{1}(g^{\ell})))$ are disjoint.

\bigskip

With this parameterization $G^{\ell},V$ are positively oriented, while $g^{\ell},V$ are negatively oriented.
\end{enumerate}

Then either
$$\mathcal{G} = \Gamma \text{ or } \mathcal{G} = \Gamma \cup \bigcup_{\ell=1}^{N^{(2)}} L_{\ell}$$
where the following hold:
\begin{itemize}
 \item $\Gamma$ is a continuous Jordan arc with endpoints $\Gamma_{1}(0),\Gamma_{2}(0),$ smooth away from the collection of points $\{ x_{1},\ldots, x_{N} \},$ and so that the images of $\Gamma_{1},\Gamma_{2}$ are contained in the image of $\Gamma.$
 
 \item In the latter case, $N^{(2)} \in \N$ and for each $\ell \in \{ 1,\ldots,N^{(2)} \}$ we have that
 $L_{\ell}$ is a continuous closed Jordan curve, smooth away from $\{ x_{1},\ldots,x_{N} \},$ and with $L_{\ell} \cap \{ x_{1},\ldots,x_{N} \} \neq \emptyset.$ The curves $\Gamma,L_{1},\ldots,L_{N^{(2)}}$ are pairwise disjoint.
\end{itemize}
\end{enumerate}
\end{lemma}

{\bf Proof:} To prove (a), suppose for contradiction (and without loss of generality) that $\Gamma_{1},V$ are positively oriented, while $\Gamma_{2},V$ are negatively oriented.

\bigskip

Observe more specifically that $\Gamma_{2},V$ are negatively oriented and $\Gamma_{2}(0) \in \partial B_{1}(0).$ On the other hand, the following hold: 
$$\begin{aligned}
& \Gamma_{1}([0,\HH^{1}(\Gamma_{1}))) \neq \Gamma_{2}([0,\HH^{1}(\Gamma_{2}))); \\  & \gamma^{\ell}_{k},V \text{ are positively oriented}; \\
& \gamma^{\ell}_{k}(\HH^{1}(\gamma^{\ell}_{k})) \in B_{\sigma_{\ell}/2}(x_{\ell}) \text{ for each } \ell \in \{ 1,\ldots,N^{(1)} \},k \in \{ 1,\ldots,N_{\ell} \}; \\
& G^{\ell},V \text{ are positively oriented} \\
& G^{\ell}(\HH^{1}(G^{\ell})) \in B_{\sigma_{\ell}/2}(x_{\ell}) \text{ for each } \ell \in \{ N^{(1)}+1,\ldots,N \}.
\end{aligned}$$
Since $\Gamma_{2}(\HH^{1}(\Gamma_{2})) \in \bigcup_{\ell=1}^{N} B_{\sigma_{\ell}/2}(x_{\ell}),$ then we must have (after relabeling) $\Gamma_{2}=g_{N^{(1)}+1}$; in particular we must have $N^{(1)} < N.$

\bigskip

Now consider $G_{N^{(1)}+1}.$ Since $G_{N^{(1)}+1},V$ are positively oriented and $G^{\ell}(\HH^{1}(G^{N^{(1)}+1})) \in B_{\sigma_{N^{(1)}+1}/2}(x_{N^{(1)}+1}),$ then we cannot have $G^{N^{(1)}+1} = \gamma^{\ell}_{k}$ for any $\ell \in \{ 1,\ldots,N^{1} \},k \in \{ 1,\ldots,N_{\ell} \}.$ However, observe that:
$$\begin{aligned}
& G^{N^{(1)}+1}(0) \in (\partial B_{1}(0)) \cup \bigcup_{\ell=1}^{N} B_{\sigma_{\ell}/2}(x_{\ell}); \\ 
& G^{N^{(1)}+1} \cap \partial B_{\sigma_{N^{(1)}+1}}(x_{N^{(1)}+1}) \neq \emptyset; \\
& G^{N^{(1)}+1} \neq g^{N^{(1)}+1}.
\end{aligned}$$
We conclude that either $G_{N^{(1)}+1} = \Gamma_{1}$ or (after relabeling) $G_{N^{(1)}+1} = g_{N^{(1)}+2}.$ By considering $G_{N^{(1)}+2}$ and arguing iteratively, we see that there is $\tilde{N} \in \{ 1,\ldots,N \}$ so that (after relabeling)
\begin{equation} \label{sards1}
\begin{aligned}
\Gamma_{2} & = g_{N^{(1)}+1}, \\
G_{N^{(1)}+1} & = g_{N^{(1)}+2}, \\
& \vdots \\
G_{\tilde{N}-1} & = g_{\tilde{N}}, \\
G_{\tilde{N}} & = \Gamma_{1}.
\end{aligned}
\end{equation}

\bigskip

Now consider the collection $\{ \{ \gamma^{\ell}_{k} \}_{k=1}^{N_{\ell}} \}_{\ell=1}^{N^{(1)}},$ in particular consider $\gamma^{1}_{1}.$ Since $\gamma^{1}_{1},V$ are positively oriented and $\gamma^{1}_{1}(\HH^{1}(\gamma^{1}_{1})) \in B_{\sigma_{1}/2}(x_{1}),$ then we conclude $\tilde{N} \in \{ 1,\ldots,N-1 \}$ and that (after relabeling) $\gamma^{1}_{1} = g^{\tilde{N}+1}.$ By considering $G^{\tilde{N}+1}$ and arguing iteratively, we see that (after relabeling)
\begin{equation} \label{sards2}
\begin{aligned}
\gamma^{1}_{1} & = g^{\tilde{N}+1}, \\
G^{\tilde{N}+1} & = g^{\tilde{N}+2} \\
& \vdots \\
G^{N-1} & = g^{N}.
\end{aligned}
\end{equation}
However, the fact that $G^{N} \neq \gamma^{\ell}_{k}$ for each $\ell \in \{ 1,\ldots,N^{(1)} \},k \in \{ 1,\ldots,N_{\ell} \}$ together with \eqref{sards1},\eqref{sards2} give a contradiction.

\bigskip

The proof of (b) follows similarly. $\square$

\end{flushleft}
\end{document}